\documentclass[12pt]{article}
\usepackage{a4wide}
\usepackage{amsmath}
\usepackage{theorem}
\usepackage{amsfonts}
\usepackage{epsfig}
\usepackage{md}
\usepackage{hyperref}
\hypersetup{colorlinks=true}
\date{}
\title{Elliptic curves and Hilbert's tenth problem for algebraic function 
fields over real and $p$-adic fields}
\author{Laurent Moret-Bailly
\thanks{The author is a member of the European network `Arithmetic 
Algebraic Geometry' (contract HPRN-CT-2000-00120).
}\medskip
\\
{\small IRMAR (Institut de Recherche Math\'{e}matique de Rennes,} \\
{\small UMR 6625 du CNRS)}\\
{\small Universit\'{e} de Rennes 1}\\
{\small Campus de Beaulieu}\\
{\small F-35042 Rennes Cedex}\\
{\small laurent.moret-bailly@univ-rennes1.fr}\\
{\small{http://name.math.univ-rennes1.fr/laurent.moret-bailly/}}}
\begin{document}
\maketitle
\begin{center}{\bfseries Paper accepted for publication in {\slshape J.\ reine und angew.\ Math.}\\
 (October 2004)}\end{center}
\vfil
\begin{abstract} Let $k$ be a field of characteristic zero, $V$ a 
	smooth, positive-dimensional, quasiprojective variety over $k$, 
	and $Q$ a nonempty 
	divisor on $V$. Let $K$ be the function field of $V$, and $A\subset 
	K$ the semilocal ring of $Q$.
	
	We prove the \dio\ 
	undecidability of: (1) $A$, in all cases; (2) $K$, when $k$ is 
	real and $V$ has a real point; (3) $K$, when $k$ is a subfield of a 
	$p$-adic field, for some odd prime $p$.
	
	To achieve this, we use Denef's method: from an elliptic 
	curve $E$ over $\QQ$, without complex multiplication, one constructs 
	a quadratic twist $\calE$ of $E$ over $\QQ(t)$, which has 
	Mordell-Weil rank one. Most of the paper is devoted to proving (using 
	a theorem of R.~Noot) that one can choose $f$ in $K$, vanishing at $Q$, 
	such that the group $\calE(K)$ deduced from the field 
	extension $\QQ(t)\flis\QQ(f)\inj K$ is equal to $\calE(\QQ(t))$. Then we 
	mimic the arguments of Denef (for the real case) and of Kim and Roush 
	(for the $p$-adic case).
\end{abstract}
\vfill

\begin{center}
\end{center}

\vfil\eject

\setcounter{tocdepth}{2}\tableofcontents

\section{Introduction}\label{intro}

The aim of this paper is to prove the following  result:

\begin{thm}\label{ThVague}
	Let $k$ be a field of characteristic zero. Let $V$ be a smooth,  
	positive-dimensional, quasiprojective, irreducible $k$-scheme, 
	with function field denoted by $K$.
	\smallskip
	
	\noindent{\rm(1) (see Theorem \ref{ThIndecSemiloc})} Let $Q$ be a 
	nonempty effective divisor on $V$, 
	and let $A\subset K$ be the semilocal ring of $Q$ (the intersection 
	of the local rings of the maximal points of $Q$). Then the 
	positive-existential theory of $A$ is undecidable. In other 
	words, Hilbert's tenth problem over $A$ has a negative solution. 
	\smallskip
	
	\noindent{\rm(2) (see Theorem \ref{ThReel})} Assume that $K$ is 
	formally real. Then the positive-existential theory of $K$ is undecidable.
	\smallskip
	
	\noindent{\rm(3) (see Theorem \ref{Thpadic})} Assume that $k$ is 
	a subfield of a finite 
	extension of $\QQ_{p}$, for some odd prime $p$. Then the 
	positive-existential theory of $K$ is undecidable.
	\smallskip
\end{thm}

\begin{subrem}\label{RemThVague1}
Thus, for instance, the conclusion of (2) and (3) means that there is no 
algorithm taking as input a polynomial $F\in K[X_{1},\ldots,X_{n}]$ 
(for some $n$) and giving a `yes/no' output 
according as $F$ has a zero in $K^n$ or not.

In all statements, the positive-existential theory is considered 
in the language of rings, augmented by a suitable set of constants 
which can be described. 

In each case, a more precise result will be that there is a \dio\ 
(that is, positive-existentially definable)
subset of $A^d$ (resp.~$K^d$) for some $d$ (in fact $d=2$ in cases (1) 
and (2)), with a ring structure which is also 
\dio\ and isomorphic to $\ZZ$ as a ring. By a standard argument, this 
together with the negative solution of Hilbert's tenth problem over 
$\ZZ$ (Davis-Putnam-Robinson-Matijasevich) implies the result for 
$A$ (resp.~$K$).
\end{subrem}
\begin{subrem}\label{RemThVague2}
Note that by enlarging $k$, we can assume in (1) and (2)  that $V$ 
is a \emph{curve} $C$: if $r=\dim V$, one can first assume $V$ affine 
(taking an open  
subset meeting every component of $Q$), then choose a $k$-morphism 
$V\to\Aa^{r-1}_{k}$ whose generic fibre is a smooth curve and such 
that every component of $Q$ dominates $\Aa^{r-1}_{k}$; finally, 
replace $k$ by the function field of $\Aa^{r-1}_{k}$. One can even go 
further and assume that $C$ is projective and smooth (by completing it) and 
geometrically connected over $k$ (by replacing $k$ by its algebraic 
closure in $K$). 

This reduction to the case of curves does not work in case (3): for 
instance, $\QQ_{p}(x)$ cannot be embedded in a $p$-adic field. 
However, in most of this paper, the emphasis will be on curves.
\end{subrem}
\begin{subrem}\label{RemThVague3}
Several special cases of \ref{ThVague} were known before. The case 
where $k$ is a real field and $K=k(t)$ is due to Denef \cite{Denef}, 
as well as the method used here. 

Denef's method was also used by Kim and Roush \cite{KR} to treat the 
case where $k$ is as in (3) and $K=k(t)$. We use the proof of Kim and 
Roush to prove (3), whence the restriction $p\neq2$.

Some special cases of (2) for non-rational function fields in one 
variable over real fields were obtained by Zahidi in \cite{Za}. 

Earlier versions of this paper, without part (3), were
circulated before. After completing (3), the author was informed (at the end 
of August 2004) that K.~Eisentr\"ager \cite{Eis} had independently proved the 
$p$-adic case (3), using one of these versions (specifically, Theorem 
\ref{ThIntroTer} below). She also obtained in \cite{Eis1} (with the notations of 
Theorem \ref{ThVague}) the \dio\ undecidability of $K$ when $k$ is 
algebraically closed and $\dim V\geq2$, adapting (again via our 
Theorem \ref{ThIntroTer}) the method used in \cite{KR1} for $K=\CC(t_{1},t_{2})$.

The reader may consult \cite{PhZa} for a review of other related results.
\end{subrem}
Note that (3) implies in particular:
\begin{subcor}\label{CorExtQ} Let $K$ be a finitely generated, 
transcendental extension of $\QQ$. Then $K$ is positive-existentially 
undecidable.
\end{subcor}
From now on, we shall assume in this introduction that $V=C$ is a smooth, 
projective, geometrically connected curve over $k$, and that $Q$ is a finite 
nonempty set of closed points of $C$.

To motivate our further setting, we now briefly recall Denef's method. 

\Subsection{Sketch of Denef's method. }\label{SsecMethDenef} Assume 
$C=\Pu{k}$, with standard coordinate $t$, so that $K=k(t)$. 
Take an elliptic curve $E$ over 
$k$ (defined over $\QQ$, if we wish), and `twist' it by the quadratic extension of 
$K=k(t)$ given by the usual double cover $\pi:E\to\PP^1_{k}$. The result 
is an elliptic curve $\calE$ over $K$ (the `self-twist' of $E$), with additive reduction 
at the branch points of $\pi$. An easy computation shows that $\calE(K)$ is 
`almost' the endomorphism ring of $E$. More precisely, $2\calE(K)$ is 
canonically isomorphic to $2\,{\rm End}_{k}\,(E)$, so if $E$ does not 
have complex multiplication the group $2\calE(K)$ is isomorphic to 
$\ZZ$, and the addition is clearly \dio\ since it is induced by the 
group law of $\calE$. 

Fixing an isomorphism $2\calE(K)\flis\ZZ$, the less 
obvious fact that the multiplication is also \dio\ is deduced from the `additive 
reduction' properties of $\calE$ at branch points of $\pi$, in 
particular at the point $\infty$. Specifically, we have a `reduction', 
or `specialisation' homomorphism from $2\calE(K)$ to the additive group 
$k$, which turns out to be nonzero, hence must be (up to a harmless 
constant) the inclusion of $\ZZ$ into $k$. In particular, it must be 
compatible with multiplication, allowing us to obtain a \dio\ 
definition of multiplication in $2\calE(K)$ --- \emph{provided}, 
however, that the specialisation map has good \dio\ properties, which 
is where the `real' (resp.~`$p$-adic') assumption is used; more 
precisely, this involves proving the \dio\ definability of certain 
subsets of $K$ defined by valuation conditions.

\Subsection{Extending  Denef's method to other fields. }\label{SsecExtDenef}
We want to extend this argument with $k(t)$ replaced by $K$. To do 
this, we simply take a cover $f:C\to\PP^1_{k}$ with reasonable properties, 
by means of which we identify $k(t)$ with a subfield of $K$; we then try to 
adjust the data in such a way that $\calE(K)=\calE(k(t))$ (whatever $\calE(k(t))$ may be: we 
shall forget here about the `no complex multiplication' condition).

As it turns out, we may in fact start with any elliptic curve $E$ over $k$, 
and twist it by a quadratic extension of $k(t)$, corresponding to a double 
cover $\pi:\Gamma\to\PP^1_{k}$ (here $\pi$ is any double cover of 
$\Pu{k}$ by a smooth curve $\Gamma$, not necessarily $E$ itself). The curves $C$, $E$, 
and $\Gamma$, and the morphism $\pi$, will be fixed throughout (and, therefore, so 
will the twisted elliptic curve $\calE$ over $k(t)$). The only `variable' 
piece of data is the morphism $f$; specifically, we shall allow 
ourselves to replace the initially given $f$ by $\lambda f$, for some suitable 
$\lambda\in k^\ast$. 

But now it is time to fix the notations more 
precisely. (For the rest of this introduction, we shall concentrate on 
the algebro-geometric result of the paper; the applications to 
Hilbert's tenth problem will be considered in Part \ref{Part3}).

\Subsection{Notations. }\label{DonnFond}

\Subsubsection{The ground field. }\label{NotCorps} 
In the rest of this introduction (and in most of the paper), $k$ denotes a field of
characteristic $p\geq0$. Moreover, \emph{unless otherwise specified, we shall always 
assume that $p\neq2$.} (For applications to undecidability questions, 
$p$ will be zero). We fix an algebraic closure of $k$, denoted by $\kb$. 

\Subsubsection{The fundamental curve $C$ and its function field. }\label{NotC} 
We denote by $C$ a smooth projective geometrically connected curve 
over $k$, with function field $K$ (thus, $K$ is a finitely 
generated extension of $k$, of transcendence degree $1$, and $k$ is 
algebraically closed in $K$). 

If $k'$ is an extension of $k$, the function field of 
$C_{k'}=C\times_{\Spec(k)}\Spec(k')$ will be denoted by $k'(C)$; 
thus, $\kb(C)=\kb\otimes_{k}K$, and in general $k'(C)$ is the fraction 
field of $k'\otimes_{k}K$.

We are also given a finite nonempty set $Q$ of closed points of $C$; 
we assume that their residue fields are \emph{separable\/} over 
$k$ (in other words, $Q$ is the spectrum of an \'{e}tale $k$-algebra).

\Subsubsection{The elliptic curve. }\label{NotE0} 
$E$ denotes an elliptic curve over $k$. 
It will be convenient to fix an affine equation of $E$, of the form
\begin{equation}\label{EqE0}
	y^2=P(x)
\end{equation}
for a cubic polynomial $P\in k[T]$ without multiple roots. 

\Subsubsection{The hyperelliptic curve. }\label{NotGamma} 
$\Gamma$ denotes a smooth projective geometrically connected curve over $k$, 
given as a double cover 
\begin{equation}\label{EqPi}
	\pi:\Gamma\ffl\PP^1_{k}.
\end{equation}
We shall always assume that:
\begin{sitemize}
	\item $\pi$ is \'{e}tale above $\infty\in\PP^1(k)$;
	\item $\pi$ is ramified at $0$.
\end{sitemize}
(The second assumption is made only to fix ideas and avoid some case 
discussions, and because it is the relevant case for applications to 
Hilbert's tenth problem. The first assumption, however, will be 
essential in our constructions).

Thus, $\Gamma$ can be described by an affine equation, in 
coordinates $(t,w)$:
\begin{equation}\label{EqGamma}
	w^2=R(t)
\end{equation}
where $t$, the standard coordinate on $\PP^1$, is identified with the 
rational function $t\circ\pi=\pi$ on $\Gamma$, and $R$ is a polynomial 
in $k[T]$, without multiple roots, such that $R(0)=0$, and  
$\deg R=2\,\text{genus}(\Gamma)+2$. 

Clearly, $\Gamma$ has a unique $k$-rational point 
above the point $0\in\PP^1(k)$. We denote this point by 
$0_{\Gamma}$.

The natural involution of $\Gamma$, sending $w$ to $-w$, will be 
denoted by $\sigma$.

\Subsubsection{Remarks. }\label{RemGamma}
\begin{romlist}
	\item\label{RemGamma1} We shall always think of $\Gamma$ 
	as equipped with the double 
	cover $\pi$ of (\ref{EqPi}). In other words, when using the 
	notation $\Gamma$ we shall often actually mean $\pi$. Note that 
	$\Gamma$ and $\pi$ are completely determined by the polynomial $R$ 
	of (\ref{EqGamma}); conversely, they determine $R$ up to a square 
	factor in $k^\ast$.
	\item\label{RemGamma2} An important special case is when $\Gamma=E$ 
	(thus not `hyperelliptic', strictly speaking!) and $\pi$ is the double cover given 
	by $x^{-1}$, the inverse of the coordinate $x$ in (\ref{EqE0}) . In 
	this case, the polynomial $R$ of (\ref{EqGamma}) is
	\begin{equation}\label{EqGamma=E1} R(t)=t^4\,P(1/t)
	\end{equation}
	and the functions $t,w,x,y$ on $E$ are related by 
	\begin{equation}\label{EqGamma=E2} (t,w)=(1/x,y/x^2)\quad 
	(x,y)=(1/t,w/t^2).
	\end{equation}
	In fact, this is the important case for applications to Hilbert's 
	tenth problem; however, the author feels 
	that restricting to this special case would only give a less general 
	result without any substantial simplification, while distinguishing 
	between $E$ and $\Gamma$ actually clarifies the proof.
	\item\label{RemGamma3} We could even have generalised further, by replacing $E$ by any 
	abelian variety over $k$. But this time, this would make the results 
	we need (essentially those of Section \ref{SecMD} on quadratic twists) 
	slightly less elementary.
\end{romlist}

From the data \ref{NotE0} and \ref{NotGamma} we can now perform a 
well-known construction:

\Subsubsection{The twisted elliptic curve $\calE$ over $k(t)$. 
}\label{DefMD}
We may define it as the $k(t)$-elliptic curve with affine equation
\begin{equation}\label{EqMD}
		y^2=R(t)\,P(x)
\end{equation}
(in the affine plane $\Aa^2_{k(t)}$, with coordinates $x$, $y$). We 
shall give a more intrinsic definition of twists in Section 
\ref{SecRevDouble}, and use slightly different (but of course 
equivalent) equations.

It is easy to compute the Mordell-Weil group $\calE(k(t))$ of $\calE$ in 
terms of morphisms of curves {\it over $k$}: as 
we shall see in \ref{SsecPtsMDCourbes}, there is a canonical isomorphism
\begin{equation}\label{MWMD}
	\calE(k(t))\cong\modd_{k}(\Gamma,E)
\end{equation}
where the right-hand side stands for the group of $k$-morphisms 
$h:\Gamma\to E$ compatible with involutions, i.e.~such that 
$h\circ\sigma=[-1]_{E}\circ h$.

\begin{subrem} 
	In the case $\Gamma=E$ of \ref{RemGamma}\,\ref{RemGamma2}, we get from 
	(\ref{MWMD}) an isomorphism
	\begin{equation}
		\label{EqDenef}
		E[2](k)\times{\rm End}_{k}(E)\cong\calE(k(t))
	\end{equation}
	where $E[2]$ is the kernel of multiplication by $2$ in 
	$E$.  Concretely, using the co\-or\-di\-nates in (\ref{EqE0}) and 
	(\ref{EqMD}), this sends $((\xi,0),0)$ to $(\xi,0)$ (where $\xi$ is a 
	zero of $P$), and $(0_{E},{\rm Id}_{E})$ to $(1/t,t^2 P(1/t))$ (recall 
	that $R(t)=t^4 P(1/t)$).
\end{subrem}

\Subsection{Properties of covers $C\to\PP^1_{k}$: good functions. 
}\label{PropCov}\par

Recall that our goal is to extend Denef's argument with $k(t)$ replaced 
by $K$. To do this we shall choose a suitable nonconstant rational 
function $g$ on $C$; this defines a ramified cover $g:C\to\PP^1_{k}$, 
and a corresponding field extension $k(t)\to K$, sending $t$ to $g$. 

We denote by $K_{g}$ the field $K$ viewed as an extension of $k(t)$ 
via $g$. Thus, we have an obvious inclusion of abelian groups
\begin{equation}\label{InclPtsMD}
	\calE(k(t))\inj\calE(K_{g})
\end{equation} 
both of which will turn out (see \ref{RemStrucPtsMDCourbes}) to be finitely 
generated with the same torsion subgroup, isomorphic to the kernel $E[2](k)$ of 
multiplication by $2$ in $E(k)$.

We would like (\ref{InclPtsMD}) to be an equality, for suitable $g$. If $p>0$, 
however, we can only achieve this `up to $p$-torsion', which 
motivates the following definition:
\begin{subdefi}\label{DefAlmostOnto} Let $u:A\to B$ be a morphism of 
abelian groups. We shall say that $u$ is \emph{almost bijective} 
if $u$ is injective and $\mathrm{Coker}\,u$ is a finite $p$-group.
\end{subdefi}
Of course, if $p=0$ we take this to mean that $u$ is bijective. In 
the sequel we shall only apply this notion to morphisms of 
finitely generated abelian groups.
\begin{subdefi}\label{DefGood} 
	Let  $k$, $C$, $Q$, $E$, $\Gamma$ be as in \rref{DonnFond}, and
	let $g:C\to\PP^1_{k}$ be a nonconstant $k$-morphism. 
	\smallskip
	
	\noindent{\rm(1)} We say that $g$ is \emph{admissible\/} 
	(for $\Gamma$, or for $\pi$) if:
	\begin{romlist}
	\item\label{DefGood1} $g$ has only simple branch points 
	\emph{(i.e.~no ramification index $\geq3$)};
	\item\label{DefGood2} $g$ is \'{e}tale above $\infty$ and the branch points 
	of $\pi$ \emph{(which are the zeros of $R$)};
	\item\label{DefGood3} every point of $Q$ is a zero of $g$ \emph{(automatically 
	simple, by \ref{DefGood2})}.
	\end{romlist}
	\noindent{\rm(2)} We say that $g$ is \emph{good\/} for $E$ and 
	$\Gamma$ if $g$ is admissible and the natural inclusion 
	{\rm(\ref{InclPtsMD})} is almost bijective.
	
	If $k'$ is an extension of $k$, we say that $g$ is 
	\emph{good over $k'$\/}, or  \emph{$k'$-good\/}, if the morphism 
	$g_{k'}:C_{k'}\to\Pu{k'}$ deduced from $g$ by base change is good 
	for $E_{k'}$ and 
	$\Gamma_{k'}$. 
	
	We say that $g$ is \emph{very good\/} if $g$ is $\ol{k}$-good.
	\medskip
	
	\noindent{\rm(3)} Let $f:C\to\PP^1_{k}$ be admissible. For every 
	extension $k'$ of $k$, define two subsets $\good(k')$ and 
	$\vgood(k')$  of ${k'}^\ast$ by
	$$\begin{array}{rcrcl}
		\good(k')&=&\good(E,\Gamma,f,k')&=
		&\{\lambda\in{k'}^\ast\,\mid\,\lambda\,f\hbox{\rm\ is 
	good for $E_{k'}$ and }\Gamma_{k'}\}\\
			\vgood(k')&=&\vgood(E,\Gamma,f,k')&=
			&\{\lambda\in{k'}^\ast\,\mid\,\lambda\,f\hbox{\rm\ is 
	very good for $E_{k'}$ and }\Gamma_{k'}\}.
	\end{array}
	$$
\end{subdefi} 
\Subsubsection{Remarks.}\label{RemDefGood}
\begin{romlist}
	\item\label{RemDefGood-1} By definition, $g$ is $k'$-good if and only 
	if $g$ is admissible and the natural inclusion 
	\begin{equation}\label{InclPtsMDExt}
		\calE(k'(t))\inj\calE(k'(C)_{g})
	\end{equation} 
	is  almost bijective; here $k'(C)_{g}$ is the function field $k'(C)$ of $C_{k'}$, 
	viewed as an 
	extension of $k'(t)$ via $g_{k'}$. In particular, $g$ is  very good if 
	and only if $g$ is admissible and $\calE(\kb(C)_{g})=\calE(\kb(t))$, up 
	to $p$-torsion.
	\item\label{RemDefGood0} Of course, the definition of $\vgood(k')$ 
	refers implicitly to some algebraic closure $\ol{k'}$ of $k'$, but 
	is, as usual, independent of it.
	\item\label{RemDefGood3} Let $f$ be admissible. Then for all but finitely 
	many $\lambda\in k^\ast$, the function $\lambda f$ is still admissible. 
	Thus, except for finitely many $\lambda$, the `goodness' property for 
	$g=\lambda f$ just means that {\rm(\ref{InclPtsMD})} is almost bijective.
\end{romlist}

\begin{subprop}\label{PropGood} Assume that $g:C\to\Pu{k}$ is 
admissible, and let $k'$ be an extension of $k$. Then: 
	\begin{romlist}
		\item\label{PropGood1} If $g$ is $k'$-good, then it is good. In 
		particular, very good morphisms are good.
		\item\label{PropGood2} Assume that $k$ is separably closed in $k'$.
		Then $\calE(k'(t))=\calE(k(t))$ and $\calE(k'(C)_{g})=\calE(K_{g})$; in 
		particular, $g$ is $k$-good if and only if it is $k'$-good. 
		\item\label{PropGood3} $g$ is very good if and only if $g$ is $F$-good for every 
		extension $F$ of $k$.
	\end{romlist}
\end{subprop}
\dem \ref{PropGood1} follows easily from the fact that $k(t)=k(C)\cap 
k'(t)$. 

The proof of \ref{PropGood2} will be postponed until \ref{DemPropGood2}.

The `if' part of \ref{PropGood3} is trivial. Conversely, assume $g$ is 
very good, and let $F$ be an extension of $k$, with an algebraic 
closure $\ol{F}$ containing $\kb$. Now $g_{\kb}$ is good, hence $g_{\ol{F}}$ 
is good by \ref{PropGood2}. Hence $g_{F}$ is good by \ref{PropGood1}.\qed

\begin{subcor}\label{CorPropGood} Let $f:C\to\Pu{k}$ be
admissible, and let $k'$ be an extension of $k$. Then: 
	\begin{romlist}
		\item\label{CorPropGood1} $\good(k')\cap k\subset \good(k)$, with 
		equality if $k$ is separably closed in $k'$.
		\item\label{CorPropGood2} $\vgood(k')\cap k=\vgood(k)$.\qed
	\end{romlist}
\end{subcor}
\begin{subrem}\label{RemDefGood1} The existence of admissible morphisms is 
	easy (see \ref{ExistAdm}, but note that in positive characteristic, 
	this uses our assumption that $Q$ is \'etale over $k$). The existence of very 
	good morphisms is the subject of this paper. More precisely, if $k$
	is not algebraic over a finite field, we shall prove $\vgood(k)\neq\emptyset$  
	for any admissible $f$. This of course implies $\good(k)\neq\emptyset$. 
	The reasons why we need both variants are explained in \ref{SsecEffect} below.
\end{subrem}
\Subsection{One last piece of data. }\label{Notf} In addition to the 
data of \ref{DonnFond}, we fix an admissible $k$-morphism
\begin{equation}\label{Eqf}
	f:C\to\PP^1_{k}.
\end{equation}
\begin{mainthm}\label{MainTh} 
	Let $k$, $C$, $Q$, $E$, $\Gamma$, $f$ be as above. Then:
	\begin{romlist}
		\item\label{MainTh2} Let $k'$ be an extension of $k$. 
		If $\lambda\in k'$ is transcendental over $k$, then $\lambda\in\vgood(k')$.
		\item\label{MainTh3} If $k$ is finitely 
		generated over the prime field, then $\vgood(k)$ contains a Hilbert 
		subset of $k$, in the sense of {\rm\cite{FJ}, 11.1}; in other words, 
		its complement in $k$ is a thin set  
		in the sense of {\rm\cite{Serre}}.
	\end{romlist}
\end{mainthm}
This will be proved in \ref{FinDemMainTh}. Let us now explore some 
consequences.
\begin{thm}\label{ThIntroTer}
	We keep the notations and assumptions of Theorem \rref{MainTh}. 
	\begin{romlist}
		\item\label{ThIntroTer1} Let $k_{0}$ be any subfield of $k$, 
		finitely generated over the prime field. Then $\vgood(k)$ contains a Hilbert 
		subset of $k_{0}$.
		\item\label{ThIntroTer2} If $\car k=0$, then $\vgood(k)\cap\ZZ$ is infinite.
		\item\label{ThIntroTer2,5} If $\car k=p>0$, and $u\in k$ is 
		transcendental over $\FF_{p}$, then $\vgood(k)\cap \FF_{p}[u]$ is infinite.
		\item\label{ThIntroTer3} The complement of $\vgood(k)$ in $k$ has 
		finite transcendence degree over the prime field. In particular, this 
		complement is countable.
	\end{romlist}
\end{thm}
\dem 
\ref{ThIntroTer1} If $k_{0}$ is finite, the claim is empty. 
Therefore we may assume that $k_{0}$ is either finitely generated 
over $\QQ$, or finitely generated and transcendental over a finite 
field. In both cases, $k_{0}$ is Hilbertian.

There is a subfield $k_{1}$ of $k$, containing 
$k_{0}$, finitely generated over the prime field, and such that $C$, $E$, 
$\Gamma$, and $f$ are defined over $k_{1}$. Now apply \ref{MainTh} 
with $k_{1}$ as ground field: by 
\ref{CorPropGood}\,\ref{CorPropGood2}, $\vgood(k)$ contains 
$\vgood(k_{1})$ which contains a Hilbert subset of $k_{1}$ by 
\ref{MainTh}\,\ref{MainTh3}. Since $k_{0}$ is 
Hilbertian, it follows that $\vgood(k_{1})\cap k_{0}$ contains a 
Hilbert subset of $k_{0}$, by  \cite{FJ}, 11.7 and 11.8(b). 

Assertion \ref{ThIntroTer2} then follows from \ref{ThIntroTer1} (with 
$k_{0}=\QQ$) and \cite{FJ}, Theorem 12.7, and 
\ref{ThIntroTer2,5} is similar. 

For \ref{ThIntroTer3}, take $k_{1}$ as in the proof of \ref{ThIntroTer1}: then 
\ref{MainTh}\,\ref{MainTh2}, applied over $k_{1}$, shows that 
the complement of $\good(k)$  is contained in the algebraic closure 
of $k_{1}$ in $k$, whence the result.\qed

\begin{rem} Theorem \ref{MainTh} may well be 
	true in characteristic $2$. Presumably, the arguments of this paper, 
	suitably adapted, might lead to a proof that $\good(k)$ contains a 
	Hilbert set when $k$ is finitely 
	generated over $\FF_{2}$. However, to obtain the same result for 
	$\vgood(k)$ (and hence the general result, for arbitrary $k$), the 
	present proof makes use of very strong properties of a pencil of 
	curves over $\PP^1_{k}$ considered in Section \ref{SecPencil} 
	(namely, that its fibres over $\PP^1\setminus\{0\}$ are semistable and its fibre at 
	$0$ is `tame'). Both these properties fail in general in characteristic 
	$2$, which definitely ruins the crucial Lemma \ref{LemNonRam}.
	
	Of course, to treat the characteristic $2$ case one would first have 
	to rewrite the generalities of Section \ref{SecRevDouble} on double covers and 
	twists. 
	
	In any case, the applications to Hilbert's tenth problem, which were 
	the prime motivation for this paper, work only in characteristic 
	zero.
\end{rem}

\Subsection{Outline of the proof of the Main Theorem. 
}\label{SsecOutline}
For simplicity, we assume $p=0$ (thus, `almost bijective' just means 
`bijective'). For $\lambda\in k$, we consider the inclusion 
$\calE(k(t))\inj\calE(K_{\lambda f})$ (resp.~$\calE(\kb(t))\inj\calE(\kb(C)_{\lambda 
f})$). The first group is independent of $\lambda$, while the second 
varies with $\lambda$, and clearly we have to make the groups $\calE(K_{\lambda 
f})$ and $\calE(\kb(C)_{\lambda f})$) `as small as possible'. 

\Subsubsection{}
Our first task is to `compute' all these groups in terms of `geometry 
over $k$'. For $\calE(k(t))$, and similarly for $\calE(\kb(t))$, this is 
achieved by formula (\ref{MWMD}). To generalise this, we introduce (in 
\ref{DefCg}) the $k$-curve
\begin{equation}\label{EqCTil}
	\til{C}_{\lambda f}:=C\times_{\lambda f,\Pu{k},\pi}\Gamma
\end{equation}
consisting of pairs $(c,\gamma)$ in $C\times\Gamma$ such that $\lambda 
f(c)=\pi(\gamma)$. As a double cover of $C$, it carries a natural 
involution, and just as in (\ref{MWMD}) there is a canonical 
isomorphism 
\begin{equation}\label{PtsTwEllKIntro}
	\calE(K_{\lambda f})\cong
	\modd_{k}(\til{C}_{\lambda f},E).
\end{equation}
It follows (Proposition \ref{Trad1}) that $\lambda\in\good(k)$ if 
and only if every odd $k$-morphism $\til{C}_{\lambda f}\to E$ is obtained 
from an odd $k$-morphism $\Gamma\to E$ by composition with the natural map 
$\til{C}_{\lambda f}\inj C\times\Gamma\to\Gamma$ (and of course there 
is a similar characterisation of $\vgood(k)$ in terms of 
$\kb$-morphisms). 

The next step consists in translating this condition in terms of 
Jacobians, and removing the `odd' restriction. The result is this 
(Propositions \ref{RedPicSurf} and \ref{RedPicSurfBar}): consider the 
morphism of abelian varieties
\begin{equation}\label{EqMorJac0}
	\jac(C)\times\jac(\Gamma)\ffl\jac(\til{C}_{\lambda f})
\end{equation}
deduced from the natural projections from $\til{C}_{\lambda f}$ to $C$ 
and $\Gamma$. Then $\lambda$ is in $\good(k)$ (resp.~in $\vgood(k)$) if every 
$k$-morphism (resp.~$\kb$-morphism) of abelian varieties
$E\to\jac(\til{C}_{\lambda f})$ factors through (\ref{EqMorJac0}). 

\Subsubsection{} So, putting $J_{\lambda}:=\jac(\til{C}_{\lambda f})$, we are 
led to investigate how the groups $H_{\lambda}:=\Hom_{k}\,(E,J_{\lambda})$ 
and $\ol{H}_{\lambda}:=\Hom_{\kb}\,(E,J_{\lambda})$ vary with $\lambda$. 

Now $J_{\lambda}$ is the fibre at $\lambda$ of a pencil of abelian 
varieties parametrised by an open subset $U\subset\Pu{k}$. Denote by 
$\eta=\Spec(k(z))$ the generic point of $U$ (here $z$ denotes the 
natural coordinate on $\Pu{k}$): the groups in question have `generic' 
values $H_{\eta}:=\Hom_{k(z)}\,(E,J_{\eta})$ 
and $\ol{H}_{\eta}:=\Hom_{\ol{k(z)}}\,(E,J_{\lambda})$. 
There are injective `specialisation' maps $H_{\eta}\inj H_{\lambda}$ 
and $\ol{H}_{\eta}\inj\ol{H}_{\lambda}$, and 
a `specialisation theorem' due to R.~Noot asserts that when $k$ is 
finitely generated over $\QQ$, these specialisation maps are 
isomorphisms for every $\lambda$ in a Hilbert subset of $k$. 

Hence, to prove the Main Theorem, it suffices to show that our 
generic groups $H_{\eta}$ and $\ol{H}_{\eta}$ are  isomorphic to 
$\Hom_{k}\,(E,\jac(C)\times\jac(\Gamma))$ and 
$\Hom_{\kb}\,(E,\jac(C)\times\jac(\Gamma))$, respectively. 

To achieve this, we have to go back to the definition of the curves 
$\til{C}_{\lambda f}$: as curves on the surface $C\times\Gamma$, they 
are the fibres of the rational map 
$C\times\Gamma\cdots\!\!\fl\Pu{k}$ sending $(c,\gamma)$ to 
$\pi(\gamma)/f(c)$. The result for $H_{\eta}$ then follows from more or 
less standard facts (Theorem \ref{ThSpecJac}) on Jacobians of pencils 
of curves. The analogous result for $\ol{H}_{\eta}$ (Theorem \ref{deltaIsom}) is more 
delicate and requires a detailed analysis of the degenerations of the 
pencil, carried out in Section \ref{SecPencil}.

\Subsection{Effectivity questions. }\label{SsecEffect}
	Observe that in the proof of \ref{ThIntroTer}\,\ref{ThIntroTer1}, we 
	have used the inclusion $\vgood(k_{1})\subset\vgood(k)$ in an essential 
	way. This explains why we need to consider `$\vgood$' sets, even to obtain the 
	result for $\good(k)$ for arbitrary $k$. 
	
	Unfortunately, the proof that $\vgood(k)\neq\emptyset$ relies on a highly 
	nonconstructive argument involving infinite Galois groups (this occurs 
	in the proof of the specialisation theorem \ref{ThSpec}). 

On the other hand, if we limit ourselves to fields $k$ finitely generated over 
the prime field (and to proving that $\good(k)\neq\emptyset$), then there is a more 
effective result, stated below and proved in  \ref{FinDemMainTh} (here $\jac(X)$ 
denotes the Jacobian of a curve $X$). We refer to \cite{FJ}, Chapter 
17 for presented fields and related notions; in particular, recall 
that $k$ is \emph{presented} over its prime field $\kappa$ if it is 
described as $k=\kappa(x_{1},\ldots,x_{n})$ where, for each $i\geq1$, the 
minimal polynomial of $x_{i}$ over $k_{i-1}:=\kappa(x_{1},\ldots,x_{i-1})$ is 
explicitly given (and understood to be zero if $x_{i}$ is 
transcendental over $k_{i-1}$). Many standard algebro-geometric 
constructions over $k$ can then be carried out `effectively'; see \cite{FJ} 
for details. 
\begin{thm}\label{MainThEff} {\rm(Effective version of the Main Theorem)} 
	Assume that $k$ is presented over the prime field, and that the rank of the 
	finitely generated abelian groups $\Hom_{k}\,(E,\jac(\Gamma))$ and 
	$\Hom_{k}\,(E,\jac(C))$ are known.
	
	Then $\good(k)$ contains an `effective' Hilbert subset of $k$, 
	in the following sense: 
	$z$ and $y$ denoting indeterminates, there is an effectively 
	computable $\Phi(z,y)\in k(z)[y]$, with no root in $k(z)$ (as a 
	polynomial in $y$), with the property that for all $\lambda\in 
	k$, if $\Phi(\lambda,y)\in k[y]$ has no root in $k$ 
	then $\lambda\in\good(k)$.
\end{thm}

\Subsection{Organisation of the paper}\label{SsecOrga}
Apart from this introduction, the paper is divided into three parts.
\medskip

\noindent Part I exposes background material, mostly from algebraic geometry; 
nothing in this part is really new. 

Section \ref{SecNot} contains miscellaneous 
(and more or less well-known) results and basic definitions. 

Section \ref{SecThSpec} is devoted to Noot's specialisation theorem; we give a 
proof there because 
the statement we use is actually a variant of Noot's original result. 
We also give a proof of the `effective' variant we need (see 
\ref{SsecEffect} above). 

In section \ref{SecJacFib}, we give some important (and perhaps 
not so familiar) properties of the relative Jacobian of a  
surface fibered over the projective line. These properties are at the 
heart of our proof of the Main Theorem. 

Finally, Section \ref{SecRevDouble} presents 
the basic facts on quadratic twists, especially of elliptic curves. 
\medskip

\noindent In Part II, we prove the Main Theorem, as outlined in \ref{SsecOutline} 
above.
\medskip

\noindent Part III contains the applications to model theory, and the proof of 
Theorem \ref{ThVague}. In this part, the Main Theorem is applied in 
the special case where $E=\Gamma$ (a fixed elliptic curve over $k$), 
as explained in \ref{RemGamma}\,\ref{RemGamma2}. The resulting twist $\calE$ is the 
`self-twist' of  $E$; generalities on this construction are exposed in 
Section \ref{SecSelfCan}.

In Section \ref{SecLambda}, we define the (hopefully \dio) model of $\ZZ$ deduced 
from $\calE$; this ring is denoted by $\Lambda$. This is a subset of 
$K^2$. To show the \dio\ undecidability of $K$, we have to prove two 
things: that $\Lambda$ is a \dio\ set (this is where we use the Main
Theorem), and that the multiplication of $\Lambda$ is relatively 
\dio\ (a notion explained in \ref{RelDio}); how we prove the latter 
depends on the field (or ring) $K$. 

Section \ref{SecIndecSemiloc} contains the proof of part (1) of 
Theorem \ref{ThVague}, as well as general notations used in the sequel.

Section \ref{SecReel} contains the proof of part (2) of \ref{ThVague}, following 
Denef.

Finally, in Section \ref{Secpadic} we prove part (3) of \ref{ThVague}, adapting 
the method of Kim and Roush.

\Subsection{Acknowledgments. }
The author is grateful to Karim Zahidi, Luc B\'{e}lair, Bas Edixhoven for 
discussions on the subject of this paper, and most especially to 
Rutger Noot for discussions on the specialisation theorem.
\part{Geometric background}
\section{Basic material}\label{SecNot} 

Throughout this section, $F$ denotes a field, $\ol{F}$ 
an algebraic closure of $F$, and ${F}^{\rm s}$ the 
separable closure of $F$ in  $\ol{F}$.

\Subsection{Rings, varieties, morphisms. }\label{SsecVar}
All rings are commutative with unit.

If $S$ is a scheme, and $X$, $Y$ are $S$-schemes, we shall 
denote by $\mor_{S}\,(X,Y)$ the set of $S$-scheme morphisms from $X$ to 
$Y$. If $X=\Spec(R)$ is affine, we also use the notation 
$Y(R)$.

In general we use subscripts to denote base change: thus, if $S'$ is 
an $S$-scheme, we write $X_{S'}$ for $X\times_{S}S'$. 
By abuse, we sometimes omit some of the subscripts, writing for 
instance $\mor_{S'}\,(X,Y)$ for $\mor_{S'}\,(X_{S'},Y_{S'})$.

In these notations,  affine base schemes
are often denoted by the corresponding ring: thus, if  $S=\Spec(F)$, we may write 
$\mor_{F}\,(X,Y)$ rather than $\mor_{S}\,(X,Y)$.

\Subsection{Involutions, odd morphisms, algebraic groups. }
\label{ParInvol}
If $X$ and $Y$ as above are provided with $S$-involutions $\sigma$ 
and $\tau$ respectively, we shall denote by $\modd_{S}(X,Y)$ the set of 
$S$-morphisms $\varphi:X\to Y$ such that 
$\varphi\circ\sigma=\tau\circ\varphi$.

The involutions considered will in general be clear from the context. 
In particular, if $Y$, say, is a commutative $S$-group scheme, 
written additively, the involution on $Y$ will be multiplication by $-1$, 
unless otherwise specified.

If $X$ and $Y$ are commutative $S$-group schemes (always assumed to be 
separated and of finite presentation as $S$-schemes), we shall denote 
by $\Hom_{S}\,(X,Y)$ the set of morphisms \emph{of  
$S$-group schemes} from $X$ to $Y$ (a subgroup of 
$\modd_{S}(X,Y)$), and we write ${\rm End}_{S}(X)$ for 
$\Hom_{S}\,(X,X)$.

If $G$ is a commutative $S$-group scheme, and $n\in\ZZ$, 
we denote by $[n]_{G}$ or $[n]$ the 
endomorphism of $G$ given 
by multiplication by $n$, and by $G[n]$ its kernel. If $n$ is 
\emph{invertible on $S$}, then $[n]$ is an unramified morphism (in fact it 
is \etale\ along the unit section of $G$), and $G[n]$ can therefore 
be written as the disjoint union of the unit section and a subscheme 
$G[n]^\ast$: in particular, if $n$ is prime (the only case we shall 
use is $n=2$), we may define $G[n]^\ast$ as the subscheme of $G$ of `points of 
exact order $n$'.

We shall use some `rigidity' properties of odd morphisms to a 
commutative group scheme:

\Subsubsection{Rigidity of odd morphisms: notations. }\label{NotRigOdd}
	We assume that $2$ is invertible on $S$ (equivalently, all residue 
	characteristics of points of $S$ are different from $2$).
	
	Let $f:X\to S$ be a morphism of schemes. Assume that $f$ is flat, 
	proper, of finite presentation and that 
	$$\calO_{S}\fflis f_{\ast}{\calO_{X}}\text{ universally}$$
(that is, $\calO_{S'}\fflis (f_{S'})_{\ast}{\calO_{X_{S'}}}$ for every 
$S$-scheme $S'$; these conditions are satisfied in particular when 
$S$ is Noetherian, $f$ is projective and flat, and all its geometric 
fibres are irreducible and reduced). Let $\sigma$ be an $S$-involution of $X$.

We also fix a commutative $S$-group scheme $G\to S$, separated and of 
finite presentation as an $S$-scheme. We write $G$ additively.

\begin{subprop}\label{RigOdd1} With the assumptions of\/ 
\rref{NotRigOdd},
	let $u:X\to G$ be an \emph{odd} $S$-morphism. Then, the following 
	conditions are equivalent:
	\begin{romlist}
		\item\label{RigOdd11} $u=0$;
		\item\label{RigOdd12} $u=0$ set-theoretically (that is, $u$ maps the underlying 
		space of $X$ to the unit section of $G$);
		\item\label{RigOdd13} for every point $s$ of $S$, $u(X_{s})$ is contained in an 
		affine open subset of $G_{s}$ disjoint from $G_{s}[2]^\ast$.
	\end{romlist}
	
	Moreover, the set $\Sigma:=\{s\in S\mid u_{s}:X_{s}\to G_{s}\text{ 
	\rm is 
	zero }\}$ is open and closed in $S$.
\end{subprop}
\dem It is trivial that 
\ref{RigOdd11}$\Impl$\ref{RigOdd12}$\Impl$\ref{RigOdd13}. Let us 
prove \ref{RigOdd13}$\Impl$\ref{RigOdd12}. We need to show that, for 
every $s\in S$ (with residue field $\kappa(s)$), $u$ maps $X_{s}$ to 
the unit $0_{s}$ of $G_{s}$. By \ref{RigOdd13}, $u$ maps $X_{s}$ to an affine 
scheme; but any morphism from $X_{s}$ to an affine scheme 
must factor through $\Spec(\Gamma(X_{s},\calO_{X_{s}}))$ which is 
$\Spec\kappa(s)$ by our assumptions on $X$. In other words,  $u$ maps 
$X_{s}$ to a rational point $\gamma$ of $G_{s}$ which must be fixed by $[-1]$ 
since $u$ is odd; by the assumption \ref{RigOdd13} we must have 
$\gamma=0_{s}$.

To prove \ref{RigOdd12}$\Impl$\ref{RigOdd11} we may assume that $S$ 
is affine, and (localising further if necessary) that there exists an 
affine open neighbourhood $U$ of the unit section of $G$, disjoint from 
$G[2]^\ast$. Clearly, \ref{RigOdd12} implies that $u$ factors 
through $U$. As above, this implies that $u$ must factor as 
$\gamma\circ f:X\to S\to G$, where $\gamma\in G(S)$ is a 
section which must be fixed by $[-1]$, hence equal to the unit 
section by our assumption on $U$.

Let us now prove the last claim. To see that $\Sigma$ is open, let us take 
$s$ in $\Sigma$, and show that $u=0$ over a neighbourhood of $s$. We may assume that 
there is an open subset $U$ of $G$ as in the proof of 
\ref{RigOdd12}$\Impl$\ref{RigOdd11} above. Then $u^{-1}(U)$ is an 
open subscheme of $X$ containing $X_{s}$; since $f$ is proper (hence 
closed), it must 
contain $f^{-1}(V)$ for some neighbourhood $V$ of $s$ in $S$. But 
then $f_{V}:X_{V}\to G_{V}$ factors through $U$, hence is zero.

To see that $\Sigma$ is closed, consider the inverse image in $X$ of 
the unit section of $G$: since $G$ is separated, this is a closed subscheme 
of $X$, hence its complement $W\subset X$ is open, and so is 
$f(W)\subset S$ ($f$ is flat of finite presentation, hence open). 
But the complement of $f(W)$ is just $\Sigma$.\qed

\begin{subcor}\label{RigOdd2} \emph{(Rigidity of odd morphisms)} With the assumptions of 
	\rref{NotRigOdd}, let $T\to S$ be a faithfully flat quasicompact $S$-scheme, 
	with \emph{geometrically connected fibres}. Then the natural `base change' map
	$$\modd_{S}(X,G)\ffl \modd_{T}(X_{T},G_{T})$$
	is an isomorphism.
\end{subcor}
\dem Injectivity is clear since $T\to S$ is faithfully flat. Let 
$u_{T}:X_{T}\to G_{T}$ be an odd $T$-morphism. By flat descent, it is 
enough to show that the two morphisms $u_{1}, u_{2}:X_{T\times_{S} T}\to 
G_{T\times_{S} T}$ deduced from $u_{T}$ by base change via the two 
projections $T\times_{S} T\to T$ are equal. But clearly they coincide 
along the diagonal $T\to T\times_{S}T$, hence also, by \ref{RigOdd2}, 
along an open and closed subscheme of $T\times_{S} T$ containing the 
diagonal. By our assumptions, the only such subscheme is 
$T\times_{S}T$, and the corollary is proved.\qed

\begin{subrem}\label{RemRigOdd}
	An interesting special case is when $S=\Spec(F)$ and $T=\Spec(F')$ 
	where $F'$ is an extension of $F$. Then $T$ is geometrically 
	connected if and only if $F$ is \emph{separably closed in} $F'$. 
	Thus, in this case, we have $\modd_{F'}(X,G)=\modd_{F}(X,G)$.
\end{subrem}

\Subsection{Existence of admissible morphisms on curves, and of odd 
projections on varieties. 
}\label{BasicCurves}
Let $C$ be a projective, smooth, 
geometrically connected $F$-curve, of genus $g$. If $D$ is a divisor 
on $C$, we denote by $\vert D\vert$ the linear system associated to 
$D$, i.e.~the space of effective divisors linearly 
equivalent to $D$. In other words, $\vert D\vert=\PP({\rm 
H}^0(C,\mathcal{O}_C(D))^{*}$ is the projective space of lines 
of the $F$-vector space $\mathcal{L}(D)={\rm H}^0(C,\mathcal{O}_C(D))$.

Put $d=\deg D$. By Riemann-Roch, $\dim\vert D\vert\geq d-g$, with 
equality if $d\geq2g-1$.

\begin{subprop}\label{ExistAdm} With the above assumptions, assume 
that $F$ is infinite, and that $d=\deg D\geq2g+2$. Then there exists 
an $F$-morphism $f:C\to\Pu{F}$, having only simple ramification and 
such that the divisors $f^{-1}(\xi)$, for $\xi\in\Pu{F}$, belong to $\vert 
D\vert$ (equivalently, the invertible sheaf 
$f^*\mathcal{O}_{\Pu{F}}(1)$ is isomorphic to $\mathcal{O}_C(D)$). In 
particular, $\deg f=d$.

Moreover, assume that some $F$-rational $E_{0}\in|D|$ is fixed, without triple 
points (over $\ol{F}$). Then one can choose $f$ as above such that 
$f^{-1}(0)=E_{0}$ (as divisors).
\end{subprop}
\dem The case $g=0$ is left to the reader; we assume $g>0$ and 
in particular $d\geq4$. 

For $n\in\NN$, denote by $C^{(n)}$ the $n$-th symmetric power of 
$C$. It is a smooth projective $F$-scheme, whose $\ol{F}$-points 
correspond canonically to effective divisors of degree $n$ on 
$C_{\ol{F}}$. We define a closed subvariety $W$ of $C\times 
C^{(d-3)}$ by
$$W:=\{(P,D_{1})\,\mid\,3P+D_{1}\sim D\}$$
where $\sim$ denotes linear equivalence. There is an obvious 
morphism $W\to|D|$ sending $(P,D_{1})$ to $3P+D_{1}$, whose image $W'$ 
consists of divisors having a triple point. 

I claim that $\dim W=d-g-2$. Indeed, consider the natural 
projection $W\to C$ sending $(P,D_{1})$ to $P$. The fibre of 
a point $P\in C$ is canonically 
isomorphic to the projective space $|D-3P|$, which has dimension 
$d-3-g$ since by assumption $d-3\geq2g-1$. Hence $\dim W=\dim 
C+d-g-3=d-g-2$, as claimed. In particular, $W'$ has codimension $\geq2$ in $|D|$. 

Since $F$ is infinite, we can choose an $F$-rational point $E_{0}$ in 
$|D|\setminus W'$ (for instance the given one, if necessary). Write 
$E_{0}=P_{1}+\cdots+P_{d}$ (over $\ol{F}$). For each 
$i\in\{1,\ldots,d\}$, divisors in $|D|$ containing $P_{i}$ form a 
linear subspace isomorphic to $|D-P_{i}|$, hence a hyperplane 
$H_{i}$ of $|D|$ (not necessarily defined over $F$). The divisors 
meeting $E_{0}$ thus form a hypersurface $H$ (the union of the $H_{i}$'s, 
which is defined over $F$). Again, since $F$ is infinite, there is a 
line $\Delta\subset|D|$ containing $E_{0}$, disjoint from $W'$ and 
not contained in $H$. Take an $F$-rational point $E_{\infty}$ on $\Delta$, 
distinct from $E_{0}$ (hence not in $H$, because $\Delta\cap 
H_{i}=\{E_{0}\}$ since $\Delta$ is a line). There is a rational 
function $f$ on $C$ with divisor $E_{0}-E_{\infty}$, whence an 
$F$-morphism $\varphi:C\to\Pu{F}$ whose fibres are precisely the 
points of $\Delta$. This $f$ satisfies the required conditions.\qed

\begin{subrem} The result should also be true if $F$ is finite, 
possibly with a stronger condition on the degree.
\end{subrem}

\begin{subrem} Proposition \ref{ExistAdm} clearly implies the 
existence of admissible morphisms, in the sense of \ref{DefGood}. The 
extra information on the degree will be used in Section 
\ref{Secpadic}, via the following special case: if $C$ admits a divisor of odd 
degree, then there is an admissible morphism $C\to\Pu{F}$ of odd 
degree. In the same vein, we shall also need the next proposition.
\end{subrem}

\begin{subprop}\label{ProjectImpaire} Assume that $\car F=0$, and let $K$ 
	be a finitely generated, regular, transcendental extesion 
	of $F$. Then there is a transcendence basis $(z_{1},\ldots,z_{n})$ of 
	$K$ over $F$ such that 
	$K/F(z_{1},\ldots,z_{n-1})$ is a regular exension. 
	
	If, moreover, $F$ is algebraically closed, then $(z_{1},\ldots,z_{n})$ 
	may be chosen such that, in addition, $[K:F(z_{1},\ldots,z_{n})]$ is odd.
\end{subprop}
\dem Let $V\subset\PP^{n+1}_{F}$ be a projective hypersurface with function 
field $K$, and put $d:=\deg V$. Since $K/F$ is regular, $V$ is 
geometrically integral (i.e.~$V_{\ol{F}}$ is irreducible and reduced).

By Bertini's theorem, there is a plane $\Pi\subset\PP^{n+1}_{F}$ such that 
$\Pi\cap V$ is a geometrically  integral curve; in fact, this property holds 
for all $\Pi$ in a dense open subset of the Grassmannian of planes.

Take an $F$-rational point $A\in\Pi$, not 
in $V$. Consider the projection $\pi$ from $V$ to the  
space $S_{A}$ (isomorphic to $\PP^n_{F}$) of lines through $A$. This 
$\pi$ sends a point $P\in V$ to 
the line through $P$ and $A$, and is clearly a finite surjective morphism of degree 
$d$. Moreover, $\Pi\cap V=\pi^{-1}(L)$ where the line $L\subset S_{A}$ is the 
image of $\Pi$. We can choose coordinates $z_{1},\ldots,z_{n}$ on $S_{A}$ 
such that $L$ is defined by, say, $z_{1}=\ldots=z_{n-1}=0$. The 
rational map $\varphi:=(z_{1},\ldots,z_{n-1}): V\cdots\fl\Aa^{n-1}_{F}$ then 
has the property that for all $\xi=(\xi_{1},\ldots,\xi_{n-1})\in \ol{F}^{n-1}$ 
except in a proper Zariski closed subset, $\varphi^{-1}(\xi)$ is an 
integral curve. This implies that the generic fibre of $\varphi$ is 
geometrically integral, i.e.~that the extension $K/F(z_{1},\ldots,z_{n-1})$ is regular.
\smallskip

Note that in the previous construction, $[K:F(z_{1},\ldots,z_{n})]=d$. Hence 
to prove the last assertion, we assume $d$ even and $F$ algebraically 
closed. Again we use a projection, but this time we take a smooth 
point $A\in\Pi\cap V$ and project from $A$. This defines a morphism 
$\pi:V\setminus\{A\}\to S_{A}$; for a general point of $S_{A}$, 
corresponding to a line $l$ through $A$, the fibre $\pi^{-1}(l)$ 
consists of the $d-1$ points of $l\cap V$ distinct from $A$, hence 
$\pi$ has degree $d-1$, which is odd. 

Next, defining $L\subset S_{A}$ as above, we have 
$\pi^{-1}(L)=(\Pi\cap V)\setminus\{A\}$ which is geometrically integral, and the 
same property holds for all lines $L'$ in a dense open subset of the 
Grassmannian of lines in $S_{A}$. Choosing coordinates as above, we 
again conclude that $K/F(z_{1},\ldots,z_{n-1})$ is regular.\qed

\Subsection{Abelian varieties and schemes. }\label{ParVarAb}
An \emph{abelian scheme\/} over a scheme $S$ is a smooth proper 
$S$-group scheme $A\to S$, with connected fibres. An abelian scheme 
over $\Spec(F)$ is called an abelian variety over $F$. Abelian schemes 
are automatically commutative, and abelian varieties are projective.

\begin{subprop}\label{RigSousVarAb} Let $A$ be an abelian variety over $F$, $F'$ an 
extension of $F$, and $B\subset A_{F'}$ an abelian subvariety. Then 
$B$ can be defined over a finite extension of $F$.
\end{subprop}
\dem We assume $\car(F)\neq2$, which is sufficient for our purposes 
(but the result holds in general). We may assume $F$ algebraically closed. 
By standard arguments, there is a finitely 
generated $F$-subalgebra $R\subset F'$ and an abelian subscheme 
$\mathcal{B}\subset A_{S}$ over $S=\Spec(R)$ such that 
$\mathcal{B}_{F'}=B$. Now $S$ has an $F$-rational point $x$; let 
$B_{0}$ be the fibre of $\mathcal{B}$ at $x$, an abelian subvariety 
of $A$. Consider the natural $S$-morphism $\mathcal{B}\inj A_{S}\surj 
(A/A_{0})_{S}$: it is zero at $x$, hence zero by \ref{RigOdd1}. Thus, 
$\mathcal{B}\subset A_{0,S}$. Since  the opposite inclusion is proved 
similarly, we have equality, and in particular $B=A_{0,F'}$.\qed

\Subsubsection{Torsion points and Tate modules. }\label{TorsVarAb}
If $A$ is an abelian scheme over $S$, of relative dimension $g$, then
$[n]_{A}$ is a finite locally free morphism of degree $n^{2g}$, 
\'{e}tale above all points of $S$ whose residue characteristic does not 
divide $n$. 

If $S=\Spec(F)$ and ${\rm char\,}(F)\nmid n$, then $A[n](\ol{F})$ 
is isomorphic to $(\ZZ/n\ZZ)^{2g}$. If $l\neq{\rm char\,}(F)$ is a prime 
number, we define the \emph{$l$-adic Tate module\/} $T_{l}(A)$ of $A$ by 
\begin{equation}\label{EqModTate}
	T_{l}(A):=\varprojlim_{n\geq1}A[l^n](\ol{F})
\end{equation}
with transition maps induced by $A[l^{n+1}]\varfl{\times l}A[l^{n}]$. 
This is a free $\ZZ_{l}$-module of rank $2g$, with a continuous 
action of $G_{F}={\rm Gal\,}({F}^{\rm s}/F)$; by the way, note that 
$A[l^n](\ol{F})=A[l^n]({F}^{\rm s})$ since $A[l^n]$ is 
\'{e}tale over $F$. 

It is often convenient to use the corresponding $\QQ_{l}$-vector 
space:
\begin{equation}\label{EqVecTate}
	V_{l}(A):=\QQ_{l}\otimes_{\ZZ_{l}}T_{l}(A)
\end{equation}
which defines a $2g$-dimensional continuous representation of $G_{F}$ over 
$\QQ_{l}$.

\Subsubsection{Homomorphisms. }\label{HomVarAb}
$T_{l}(A)$ is obviously functorial in $A$: thus, if $B$ is another 
abelian variety over $F$, there is a canonical homomorphism
\begin{equation}\label{EqHomTate}
	\Hom_{F}\,(A,B)\otimes_{\ZZ}\ZZ_{l}\ffl
	\Hom_{\ZZ_{l}[G_{F}]}\,(T_{l}(A),T_{l}(B)).
\end{equation}
We have the following properties (for the first three, see for instance \cite{Mum}, 
\S19):
\begin{romlist}
	\item\label{PropHom1} The homomorphism (\ref{EqHomTate}) is 
	injective.
	\item\label{PropHom2} $\Hom_{F}\,(A,B)$ is a free finitely 
	generated abelian group.
	\item\label{PropHom3} If $S$ is a geometrically connected $F$-scheme, 	
	then $\Hom_{S}\,(A,B)=\Hom_{F}\,(A,B)$. (If $\car(F)\neq2$, this can 
	also be deduced from \ref{RigOdd2}). In particular, if $\Omega$ is any extension of $F^{\rm s}$, 
	then $\Hom_{\Omega}\,(A,B)=\Hom_{F^{\rm s}}\,(A,B)$. 
	\item\label{PropHom4} If $F$ is 
finitely generated over the prime field, then  (\ref{EqHomTate}) is an isomorphism. 
\pauseromlist
Property \ref{PropHom4} is Tate's conjecture for homomorphisms of abelian 
varieties, proved by Fal\-tings: for a proof, see  
\cite{FW}, VI, \S3, Theorem 1 (where it is stated for $A=B$, which 
implies the general case by considering products). 

We shall also use \ref{PropHom4} in the following (seemingly) weaker form: consider 
the natural homomorphism
\begin{equation}\label{EqHomTate2}
	\Hom_{\ol{F}}\,(A,B)\otimes_{\ZZ}\QQ_{l}\ffl
	H_{l}(A,B):=\Hom_{\QQ_{l}}\,(V_{l}(A),V_{l}(B))
\end{equation}
which is the map (\ref{EqHomTate}), taken  over $\ol{F}$ and 
tensored with $\QQ_{l}$, and is therefore injective. Faltings' 
theorem implies (and, in fact, is equivalent to):
\finpauseromlist
	\item\label{Faltings} Assume $F$ is finitely generated over the 
	prime field. Then the image of (\ref{EqHomTate2}) consists of 
	those elements of $H_{l}(A,B)$ whose stabiliser in $G_{F}$ 
	is open. In particular, it is determined by the image of $G_{F}$ 
	in ${\rm Aut}_{\QQ_{l}}(H_{l}(A,B))$.
\end{romlist}

\Subsubsection{Elliptic curves. }\label{ParEll}
An \emph{elliptic curve\/} over a scheme $S$ is a pair $(E,\omega)$ where $E$ 
is a smooth proper $S$-scheme whose fibres are curves of genus $1$ 
and $\omega$ is a section of $E$ over $S$. We shall often drop $\omega$ from the
notation. Recall that there is a unique (commutative) $S$-group scheme structure 
on $E$ with unit section $\omega$: thus, we may 
equivalently define an elliptic curve to be an abelian scheme of 
relative dimension $1$.

\Subsection{Picard groups and schemes. }\label{Pic}
For any scheme $X$, the Picard group of $X$, denoted by ${\rm 
Pic\,}(X)$, is the group of isomorphism classes of invertible ${\cal 
O}_{X}$-modules. In good cases, this coincides with the group of 
Cartier (i.e.~locally principal) divisors on $X$, modulo principal 
divisors: this is the case in particular if $X$ is 
quasiprojective over $F$, by (\cite{ega4}, 21.3.4). If $X$ is regular,  
${\rm Pic\,}(X)$ is just the 
usual group of divisor classes on $X$.

Now let $f:X\to S$ be a morphism of schemes. We assume that $f$ is 
proper and flat, and $f_{\ast}{\cal O}_{X}\cong{\cal O}_{S}$ 
universally. In our applications, $S$ will be either an integral regular 
scheme of dimension $1$ (e.g.~a nonsingular curve), or the spectrum of 
a field $F$. In the latter case, the condition on $f_{\ast}{\cal O}_{X}$ 
means that ${\Gamma}(X,{\cal O}_{X})=F$; this holds whenever $X$ is 
geometrically integral over $F$.

One can then define the \emph{Picard functor\/} 
$\soul{\rm Pic}_{X/S}$ of $f$. As a `first approximation', we define the 
\emph{naive Picard functor\/} from  $S$-schemes to abelian groups, by
\begin{equation}\label{PicNaif}
	\picn_{X/S}(T):=\pic\,(X\times_{S}T)/\mathrm{pr}_{2}^{\ast}\,(\pic\,(T)).
\end{equation}
Now $\soul{\rm Pic}_{X/S}$ is another contravariant 
functor from $S$-schemes to abelian groups, with the following 
properties (for which we refer to \cite{BLR}, Chapter 8):
\begin{romlist} 
	\item\label{Pic0} For any $S$-scheme $T$, there is an  
	injective homomorphism, functorial in $T$:
	\begin{equation}\label{EqPicnToPic}
		a_{X/S}(T):\quad\picn_{X/S}(T)\ffl\soul{\rm Pic}_{X/S}(T).	
	\end{equation}
	\item\label{Pic1} $\soul{\rm Pic}_{X/S}$ `commutes 
	with any base change' $S'\to S$, in the sense that if $T$ is an 
	$S'$-scheme and $X'=X\times_{S}S'$ then $\soul{\rm 
	Pic}_{X/S}\,(T)\cong\soul{\rm Pic}_{X'/S'}\,(T)$ (where, in the 
	left-hand side, $T$ is viewed as an $S$-scheme in the natural way).  
	\item\label{Pic2} If $f$ has a section $\varepsilon:S\to X$, then 
	the morphism $a_{X/S}$ is an isomorphism of functors on $S$-schemes. 
	Moreover, in this case the functors $\picn_{X/S}$ and $\soul{\rm Pic}_{X/S}$ 
	are isomorphic to the functor
	$$T\mapsto \pic_{X/S}^{\varepsilon}(T):=
	\ker{\left[{\rm Pic\,}(X\times_{S}T) 
	\varfl{(\varepsilon\times{\rm Id}_{T})^{\ast}} {\rm 
	Pic\,}(T)\right]}$$ 
	of `invertibles sheaves on $X$ which are trivial along $\varepsilon$'.
	\item\label{Pic3} If $S=\Spec(F)$, then 
	$\soul{\rm Pic}_{X/F}$ is representable by an $F$-group scheme 
	locally of finite type, whose connected component (denoted by 
	$\soul{\rm Pic}^0_{X/F}$) is an algebraic group over $F$, which is 
	smooth if ${\rm char\,}(F)=0$ or $\dim X=1$.  
	\item\label{Pic4} If $S=\Spec(F)$ and $X$ is smooth over $F$, then 
	$\soul{\rm Pic}^0_{X/F}$ is proper (hence an abelian variety if 
	${\rm char\,}(F)=0$ or $\dim X=1$).
	\item\label{Pic5} If $S=\Spec(F)$ and $X$ is a semistable 
	curve (that is, $X_{\ol{F}}$ has only ordinary double points as 
	singularities), then $\soul{\rm Pic}^0_{X/F}$ is \emph{semiabelian\/}, 
	i.e.~an extension of an abelian variety by a torus.  
\pauseromlist
Observe that if $S=\Spec(F)$, and $X$ has an $F$-rational point 
(which is automatic when $F$ is algebraically closed), \ref{Pic2} implies 
$\soul{\rm Pic}_{X/F}\,(F)\cong {\rm Pic\,}(X)$. In this case, 
$\soul{\rm Pic}^0_{X/F}\,(F)$ is the subgroup consisting of (classes of) 
invertible sheaves algebraically equivalent to zero; in particular, if $X$ is a 
(not necessarily irreducible) curve, then $\soul{\rm Pic}^0_{X/F}\,(F)$ is 
the group of invertible sheaves having degree $0$ on each component 
of $X$.

Property \ref{Pic3}, in the general case, is due to Murre; we shall use it only 
when $X$ is a curve or a nonsingular surface. The situation is more 
delicate over a more general base $S$. But of course, we can apply 
\ref{Pic3} to the fibres of $f$, which at least allows us to \emph{define\/} a 
subfunctor $\soul{\rm Pic}^0_{X/S}$ of $\soul{\rm Pic}_{X/S}$ by 
$$\soul{\rm Pic}^0_{X/S}\,(T)=\{x\in\soul{\rm 
Pic}_{X/S}\,(T)\,|\,\forall 
t\in T, x_{t}\in\soul{\rm Pic}^0_{X_{t}/\kappa(t)}(\kappa(t))\}$$
where, as usual, $\kappa(t)$ is the residue field of $t$ and the 
subscript $t$ 
means base change by $\Spec\kappa(t)\to T\to S$. So, loosely 
speaking, 
$\soul{\rm Pic}^0_{X/S}$ parametrises invertible sheaves on $X$ 
which are algebraically equivalent to zero in the fibres of $f$.

There are deep representability results for $\soul{\rm Pic}$ and 
$\soul{\rm Pic}^0$; 
we shall only need the following special case, due to Raynaud: 
\finpauseromlist
\item\label{Pic6} In addition to our general hypotheses, assume that 
	$S$ is a regular 
	integral scheme of dimension $1$, that $X$ is normal, and that each 
	geometric fibre of $f$ is a curve with at least one reduced 
	irreducible component. 
	Then $\soul{\rm Pic}^0_{X/S}$ is representable by a smooth separated 
	$S$-group scheme. 
\pauseromlist
The representability is a special case of Theorem 2 of \cite{BLR}, 
9.4 
(which is stated over a discrete valuation ring, but the extension to 
our case 
is standard). Smoothness is automatic for the Picard functor of a 
curve, 
as explained in \cite{BLR}, 8.4, Proposition 2 (essentially, the 
reason is that ${\rm H}^2(X_{s},{\cal O}_{X_{s}})=0$ for all $s\in 
S$).
\medskip

Finally we shall need two well-known facts about Picard groups of 
surfaces. The first one is the birational invariance of $\soul{\rm 
Pic}^0$:
\finpauseromlist
	\item\label{Pic7} Let $Z$ be a smooth projective 
	geometrically connected surface 
	over $F$, and let $\rho:Z'\to Z$ be the blowing-up of finitely many 
	(reduced) points. Then $\rho^\ast$ induces an isomorphism 
	$\soul{\rm Pic}^0_{Z/F}\flis\soul{\rm Pic}^0_{Z'/F}$.
\pauseromlist
The other result we need is the structure of $\soul{\rm 
Pic}^0$ of a product surface:
\finpauseromlist
	\item\label{Pic8} Let $X$, $Y$ be two smooth projective 
	geometrically connected 
	varieties over $F$. Then the natural morphism
	$${\rm pr}_{1}^\ast\oplus{\rm pr}_{2}^\ast:
	\soul{\rm Pic}^0_{X/F}\times_{F}\soul{\rm Pic}^0_{Y/F}
	\ffl \soul{\rm Pic}^0_{(X\times_{F}Y)/F}$$
	is an isomorphism.
\end{romlist}
(Note that the analogues of \ref{Pic7} and \ref{Pic8} where 
$\soul{\rm 
Pic}^0$ is replaced by $\soul{\rm Pic}$ are false).

\Subsection{Jacobians. }\label{Jac} 
If $X$ is a projective geometrically connected curve over $F$, 
we denote by $\jac(X)$ the  Jacobian of $X$, which is by 
definition $\soul{\rm Pic}^0_{X/F}$. 

If $X$ is smooth, it follows from \ref{Pic} that $\jac(X)$ is an abelian 
variety over $F$, and that if $L$ is an extension of 
$F$ such that $X(L)\neq\emptyset$, then $\jac(X)(L)=\jac(X_{L})(L)={\rm 
Pic}^0(X_{L})$, the group of divisor classes of degree zero on 
$X_{L}$. 

If $X$ is semistable, then $\jac(X)$ is semiabelian, by 
\ref{Pic}\,\ref{Pic5}. 

If $(E,\omega)$ is an elliptic curve, there is a canonical 
isomorphism 
$E\flis \jac(E)$ sending a point $x$ to the divisor class 
$[x]-[\omega]$. We shall henceforth identify $E$ and $\jac(E)$ in 
this 
way.

If $f:X\to Y$ is a $k$-morphism of smooth projective geometrically 
connected 
curves, there is a pullback morphism $f^\ast:\jac(Y)\to \jac(X)$, 
corresponding to the usual pullback of divisors. 
If $Y=E$ is an elliptic curve, the map $f\mapsto f^*$ induces a 
\emph{group homomorphism}
\begin{equation}
	\label{MorJac}
	\mor_{F}\,(X,E)/E(F)\ffl\Hom_{F}\,(E, \jac(X))
\end{equation}
where of course we identify $E(F)$ with the group of constant 
morphisms from $X$ to $E$ (or, alternatively, with the group of 
translations on $E$ acting on $\mor_{F}\,(X,E)$). This map is 
always injective (to see this, use the fact that a morphism $f:X\to 
E$ 
also induces $f_{\ast}:\jac(X)\to E$ satisfying $f_{\ast}\circ 
f^{\ast}=(\deg f)\,{\rm Id}_{E}$). Furthermore, if $X$ has a rational 
point, then (\ref{MorJac}) is bijective: this can be seen using 
self-duality of the Jacobian and the embedding of $X$ into $J(X)$ 
attached to a rational point of $X$ (if $X$ has genus $\geq1$; 
otherwise both sides of (\ref{MorJac}) are zero).

With $X$ and $E$ as above, assume now that $X$ is endowed with an 
involution $\sigma$ (and $E$ with the involution $[-1]$). Then 
(\ref{MorJac}) induces an injective homomorphism
\begin{equation}
	\label{MorJacOdd}
	\modd_{F}(X,E)/E[2](F)\ffl\Hom_{F}^{\rm odd}(E, 
	\jac(X))
\end{equation}
where, in the right-hand side, $E$ (resp.~$\jac(X)$) is given the 
involution $[-1]$ (resp.~$\sigma^\ast$). 
\begin{subprop}\label{OddEtJac} If $X$ has an $F$-rational point
 fixed by $\sigma$, then {\rm(\ref{MorJacOdd})} is an isomorphism.
\end{subprop}
\dem Let $P\in X(F)$ be such a point, and let $v:E\to\jac(X)$ be an 
odd morphism. Since (\ref{MorJac}) is bijective, there is a morphism 
$u:X\to E$ such that $u^*=v$. Changing $u$ by a translation on $E$ 
we may assume $u(P)=0$ (the origin of $E$). On the other hand, the fact that 
$u^*$ is odd means that $(u+u\circ\sigma)^*=0$, hence $u+u\circ\sigma: X\to E$ 
is constant. By our assumptions it sends $P$ to $0$, hence $u+u\circ\sigma=0$ 
and $u$ is odd.\qed
\begin{subrem}\label{RemStructMorodd} It follows in particular that 
	$\modd_{F}(X,E)$ is a finitely generated abelian group, 
	with torsion subgroup $E[2](F)$ and rank $\leq4\,\text{genus}(X)$.
\end{subrem}

\Subsection{Affine \dio\ sets. }\label{SsecAffDio}

Throughout this section, we denote by $R$ a ring and by $\calO$ an 
$R$-algebra. 

We denote by $\lr=\{+,-,.,0,1\}$ the language of 
rings, and by $\lr(R)$ the language $\lr$ augmented by the set of 
constants $R$. 

\Subsubsection{The case of affine $n$-space. }\label{SsecDioEspAff}
Let $n$ be a natural integer. A subset 
$X$ of $\calO^n$ will be called \emph{primitive \dio} (with respect to 
$R$) if there is an integer $q$ and a finite sequence of polynomials 
$F_{1},\ldots,F_{r}\in R[T_{1},\ldots,T_{n},Y_{1},\ldots,Y_{q}]$ such 
that, for  $\soul{t}=(t_{1},\ldots,t_{n})\in \calO^n$, we have the 
equivalence
\begin{equation}\label{EqAffine2}
	\soul{t}\in X\;\Leftrightarrow\;\hfill
	\exists\,\soul{y}=(y_{1},\ldots,y_{q})
	\in \calO^{q}\hbox{ such that }F_{1}(\soul{t},\soul{y})=\cdots=
	F_{r}(\soul{t},\soul{y})=0.
\end{equation}
With these notations, consider the $R$-scheme 
$W:=\Spec(R[\soul{T},\soul{Y}]/(F_{1},\ldots,F_{r}))\inj\Aa^{n+q}_{R}$, 
and the $R$-morphism $\varphi:W\to\Aa^n_{R}$ deduced from the 
projection $(\soul{t},\soul{y})\mapsto\soul{t}$: then $X$ is simply the 
image of the map $\varphi(\calO):W(\calO)\to\Aa^n_{R}(\calO)=\calO^n$ given by $\varphi$. 
Conversely, it is easy to see that if $W$ is any affine $R$-scheme of 
finite presentation (that is, of the form $\Spec(A)$ where $A$ is a 
finitely presented $R$-algebra), and $\varphi:W\to\Aa^n_{R}$ is an 
$R$-morphism, the image of $\varphi(\calO):W(\calO)\to \calO^n$ is primitive \dio.

A \emph{\dio\ subset} of $\calO^n$ is by definition a finite union of 
primitive \dio\ subsets. It is well known, and easy to see, that 
these subsets are precisely the \emph{positive-existentially definable} 
subsets of $\calO^n$, in the language $\lr(R)$.
\bigskip

The above remarks justify the following definition: 

\begin{subdefi}\label{DefDioAff}
	Let $V$ be an 
	affine $R$-scheme of finite presentation, and let $\calO$ be an $R$-algebra.
	
	A subset $X$ of $V(\calO)$ is called \emph{primitive \dio} (with respect to 
	$R$) if there is an affine $R$-scheme $W$ of finite presentation, and 
	an $R$-morphism $\varphi:W\to V$, such that $X$ is the image of 
	$\varphi(\calO):W(\calO)\to V(\calO)$. 
	
	A \emph{\dio} subset of $V(\calO)$ is a finite union 
	of primitive \dio\ subsets. 
\end{subdefi}

The class of \dio\ subsets is closed under finite intersections and 
unions, images and inverse images by $R$-morphisms, and various other 
operations (such as fibre products); these properties are easy to check. 

In particular, if $V\subset V'$ is an immersion of affine $R$-schemes 
of finite presentation, then $X\subset V(\calO)$ is \dio\ if and only if 
it is \dio\ as a subset of $V'(\calO)$. The most important case is, of 
course, when $V'=\Aa^n_{R}$: thus, in this case, the \dio\ subsets of  $V(\calO)$ 
are those which are positive-existentially definable in $\calO^n$.

\begin{subrem}\label{DioPasPrimDio}
	It is not true in general that \dio\ sets are primitive \dio; however, 
	this does hold if $\Spec(\calO)$ is connected (or, equivalently, if $\calO$ 
	has no idempotent element other than $0$ and $1$). In particular, 
	this is true if $\calO$ is a domain, or a local ring.
\end{subrem}

\begin{subrem}\label{RemDioNonAffine}
	It is of course natural to ask whether Definition \ref{DefDioAff} 
	generalises to schemes $V$ which are not necessarily affine.  
	
	First, one should probably keep the `finite presentation' 
	restriction on $V$: loosely speaking, $R$-schemes of finite presentation 
	are those which can be defined by a finite set of data from $R$.
	
	Next, of course the extended notion should specialise to the previous 
	one for affine $V$; hence, it also seems reasonable to `use only affine $W$'s' in the 
	definition. 
	
	So, our definition of a primitive \dio\ subset of $V(\calO)$ (for an 
	$R$-scheme $V$ of finite presentation) would be a subset which is the 
	image of the map $W(\calO)\to V(\calO)$ induced by an $R$-morphism $W\to V$, 
	where $W$ is some \emph{affine} $R$-scheme of finite presentation. 
	Of course, a \dio\ subset is a finite union of primitive \dio\ subsets.
	
	The above definition differs from Mazur's (\cite{mazur}, Definition 1), 
	who defines a \dio\ subset of $V(\calO)$ as one 
	which is the image of $W(\calO)$ for some $R$-morphism $W\to V$ of 
	$R$-schemes of finite type. 
	
	In any case, we have refrained from including the basic properties of 
	these generalised \dio\ sets in this paper. Of course, in a sense, this 
	would have been the natural framework when working with elliptic curves; 
	but as it turns out, elliptic curves contain very nice affine open subsets 
	which are sufficient for our needs.
\end{subrem}

\Subsubsection{\dio\ relations and maps. }\label{ApplDio}
If $V$, $V'$ are affine $R$-schemes of finite presentation, and $X$, 
$X'$ are \dio\ subsets of $V(\calO)$ and $V'(\calO)$ respectively, a binary 
relation $Z\subset X\times X'$ is said to be \emph{\dio} if it is \dio\ 
as a subset of $(V\times V')(\calO)$. (Here products are in the category of of 
$R$-schemes, i.e.~fibered over $\Spec(R)$). 

In particular, a \emph{map} $f:X\to X'$ is \dio\ if its graph is \dio\ in 
$(V\times V')(\calO)$. Compositions of \dio\ maps are \dio, and images 
(resp.~inverse images) of \dio\ sets by \dio\ maps are again \dio\ sets.

\Subsubsection{Relative \dio\ sets. }\label{RelDio}
If $V$ is an affine $R$-scheme of finite presentation, and $X\subset 
Y$ are subsets of $V(\calO)$, we say that $X$ is \emph{relatively \dio} 
in $Y$ if it is of the form $D\cap Y$, where $D\subset V(\calO)$ is \dio. 
Of course, if $Y$ is \dio, this is equivalent to $X$ being \dio\ in 
$V(\calO)$.

This notion will be convenient in the following situation: if $Z$ is a 
subset of $V(\calO)$, an $n$-ary relation on $Z$ will be called 
relatively \dio\ if it is a relatively \dio\ subset of $Z^n$ (wiewed 
as a subset of $V^n(\calO)$). In particular, we can speak of a relatively 
\dio\ group (or ring) structure on $Z$, even if $Z$ is not known to be \dio.

\Subsubsection{\dio\ structures. }\label{DioStruc}
Let $V$ be an $R$-scheme of finite presentation, and $X\subset V(\calO)$ a \dio\ 
set. If $\calL$ is a first order language, a \emph{\dio\ $\calL$-structure\/} 
relative to $R$ and $\calO$, or $(R,\calO)$-\dio\ $\calL$-structure, with underlying 
set $X$ is an $\calL$-structure on $X$ such that all subsets of the various 
product sets $X^r$, and all maps $X^r\to X$, relevant to the structure are \dio. 

As an example, take $\calL=\lr$, the language of rings. Then an $(R,\calO)$-\dio\ 
$\lr$-structure consists of:
\begin{romlist}
	\item a \dio\ set $X$ (w.r.t.~$R$, in some $V(\calO)$);
	\item three \dio\ maps $+,-,.:X\times X\to X$;
	\item two elements $0_{X}$ and $1_{X}$ of $X$, which are \dio\ 
	(i.e.~$\{0_{X}\}$ and $\{1_{X}\}$ are \dio\ in $V(\calO)$).
\end{romlist}
Of course, an $(R,\calO)$-\dio\ \emph{ring\/} is an $(R,\calO)$-\dio\ 
$\lr$-structure which satisfies the axioms of rings (commutative with 
unit), in the obvious sense. Equivalently, it is a \dio\ set $X$ with 
a ring structure such that addition and  multiplication are \dio\ 
maps, and the unit is a \dio\ element (the other conditions easily 
follow from these using the ring axioms).

\begin{subprop}\label{DioEtDec}
	Let $V$ be an affine
	$R$-scheme of finite presentation. Let $\calL$ be any first 
	order language, and $X\subset V(\calO)$ an $(R,\calO)$-\dio\ $\calL$-structure. 
	\begin{romlist}
		\item\label{DioEtDec1} Let $r$ a nonnegative integer, and 
		$Z\subset X^r$ 
		a subset which is positive-existentially definable in $\calL$. Then $Z$ 
		is \dio\ as a subset of $V^r(\calO)$.
		\item\label{DioEtDec2} Assume that the positive-existential theory 
		of $\calO$ in $\lr(R)$ is decidable. Then the positive-existential 
		theory of $X$ in $\calL$ is decidable. 
	\end{romlist}
\end{subprop}
\dem \ref{DioEtDec1} By definition, there is an integer $m$ and a 
quantifier-free 
and negation-free formula $\phi$ in $\calL$, in $r+m$ variables, such that
$$
Z=\{\;(x_{1},\ldots,x_{r})\in 
X^r\,\mid\,\exists(y_{1},\ldots,y_{m})\in 
X^m\hbox{ such that 
}\phi(\soul{x},\soul{y})\hbox{ holds }\}.
$$
This is the image, by the projection $V^m\to V^r$ to the first $r$ 
factors, of the set $Z'\subset X^{r+m}$ defined by $\phi$. It suffices to prove 
that $Z'$ is \dio, which is done by an easy induction on the length of 
$\phi$, using elementary properties of \dio\ sets and maps.
\smallskip

Now \ref{DioEtDec2} is an easy consequence of \ref{DioEtDec1}. 
The assumption means that there is a procedure $P$ to decide, for any 
given 
affine $R$-scheme $W$ of finite presentation, whether the set $W(\calO)$ 
is empty or 
not (indeed, $W(\calO)$ may be described as the set of solutions in some 
$\calO^N$ of a finite  system of polynomial equations with coefficients 
in $R$). 
We now need to find another procedure which does the same for 
subsets of $X^r$ which are $\calL$-positive-existentially definable. But 
by \ref{DioEtDec1} such a set $Z$ is \dio\ in $V^r(\calO)$, hence there 
are affine $R$-schemes $W_{1},\ldots,W_{s}$ of finite presentation 
and 
morphisms $\varphi_{i}:W_{i}\to V^r$ such that $Z=\bigcup_{i=1}^s 
\varphi_{i}(W_{i}(\calO))$. So $Z$ is empty if and only if each 
$W_{i}(\calO)$ is empty, which can be detected by applying $P$.\qed
\medskip

In particular, from \ref{DioEtDec2} and Matijasevich's theorem, we get:
\begin{subcor}\label{DioEtDecBis}
	Let $R$ be a ring and $\calO$ an $R$-algebra. Assume that there 
	exists an $(R,\calO)$-\dio\ ring $\Lambda\subset V(\calO)$, for some $R$-scheme 
	$V$ of finite presentation, such that $\Lambda$ is isomorphic to $\ZZ$ as 
	a ring. 

	Then the positive-existential theory 
	of $\calO$ in $\lr(R)$ is undecidable.\qed
\end{subcor}

\section{The specialisation theorem}\label{SecThSpec}
\Subsection{The specialisation map. }\label{SpecVarAb}
Let $R$ be a discrete valuation ring with fraction field $F$, and put 
$S=\Spec(R)$. Denote by $k$ the residue field of $R$. 
Choose an algebraic closure $\ol{F}$ of $F$, and 
a prime $l\neq{\rm char}(k)$. Fix a valuation $\ol{v}$
of $\ol{F}$ extending the valuation $v$ defined by $R$. The 
residue field of $\ol{v}$ is an algebraic closure of $k$, which 
we 
denote by $\kb$. The corresponding decomposition group and inertia 
group are denoted by $D=D_{\ol{v}}$ and $I=I_{\ol{v}}$ 
respectively. 

If $A$ is an abelian scheme over $S$, then $I$ acts 
\emph{trivially\/} 
on $T_{l}(A_{F})$ (because $A[l^n]$ is finite \'{e}tale over $S$, for 
all 
$n$), and we have an isomorphism
\begin{equation}\label{Isomladique}
	T_{l}(A_{F})\fflis T_{l}(A_{k})
\end{equation}
(which depends on the choice of $\ol{v}$: observe that 
$\kb$ is implicit in $T_{l}(A_{k})$). 

On the other hand, if $B$ is another abelian $S$-scheme, every 
$F$-homomorphism $A_{F}\to B_{F}$ extends (uniquely, of course) to 
an $S$-homomorphism $A\to B$ (\cite{BLR}, 1.2, Proposition 8). So we 
get a `specialisation' homomorphism 
\begin{equation}\label{EqMorphSpec1}
	\Hom_{F}\,(A_{F},B_{F})\ffl \Hom_{k}\,(A_{k},B_{k})
\end{equation}
which is \emph{injective}: this can be deduced from the isomorphism 
on Tate modules defined above, or from \ref{RigOdd1}, or from the `rigidity lemma' of 
\cite{GIT}, Proposition 6.1. Of course, this also applies over 
finite extensions of $F$, so that we have an injective homomorphism
of free finitely generated $\ZZ$-modules
\begin{equation}\label{EqMorphSpec2}
	\Hom_{\ol{F}}\,(A_{\ol{F}},B_{\ol{F}})
	\ffl \Hom_{\kb}\,(A_{\kb},B_{\kb})
\end{equation}
which is compatible with  the action of $D_{\ol{v}}$. One 
recovers (\ref{EqMorphSpec1}) from (\ref{EqMorphSpec2}) by taking 
Galois invariants on the left, and $D_{\ol{v}}$-invariants on the 
right (note that morphisms of abelian varieties are always defined 
over finite separable extensions of the ground field, hence we are 
safe from inseparability problems). 
Changing the choice of $\ol{v}$
does not change the cokernel of (\ref{EqMorphSpec2}), up to an 
isomorphism of abelian groups. Moreover:

\begin{subprop}\label{CokSpec}
	The cokernel of the specialisation map {\rm(\ref{EqMorphSpec1})} has 
	no torsion prime to the characteristic of $k$. Consequently, the 
	same holds for {\rm(\ref{EqMorphSpec2})}.
\end{subprop}
\pf Assume we have homomorphisms $u:A\to B$ and $v_{k}:A_{k}\to 
B_{k}$ such that $n\,v_{k}=u_{k}$, where $n$ is prime to ${\rm 
char\,}(k)$. We need to show that there is a $v:A\to B$ such that 
$nv=u$ (this $v$ will automatically lift $v_{k}$). 
Now, $u$ induces a morphism $u[n]:A[n]\to B[n]$ of finite \'{e}tale 
group schemes. Its kernel must then be open and closed in $A[n]$, but 
the assumption implies that it contains $A_{k}[n]$, so it is equal to 
$A[n]$. 
This means that $u$ factors through $[n]_{A}$, as desired.\qed

\begin{subrem}\label{RemExtVal}
	If $R'$ is a discrete valuation ring dominating $R$, with fraction 
	field $F'$ and residue field $k'$, we obtain a specialisation map 
from 
	$\Hom_{\ol{F'}}\,(A_{\ol{F'}},B_{\ol{F'}})$ to 
	$\Hom_{\ol{k'}}\,(A_{\ol{k'}},B_{\ol{k'}})$, 
	with obvious notations. 
	Due to \ref{HomVarAb}\,\ref{PropHom3}, this map is isomorphic to 
	(\ref{EqMorphSpec2}), in the obvious sense; 
	in particular, they have isomorphic cokernels.
	
	These constructions also apply when $R=F$ and $R'$ is 
	a discrete valuation ring containing $F$: of course, 
	(\ref{EqMorphSpec1}) and (\ref{EqMorphSpec2}) are then just identity 
	maps, and therefore the specialisation map `over $R'$' is an 
	isomorphism. Believe it or not, this trivial remark will be used 
	below.
\end{subrem}

\Subsection{The specialisation theorem: notations. }\label{NotThSpec}
Let $k$ be a field of characteristic $p\geq0$. Let 
$U\subset\PP^1_{k}$ 
be a nonempty open set, and let $A\to U$ and $B\to U$ be two abelian 
schemes over $U$. Let $z$ denote the standard coordinate on 
$\PP^1_{k}$, and put $F=k(z)$ (the function field of $U$). 
For $x\in U(k)$, we have, as special cases of (\ref{EqMorphSpec1}) 
and (\ref{EqMorphSpec2}), specialisation maps 
\begin{eqnarray}
	{\spe}_{x}:\Hom_{{F}}\,(A_{F},B_{{F}}) & 
	\ffl & \Hom_{{x}}\,(A_{{x}},B_{{x}}) \label{EqMorphSpec2,5}\\
	\spb_{x}:\Hom_{\ol{F}}\,(A_{\ol{F}},B_{\ol{F}}) 
	& \ffl & \Hom_{\ol{x}}\,(A_{\ol{x}},B_{\ol{x}})
	\label{EqMorphSpec3}
\end{eqnarray}
where (\ref{EqMorphSpec3}) 
actually depends on some choices (in particular, ${\ol{x}}$ is a 
geometric point above $x$). Of course, if $k'$ is any extension of 
$k$, we can do the same for points $x\in U(k')$, using the abelian 
schemes $A_{U'}$ and $B_{U'}$ over $U':=U\times_{\Spec(k)}\Spec(k')$. 
Now define two `regular sets' in $U(k')$:
\begin{equation}\label{DefExc0}
	\begin{array}{rrl}
		\reg(A,B,k')=\reg(k') & := & 
		\{x\in U(k')\,\mid\, {\rm Coker\,}({\spe}_{x}) 
		\hbox{ is finite}\}\cr
		& = & \{x\in U(k')\,\mid\, \spe_{x}\text{ is almost bijective 
		(\ref{DefAlmostOnto})}\}\cr
		& = & \{x\in U(k')\,\mid\, {\rm rk\,}{\rm 
		Hom}_{F}(A_{F},B_{F})
		= {\rm rk\,}\Hom_{x}\,(A_{x},B_{x})\}
	\end{array}
\end{equation}
(for the second equality we use \ref{CokSpec}). And we have the 
`geometric' version of $\reg$:
\begin{equation}\label{DefExc}
		{\Reg}(A,B,k')={\Reg}(k')  := 
		\{x\in U(k')\,\mid\, {\rm Coker\,}({\spb}_{x}) 
		\hbox{ is finite}\}
\end{equation}
with similar equivalent formulations.

\begin{subrem}\label{RemRegAlgClos} Observe that the definition of 
$\Reg$ involves the algebraic closure of $k(z)$, which is a much 
bigger field than $\kb(z)$. 

	Accordingly, if $k$ is algebraically closed, it is easy to see that 
	$\reg(k)\subset\Reg(k)$ (the point is that, with the above 
	notations, we have $x=\ol{x}$) but in general the inclusion is 
	strict. For instance, let $E$ be a $k$-elliptic curve without 
	complex multiplication, and take for $A$ the constant 
	abelian scheme $E\times\Pu{k}\to\Pu{k}$. Then take for 
	$B$ the quadratic twist of $A$ by the double
	cover of $\Pu{k}$ given by $k(\sqrt{z})$ (this extends to an abelian 
	scheme over $U=\GG_{{\rm m},k}$): it is easy to see that $\Reg(k)=U(k)$ 
	while $\reg(k)=\emptyset$.
\end{subrem}

\begin{subrem}\label{RemInclReg} For general $k$, it is not true 
	that $\reg(k)\subset\Reg(k)$. 
	As an example, take $A=E\times\Pu{k}\to\Pu{k}$ as before, and
	assume there is some $k$-elliptic curve $E'$ which 
	is $\kb$-isomorphic to $E$ but has no nontrivial $k$-morphism to 
	$E$. Now take $B\to U$ such that $0\in U(k)$ and the fibre $B_{0}$ 
	is $k$-isomorphic to $E'$, while the $j$-invariant of $B$ in $k(z)$ 
	is not constant. Then $0$ is in $\reg(k)$ but not in $\Reg(k)$.
\end{subrem}

\begin{thm}\label{ThSpec}
	\begin{romlist}
		\item\label{ThSpec1} If $k'\subset k''$ are extensions of $k$, 
		then ${\Reg}(k')={\Reg}(k'')\cap k'$.
		\item\label{ThSpec2} If $x\in k'$ is 
		transcendental over $k$, then $x\in \reg(k')\cap\Reg(k')$.
		\item\label{ThSpec3} {\rm(R.~Noot \cite{Noot})} If $k$ is 
		finitely generated over the prime field, then $\reg(k)\cap{\Reg}(k)$ 
		contains a Hilbert set in $\PP^1(k)$.
	\end{romlist}	
\end{thm}
\pf parts \ref{ThSpec1} (which we shall not use) and \ref{ThSpec2} follow 
from Remark 
\ref{RemExtVal}: for \ref{ThSpec2}, observe that if $x$ is 
transcendental over $k$, then its local ring in $U_{k'}$ contains $F$.

Let us sketch the proof of \ref{ThSpec3} (the reader can find details 
in \cite{Noot}, \S1, where the context is slightly different: the 
base is 
not necessarily an open subset of $\PP^1$, but one looks at the 
specialisation map for closed points, not just rational points). 

Choose a  prime $l\neq p$, and put 
$H_{l}(A_{F},B_{F}):=\Hom_{\QQ_{l}}\,(V_{l}(A_{F}),V_{l}(B_{F}))$, 
just as in (\ref{EqHomTate2}) of \ref{HomVarAb}. 
For $x\in U(k)$, we have a similarly defined $\QQ_{l}$-vector space 
$H_{l}(A_{x},B_{x})$, which is isomorphic to $H_{l}(A_{F},B_{F})$ via 
the isomorphism (\ref{Isomladique}); the $G_{k}$-action on 
$H_{l}(A_{x},B_{x})$ is compatible with the $G_{F}$-action on 
$H_{l}(A_{F},B_{F})$, via the natural group homomorphisms $G_{k}\surjgauche 
D_{\ol{x}}\hookrightarrow G_{F}$. Let us denote both these 
$\QQ_{l}$-spaces by $H_{l}$.

Now if $k$ (hence also $F$) is finitely 
generated over the prime field, we know by Faltings' theorem 
(\ref{HomVarAb}\,\ref{Faltings}) that the rank of 
$\Hom_{\ol{F}}\,(A_{{F}},B_{{F}})$ (resp. 
of $\Hom_{\ol{k}}\,(A_{{x}},B_{{x}})$) is determined by the 
image of $G_{F}$ (resp. $G_{k}$) in ${\rm 
Aut}_{\QQ_{l}}(H_{l})$. 

In particular, we shall have $x\in{\Reg}(k)$ whenever $D_{\ol{x}}$ 
has the same image in ${\rm Aut}_{\QQ_{l}}(H_{l})$ as $G_{F}$. Moreover, 
for such an $x$ the groups $D_{\ol{x}}$ and $G_{F}$ have the same 
invariants in $H_{l}$, hence we can also conclude that $x\in\reg(k)$.

But since the image of $G_{F}$ is an $l$-adic Lie group (as a closed 
subgroup of ${\rm Aut}_{\QQ_{l}}(H_{l})$), we conclude 
from a result of Serre (\cite{SerreRibet}, or \cite{Serre}, 10.6) 
that for $x$ in a suitable Hilbert set the images of $G_{F}$ and
$D_{\ol{x}}$ are equal.\qed
\bigskip

Let us now prove an effective version of (a weaker form of) 
\ref{ThSpec3}:

\begin{thm}\label{ThSpecEff} With the notations and assumptions of 
	{\rm\ref{NotThSpec}}, assume that $k$ is finitely generated over the prime 
	field, and assume that the rank $r$ of $\Hom_{{F}}\,(A_{F},B_{{F}})$ 
	is known. Then there is an algorithm to construct a polynomial $P\in 
	k[z,u]$, with the following properties:
	\begin{romlist}
		\item\label{ThSpecEff1} as an element of $k(z)[u]$, $P$ is separable 
		with no roots in $k(z)$;
		\item\label{ThSpecEff2} for any $\lambda\in U(k)$, we have the property:
		`if $P(\lambda,u)\in k[u]$ has no root in $k$, then 
		$\lambda\in\reg(k)$'.
	\end{romlist}
\end{thm}
\dem Put $d_{A}=\dim\,(A/U)$, $d_{B}=\dim\,(B/U)$, and 
$U=\Spec\,R$, with $R=k[z,D(z)^{-1}]$ for some $D\in k[z]$.

Fix a prime $l\neq{\rm char\,}(k)$, and consider, for $n\in\NN$,
the $U$-schemes $A[l^n]$ and  $B[l^n]$. These are finite 
\'etale $U$-schemes in $\ZZ/l^n\ZZ$-modules, locally free of 
ranks $2\,d_{A}$ and $2\,d_{B}$, respectively (for the \'etale 
topology on $U$). It follows that
\begin{equation}\label{EqHn}
	H_{n}:=\soul\Hom_{\,U\text{-group schemes}}\,(A[l^n],B[l^n])
\end{equation}
is a similar group scheme, with $\ZZ/l^n\ZZ$-rank $4\,d_{A}\,d_{B}$. 

Note that, for given $n$, equations for $A[l^n]$ and 
$B[l^n]$ can be computed from equations of $A$ and $B$. Thus, one 
can describe $A[l^n]$ (resp.~$B[l^n]$) as the spectrum of a locally free 
(hence, in fact, free) $R$-bialgebra $\Lambda_{n,A}$  
(resp.~$\Lambda_{n,B}$) of rank $l^{2\,n\,d_{A}}$ 
(resp.~$l^{2\,n\,d_{B}}$). The scheme $H_{n}$ is simply the 
$U$-scheme of bialgebra morphisms $\Lambda_{n,B}\to\Lambda_{n,A}$, 
which in turn can be described explicitly from the equations. 

Faltings' theorem (\ref{HomVarAb}\,\ref{Faltings}) over $F$
can be stated as
$$ \Hom_{F}\,(A_{F},B_{F})\otimes\ZZ_{l}\fflis 
\varprojlim_{n\geq1}H_{n}(F).$$
Since $(H_{n}(F))_{n\geq1}$ is an inverse system of finite groups, 
it satisfies the Mittag-Leffler condition. It follows that, for $n$ 
large enough (say, $n\geq n_{0}$), the image of $H_{n}(F)$ 
in $H_{1}(F)$ is equal to 
the image of the projective limit, and hence, by Faltings' theorem, 
to the image of $\Hom_{F}\,(A_{F},B_{F})$. The latter is clearly 
isomorphic to $\Hom_{F}\,(A_{F},B_{F})\otimes_{\ZZ}\ZZ/l\ZZ$, 
which is a $\ZZ/l\ZZ$-vector space of dimension $r$. 

Moreover, since the image of the 
natural map $\rho_{n,F}: H_{n}(F)\to H_{1}(F)$ obviously 
decreases as $n$ grows, we now see that $n_{0}$ is computable: just compute 
the image of $\rho_{n,F}$ for increasing $n$, until it has dimension $r$ 
(the knowledge of $r$ is clearly essential here!).

Consider now $H_{n_{0}}$. As a finite \'etale $U$-scheme, it decomposes 
canonically as a disjoint sum of a trivial covering $H_{n_{0}}^{\rm triv}$ 
and a finite \'etale $U$-scheme $H'_{n_{0}}$ with no section over $U$ 
(or, equivalently, over $\Spec(F)$): thus, 
$H_{n_{0}}^{\rm triv}(F)=H_{n_{0}}(F)$, and in fact 
$H_{n_{0}}^{\rm triv}$ can be identified with $H_{n_{0}}(F)\times U$. Consider 
the natural $U$-morphism
$$\rho_{n_{0}}: H_{n_{0}}\to H_{1}.$$
This is a morphism of \etale\ covers which, by construction,  sends 
$H_{n_{0}}^{\rm triv}$ to 
the trivial subcovering of $H_{1}$ whose generic fibre is the image of 
$H_{n_{0}}(F)$ in $H_{1}(F)$. So this subcovering has degree $l^r$.

Let us now take a look at $H'_{n_{0}}$: this is a  computable, open and closed 
subscheme of $H_{n_{0}}$. There is a dense open subscheme $U_{1}$ of 
$U$ (of the form $U_{1}=\Spec\,R_{1}$, with $R_{1}=k[z,(D\,D_{1})^{-1}]$ 
for some computable $D_{1}\in k[z]$)
such that $H'_{n_{0}}\times_{U}U_{1}$ is isomorphic to $\Spec 
R_{1}[u]/(P_{1})$ for some $P_{1}\in k[z,u]$. 

Now put $P=D_{1}P_{1}$. Clearly $P$ satisfies \ref{ThSpecEff1}, 
since $H'_{n_{0}}$ is \etale\ over $U$ and $H'_{n_{0}}(F)=\emptyset$. 
Let $\lambda\in U(k)$ be such that $P(\lambda,u)$ has no root in $k$. 
Then of course $D_{1}(\lambda)\neq0$, hence $\lambda\in U_{1}(k)$. 
Hence $H'_{n_{0}}(\lambda)$ is the set of roots of 
$P_{1}(\lambda,u)$ in $k$, which is empty  by assumption.
This means that $H_{n_{0}}(\lambda)=H_{n_{0}}^{\rm triv}(\lambda)$. Hence 
the image of $H_{n_{0}}(\lambda)$ in $H_{1}(\lambda)$ has 
cardinality $l^r$. But this image contains ${\rm 
Hom}_{k}(A_{\lambda},B_{\lambda})\otimes_{\ZZ}(\ZZ/l\ZZ)$, hence 
the $\ZZ$-rank of $\Hom_{k}\,(A_{\lambda},B_{\lambda})$ is at most 
$r$, hence equal to $r$.\qed

\begin{subrem}\label{RemDefEff} We do not give here a precise definition 
of `effective'. Observe, however, that the above procedure only uses the 
standard constructions of effective algebraic geometry, as explained 
for instance in \cite{FJ}, Chapter 17. In particular they can be 
carried out effectively if $k$ is `presented' over the prime field (\cite{FJ}, Section 
17.2) and $A_{F}$, $B_{F}$ and their group laws are given by explicit 
equations.
\end{subrem}

\section{The relative Jacobian of a fibered surface}\label{SecJacFib}
\Subsection{Notations. }\label{NotFibSurf} In this section, we consider the 
following situation: $k$ is a field (with algebraic closure $\kb$ and 
separable closure $k^{\rm s}$, as usual), $S$ denotes the $k$-projective 
line $\PP^1_{k}$, with standard coordinate $z$ and generic point $\eta$ 
(thus, $\eta=\mathrm{Spec}\,k(z)$). We denote by
\begin{equation}\label{EqFib}
	\theta:X\ffl S=\Pu{k}
\end{equation}
a $k$-morphism with the following properties:
\begin{romlist} 
	\item\label{HypFib1} $X$ is a projective, smooth, geometrically 
	connected surface over $k$;
	\item\label{HypFib2} $\theta$ is surjective with geometrically 
	connected fibres and 
	smooth generic fibre;
	\item\label{HypFib3} every geometric fibre of $\theta$ has at least 
	one reduced component;
	\item\label{HypFib4} $\theta_{k^{\rm s}}:X_{k^{\rm s}}\to \PP^1_{k^{\rm s}}$  
	has a section.
\end{romlist}
(In fact, condition \ref{HypFib3} is easily seen to be a consequence 
of \ref{HypFib1} and \ref{HypFib4}). 
Observe that all these properties are invariant under ground field 
extension. Also, \ref{HypFib2} implies that 
$\calO_{S}\fflis \theta_{\ast}{\calO_{X}}$ universally. We put
\begin{equation}\label{EqPicFib}
	\Pi:=\pics^{0}_{X/k}\quad\text{and}\quad J:=\pics^{0}_{X/S}.
\end{equation}
Thus, $\Pi$ is an abelian variety over $k$, and $J\to S$ is a smooth 
separated $S$-group scheme with connected fibres, by \ref{Pic}\,\ref{Pic6}.

For each extension $k'$ of $k$ and each point $\lambda\in S(k')$ we 
denote by $X_{\lambda}$ the fibre of $\theta$ at $\lambda$: this is a 
projective connected curve over $k'$, which is smooth except for 
finitely many $\lambda$. The fibre $J_{\lambda}$ of $J$ at $\lambda$ 
is the Jacobian $\jac(X_{\lambda})$ of $X_{\lambda}$. 

We have a canonical morphism of $S$-group schemes
\begin{equation}\label{FixVersPicRel}
	\xi: \Pi_{S}:=S\times_{k}\Pi\ffl J
\end{equation}
such that, for $k'$ and $\lambda$ as above, $\xi_{\lambda}:\Pi_{k'}\to 
J_{\lambda}$ is the `restriction to $X_{\lambda}$' morphism from 
$\pics^{0}_{X_{k'}/k'}$ to $\pics^{0}_{X_{\lambda}/k'}$.

\begin{prop}\label{NeronJacFib} 
	{\rm(`connected N\'{e}ron mapping property') } For every smooth group 
	scheme $G$ over $S$, with \emph{connected fibres}, and every 
	nonempty open subscheme $U$ of $S$, the obvious 
	restriction homomorphism
	$$\Hom_{U}\,(G_{U},J_{U})\ffl\Hom_{\eta}\,(G_{\eta},J_{\eta})$$
	is an isomorphism.
\end{prop}
\pf This follows from \cite{BLR}, 9.5, Theorem 4 (b) which says that $J$ 
is the connected component of the N\'{e}ron model of $J_{\eta}$. 
\qed
\begin{prop}\label{PartieFixe} {\rm(`fixed part property') } 
	The morphism $\xi$ of {\rm(\ref{FixVersPicRel})} is universal for 
	morphisms from constant abelian $S$-schemes to $J$. 
	
	In precise terms, if $A$ is an abelian variety over $k$,
	every $S$-morphism $u:A_{S}\to J$ of $S$-group schemes equals 
	$\xi\circ ({\rm Id}_{S}\times v)$, for a uniquely defined 
	$k$-morphism $v: A\to \Pi$ of abelian varieties.
\end{prop}
\pf By descent it suffices to prove the corresponding universal 
property for the morphism $\xi_{k^{\rm s}}:\Pi_{\Pu{k^{\rm s}}}\to J_{k^{\rm s}}$. 
In other words, we may assume that $k$ is separably closed. 

Now let $A$ be an abelian variety over $k$. We have to prove that the natural map
\begin{equation}\label{EqHomPartieFixe}
	\Hom_{S}\,(A_{S},\Pi_{S})\ffl \Hom_{S}\,(A_{S},J)
\end{equation}
obtained by composition with $\xi$ is an isomorphism.

To prove injectivity, let  us first remark that by \ref{HomVarAb}\,\ref{PropHom3}, 
 every $S$-morphism  $A_{S}\to\Pi_{S}$ is constant (in other 
words, we have 
$\Hom_{k}\,(A,\Pi)\flis\Hom_{S}\,(A_{S},\Pi_{S})$). Hence if such a 
morphism is nonzero, it comes from a $k$-morphism $u:A\to\Pi$ such 
that $y:=u(x)\neq0$ for some $x\in A(k)$. It is enough to prove that the 
image of the constant section $y_{S}$ by $\xi$ is a nonzero section 
of $J$ over $S$. But since $k$ is separably closed, $X(k)$ is 
nonempty, so we can view $y$ as a (nontrivial) invertible sheaf $L_{y}$ on $X$, 
algebraically equivalent to zero. By inspecting the definition of 
$\xi$, one then checks that $\xi(y_{S})$ is just the image in $J(S)$ 
of the class of $L_{y}$ in $\picn_{X/S}(S)=\pic(X)/\theta^*\pic(S)$, 
via the morphism (\ref{EqPicnToPic}). Recall that the latter is 
injective: hence, if $\xi(y_{S})$ were zero, then $L_{y}$ would be the 
pullback of an invertible sheaf on $S$, necessarily of degree zero 
(because $L_{y}$ is algebraically equivalent to zero), hence trivial 
(because $S=\PP^1$), a contradiction. 

(Remark: until this last argument, $S$ could have been any projective 
smooth curve over $k$, not just $\PP^1$. To generalise the above to 
such an $S$, just redefine $\Pi$ to be the cokernel of 
$\theta^*:\jac(S)\to\pics^0_{X/k}$).
\smallskip

For surjectivity, we shall use the assumption 
\ref{HypFib4} of \ref{NotFibSurf}, which means (since $k$ is 
separably closed) that $\theta$ has a section $\varepsilon:S\to X$. 
Hence by \ref{Pic}\,\ref{Pic2}, we can view a morphism $u:A_{S}\to J$ 
as an invertible sheaf on $(A\times_{k}S)\times_{S}X$ satisfying certain 
conditions (in particular, the condition of being trivial on 
$(\{0_{A}\}\times_{k}S)\times_{S}X$ since $u$ is a group morphism). But since 
$A_{S}\times_{S}X=A\times_{k}X$, we end up with
an invertible sheaf on $A\times_{k}X$, trivial on $\{0_{A}\}\times_{k}X$ 
and hence (since $A$ is connected) algebraically equivalent to $0$ 
in every fibre of the projection $A\times_{k}X\to A$. In other words, we 
have found a $k$-morphism $v:A\to \soul{\rm Pic}^0_{X/k}=\Pi$, which sends 
$0$ to $0$, hence is a morphism of abelian varieties; checking that 
$\xi\circ v=u$ is then routine.\qed
\medskip

\Subsection{The specialisation map.}\label{SsecSpecJac}
We now fix a nonempty open subscheme $U\subset S$ over which 
$\theta$ is smooth, and an abelian variety $A$ over $k$. We denote 
by $A_{U}$ the constant $U$-abelian scheme $A\times_{k}U$, and by 
$J_{U}$ the restriction of $J$ above $U$ (which is also an abelian 
scheme since $X_{U}$ is a smooth proper $U$-curve).

For each $\lambda\in U(k)$, we have an inclusion 
$j_{\lambda}:X_{\lambda}\inj X$, whence a $k$-morphism 
$j_{\lambda}^\ast:\Pi\to J_{\lambda}$ of abelian varieties, and a 
group homomorphism
\begin{equation}\label{HomAPicInclCourSur}
	\begin{array}{rcl}
		H(\lambda):\quad
		\Hom_{k}\,(A,\Pi) & \ffl &
		\Hom_{k}\,(A,J_{\lambda})\\
		u & \longmapsto & j_{\lambda}^\ast\circ u.
	\end{array}
\end{equation}
For later use, we also have, for any extension $k'$ of $k$, a 
homomorphism 
\begin{equation}\label{HomAPicInclCourSurk'}
		H(\lambda,k'):\quad
		\Hom_{k'}\,(A,\Pi)  \ffl 
		\Hom_{k'}\,(A,J_{\lambda})
\end{equation}
obtained by base field extension (thus, $H(\lambda)=H(\lambda,k)$).

On the other hand, we have the specialisation map 
\begin{equation}\label{EqMorphSpec3Bis}
	\spe_{\lambda}:\Hom_{\eta}(A_{\eta},J_{\eta})\to 
	\Hom_{k}(A,J_{\lambda})
\end{equation}
defined in (\ref{EqMorphSpec3}) (here, $\eta=\Spec(k(z))$ is the 
common generic point of $S$ and $U$). We shall now connect the maps 
(\ref{HomAPicInclCourSur}) and (\ref{EqMorphSpec3Bis}):
\begin{subthm}\label{ThSpecJac} There is a natural group isomorphism
$$\nu:\Hom_{k}\,(A,\Pi)\flis\Hom_{\eta}(A_{\eta},J_{\eta})$$
with the property that, for every $\lambda\in U(k)$, we have 
$\spe_{\lambda}\circ \nu=H(\lambda)$.
\smallskip

In particular, for each such $\lambda$, we have the 
equivalence:
$$\lambda\in\reg(A_{U},J_{U},k)\quad\iff\quad\text{$H(\lambda)$ is almost bijective.}
$$
\end{subthm}
\dem 
Recall from \ref{NotFibSurf} that $j_{\lambda}^\ast:\Pi\to J_{\lambda}$ 
is the fibre at $\lambda$ of the morphism $\xi_{U}:  
\Pi_{S}:=S\times_{F}\Pi\ffl J$ of (\ref{FixVersPicRel}). This allows 
us to insert both $H(\lambda)$ and $\spe_{\lambda}$ in a commutative 
diagram of groups:
\newcommand{\largeur}{12em}
\begin{equation}\label{GrandDiag}
	\begin{array}{lcc}
		{\vtop to 7pt{}}\Hom_{k}\,(A,\Pi) & 
		\varflman%
		{\largeur}{H(\lambda) \text{\scriptsize\ (see 
		(\ref{HomAPicInclCourSur}))}} 
		& \Hom_{k}\,(A,J_{\lambda})\\
		\qquad\Mapdown{\alpha} & & \\
		\Hom_{S}\,(A_{S},\Pi_{S}) &  & 
		\smash{\rlap{$\bigg\vert$}\hskip.2em\bigg\vert}\\
		\qquad\Mapdown{\beta\hbox{\scriptsize\ (composition with $\xi$)}} 
		& & \\
		{\vtop to 7pt{}}\Hom_{S}\,(A_{S},J) & 
		\varflman{\largeur}{\hbox{\scriptsize restriction to fibre at 
		$\lambda$}}
		& \Hom_{k}\,(A,J_{\lambda})\\
		\qquad\Mapdown{\gamma\hbox{\scriptsize\ (restriction to $U$)}}
		& & {\rlap{$\big\vert$}\hskip.2em\big\vert}\\
		{\vtop to 7pt{}}\Hom_{\PP^1}\,(A_{U},J_{U}) & 
		\varflman{\largeur}{\hbox{\scriptsize restriction to fibre at 
		$\lambda$}}
		& \Hom_{k}\,(A,J_{\lambda})\\
		\qquad\Mapdown{\delta\hbox{\scriptsize\ (restriction to $\eta$)}}
		& & {\rlap{$\big\vert$}\hskip.2em\big\vert}\\
		{\vtop to 7pt{}}\Hom_{\eta}\,(A_{\eta},J_{\eta}) & 
		\varflman{\largeur}{\spe_{\lambda}\hbox{\scriptsize\ (see 
		(\ref{EqMorphSpec3Bis}))}} & \Hom_{k}\,(A,J_{\lambda}).\\
	\end{array}
\end{equation}

The commutativity of the diagram is immediate from the definitions 
of $\xi$ and $H(\lambda)$ (for the top square), and from the definition 
of the specialisation map `${\rm sp}$' (for the bottom square). Let 
us now prove that the four left vertical arrows are isomorphisms.

For $\alpha$, this is clear: a morphism of constant abelian schemes 
over a connected $k$-scheme must be constant (this may be seen as a 
special case of \ref{RigOdd2}). For $\beta$, this is the fixed part property 
\ref{PartieFixe}; for both $\gamma$ and $\delta$, this is the N\'{e}ron property 
\ref{NeronJacFib}. This completes the proof.\qed
\Subsection{The geometric specialisation map.}\label{SsecSpecGeomJac}
\Subsubsection{Notations. }\label{NotSpecGeomJac}
Keeping the notations of \ref{SsecSpecJac}, we fix in addition an 
algebraic closure $\kzb$ of $\kb(z)$, and we put 
$\etab:=\Spec(\kb(z))$ and $\etabb:=\Spec(\kzb)$. For $\lambda\in 
U(k)$ (or even in $U(\kb)$), we consider the diagram
\renewcommand{\largeur}{8em}
\begin{equation}\label{GrandDiag2}
	\begin{array}{lcc}
		{\vtop to 7pt{}}\Hom_{\kb}\,(A,\Pi) & 
		\varflman%
		{\largeur}{H(\lambda,\kb)\hbox{\scriptsize\ (\ref{HomAPicInclCourSurk'})}} 
		& \Hom_{\kb}\,(A,J_{\lambda})\\
		\qquad\Mapdown{
\overset{\nu}{\cong}
} & & 
		{\rlap{$\big\vert$}\hskip.2em\big\vert}\\
		{\vtop to 7pt{}}\Hom_{\etab}\,(A_{\etab},J_{\etab}) & 
		\varflman{\largeur}{{\rm sp}_{\lambda}}
		& \Hom_{\kb}\,(A,J_{\lambda})\\
		\qquad\Mapdown{\iota\hbox{\scriptsize\ (field extension } 
		\kb(z)\to\kzb\hbox{\scriptsize)}}
		& & {\rlap{$\big\vert$}\hskip.2em\big\vert}\\
		\Hom_{\etabb}\,(A_{\etabb},J_{\etabb}) & 
		\varflman{\largeur}{\spb_{\lambda}\hbox{\scriptsize\  
		(\ref{EqMorphSpec3})}} & \Hom_{\kb}\,(A,J_{\lambda})\\
	\end{array}
\end{equation}
in which the top square is obtained by applying Theorem 
\ref{ThSpecJac} after replacing $k$ by $\kb$. The commutativity of 
the bottom square is clear from the definitions of both 
specialisation maps.

We cannot expect $\iota$ to be an isomorphism in general. However, 
this is true under additional assumptions on $\theta:X\to S$ which 
will be satisfied in the case we shall consider:
\begin{subthm}\label{deltaIsom} With the notations of \rref{NotFibSurf} 
and \rref{NotSpecGeomJac}, assume in addition that there is a point 
$s\in S(k)$ such that, putting $U=S\setminus\{s\}$:
\begin{romlist}
	\item\label{deltaIsom1} for each point $y$ of $U$, the 
	fibre $X_{y}$ of $\theta$ at $y$ is a \emph{semistable} curve 
	\emph{(i.e.~it has at worst ordinary double points, see \ref{Pic}\,\ref{Pic5})};
	\item\label{deltaIsom2} writing the geometric fibre 
	$(X_{s})_{\kb}$ as $\sum_{i=1}^s m_{i}Y_{i}$ where the $Y_{i}$'s 
	are distinct integral divisors, then none of the multiplicities $m_{i}$ 
	is divisible by $p$; moreover the reduced divisor 
	$(X_{s})_{\rm red}=\sum_{i=1}^s Y_{i}$ is a semistable curve.
\end{romlist}

Then the morphism $\iota$ in diagram {\rm(\ref{GrandDiag2})} is an 
isomorphism. In other words, every $\ol{k(z)}$-morphism from $A$ to 
$J_{\eta}$ is defined over $\kb(z)$.
\smallskip

Consequently, for each $\lambda\in k$, we have the 
equivalence:
$$\lambda\in\Reg(A_{U},J_{U},k)\quad\iff\quad\text{$H(\lambda,\kb)$ is almost bijective.}
$$
\end{subthm}
\dem We may and will assume $k=\kb$. Let us denote by $k(z)^{\rm s}$ 
the separable closure of $k(z)$ in $\kzb$, and by $\eta^{\rm s}$ its 
spectrum. By \ref{HomVarAb}\,\ref{PropHom3}, what we need to show is 
that every $k(z)^{\rm s}$-morphism from $A$ to $J_{\eta}$ is defined 
over $k(z)$. So, we consider the abelian group
$$M:=\Hom_{\eta^{\rm s}}\:(A_{\eta^{\rm s}}, J_{\eta^{\rm s}}).$$
This is a free finitely generated $\ZZ$-module 
(see \ref{HomVarAb}\,\ref{PropHom2}), with a continuous action of 
$G:={\rm Gal\:}({\eta^{\rm s}}/{\eta})$; 
we have to prove that this action is trivial. 

For every $y\in S(k)$, let us denote by $I_{y}\subset G$ one 
of the inertia groups at $y$, and by $P_{y}\subset I_{y}$ the wild 
inertia subgroup (i.e.~the maximal pro-$p$-subgroup of $I_{y}$, or 
the trivial subgroup if $p=0$). It is well known that $G$ is 
generated (as a normal subgroup) by all the $I_{y}$, for $y\in U(k)$, 
together with  $P_{s}$ (if $p=0$ this just means that $U\cong\Aa^1_{k}$ 
is simply connected; if $p>0$, there are nontrivial \'{e}tale coverings 
of $\Aa^1_{k}$ but they are wildly ramified at infinity). Hence, our 
claim will follow from:
\begin{sublem}\label{LemNonRam}
	\begin{romlist}
		\item\label{LemNonRam1} $M$ is unramified over $U$  
		(in other words, for all $y\in U(k)$, the subgroup $I_{y}$
		acts trivially on $M$).
		\item\label{LemNonRam2} $M$ is tamely ramified at $s$ (in other words, 
		$P_{s}$ acts trivially on $M$).
	\end{romlist}	
\end{sublem}
\pf Assertion \ref{LemNonRam2} follows from assumption 
\ref{deltaIsom}\,\ref{deltaIsom2}, and \cite{Saito}, Theorem (3.11) 
(in fact we only use the `easier half' of this result). Let us prove 
\ref{LemNonRam1}.

Let $l\neq p$ be a prime, and consider the $\QQ_{l}$-vector spaces 
$V_{l}(J_{\eta})$ and $V_{l}(A_{\eta})$ (see \ref{TorsVarAb}). Both 
are finite-dimensional $\QQ_{l}$-vector spaces with continuous actions 
of $G$, but in the case of $V_{l}(A_{\eta})$ this action is trivial 
since $A$ is defined over $k$ and $k$ is algebraically closed. 

Consider the natural injective map
\begin{equation}\label{HomTate3}
	M\hookrightarrow \Hom_{\QQ_{l}}\,(V_{l}(A_{\eta}),V_{l}(J_{\eta}))
\end{equation}
which is compatible with Galois actions. Since $M$ is a finitely 
generated $\ZZ$-module, there is a normal
subgroup $G'$ of $G$ of finite index which acts trivially on $M$, 
so that the image of $M$ is in fact contained in 
$\Hom_{\QQ_{l}}\,(V_{l}(A_{\eta}),V_{l}(J_{\eta}))^{G'}$. 
But since the action of $G$ on $V_{l}(A_{\eta})$ is trivial, 
we conclude that 
\begin{equation}\label{HomTate4}
	M\inj\Hom_{\QQ_{l}}\,(V_{l}(A_{\eta}),V_{l}(J_{\eta})^{G'}).
\end{equation}
Now let $y$ be a point of $U(k)$. By assumption 
\ref{deltaIsom}\,\ref{deltaIsom1}, $X_{y}$ is semistable, hence 
$J_{y}$ is semiabelian by \ref{Pic}\,\ref{Pic5}. Hence, 
the action of $I_{y}$ on 
$V_{l}(J_{\eta})$ is \emph{unipotent\/}, by  \cite{sga}, IX, 3.5 (or 
\cite{Saito}, Theorem (3.8)). 
So, we have a finite group (namely $I_{y}/I_{y}\cap G'$) acting unipotently on the 
$\QQ_{l}$-vector space $V_{l}(J_{\eta})^{G'}$: such an action must be 
trivial, so the action of $I_{y}$ on $M$ is trivial too.\qed

\section{Double covers, involutions, and twists}\label{SecRevDouble}
In this section, all rings and schemes will be over $\ZZ[1/2]$ (that 
is, $2$ is invertible in rings, and all residue characteristics in 
schemes will be $\neq2$).

\Subsection{Double covers.}\label{RevDoubles}
By a \emph{double cover} we mean a morphism $\pi:\til{S}\to S$ of 
schemes, which is finite locally free of degree $2$: in other words, 
$\pi$ is affine and the $\calO_{S}$-algebra $\pi_{\ast}\calO_{\til{S}}$ 
is locally free 
of rank $2$ as an $\calO_{S}$-module. Recall that since $\pi$ 
is affine it is completely determined by this $\calO_{S}$-algebra.

For such a morphism $\pi$, there is a canonical 
$\calO_{S}$-linear projector 
$\pi_{\ast}\calO_{\til{S}}\to\calO_{\vphantom{\til{S}}S}$ given by the 
`half-trace'. Consequently, 
$\pi_{\ast}\calO_{\til{S}}\cong\calO_{S}\oplus L$ (as 
$\calO_{S}$-module), where 
$L=\ker(\mathrm{Tr_{\calO_{\til{S}}/\calO_{\vphantom{\til{S}}S}}})$ is an 
invertible $\calO_{S}$-module.

\Subsubsection{Local description. }\label{DoublesLoc}
We can cover $S$ by open affine subsets on which $L$ is trivial. So 
let us assume $S=\Spec(A)$ is affine, and $\til{S}=\Spec(B)$ where 
$B$ has an $A$-basis of the form $(1,\delta)$ with  
$\mathrm{Tr}_{B/A}(\delta)=0$. 
An immediate computation shows that $D:=\delta^{2}\in A$, and that 
$B\cong A[\sqrt{D}]:=A[X]/(X^2-D)$. 


Conversely, every morphism which is locally of the form 
$\Spec(A[\sqrt{D}])\to\Spec(A)$ is obviously a double cover.

\Subsubsection{Globalisation. }\label{DoublesGlob}
Returning to the global isomorphism 
$\pi_{\ast}\calO_{\til{S}}\cong\calO_{S}\oplus L$, we see 
from \ref{DoublesLoc}
that multiplication in $\pi_{\ast}\calO_{\til{S}}$ induces an 
$\calO_{S}$-linear map $\mu:L^{\otimes2}\to\calO_{S}$, 
which of course completely determines the $\calO_{S}$-algebra 
structure on $\pi_{\ast}\calO_{\til{S}}$. 

In other words, the 
category of double covers of $S$ (with $S$-isomorphisms as morphisms) 
is equivalent 
to the category of pairs $(L,\nu)$ where $L$ is an invertible 
$\calO_{S}$-module and $\nu$ is a global section of 
$L^{\otimes-2}$ (with the obvious isomorphisms as morphisms).

The local description also shows that if $\pi:\til{S}\to S$ is a double 
cover, then $\til{S}$ has a \emph{canonical $S$-involution} 
$\sigma_{\pi}=\sigma_{\til{S}/S}$, given locally 
by $\delta\mapsto-\delta$, and globally on 
$\pi_{\ast}\calO_{\til{S}}$ by $x\mapsto x-\mathrm{Tr}(x)$. The 
subring of invariants is $\calO_{S}$; the $\calO_{S}$-submodule of 
$\pi_{\ast}\calO_{\til{S}}$ consisting of anti-invariant elements is 
$L$. 

In particular, the pair $(\til{S}, \sigma_{\pi})$ determines $\pi$; 
geometrically, $\pi$ is the quotient 
(in the category of schemes) of $\til{S}$ by the $\ZZ/2$-action given 
by $\sigma_{\pi}$.

\Subsubsection{Remarks. }\label{RemDoubles}We use the notations of 
\ref{DoublesLoc} and \ref{DoublesGlob}. 
\begin{romlist} 
	\item\label{RemDoubles1} The element $\delta$ used in \ref{DoublesLoc} is well defined up 
	to a unit in $A$, 
	and $D$ is well defined up to the square of a unit; the discriminant 
	of $B$ over $A$, relative to the basis $(1,\delta)$, is $4D$. The 
	scheme of zeros of $D$ in $S$ is, of course, also well defined, and 
	is the \emph{branch locus} of $\pi$.
	%
%
	\item\label{RemDoubles4} The morphism $\pi$ is \'etale if and only if (in the local 
	description) $D$ is invertible in $A$, or (globally) if $\nu$ is a 
	trivialisation of $L^{\otimes-2}$. In this case, locally for the 
	\'etale topology on $S$, $\til{S}$ is  
	isomorphic to the trivial double cover $S\somme S$. 
	\item\label{RemDoubles5} If $\nu=0$ (or, locally, if $D=0$) then 
	$\pi_{\ast}\calO_{\til{S}}$ is the $\calO_{S}$-algebra $\calO_{S}\oplus 
	L$ in which $L$ is an ideal of square zero. In particular, if 
	$L=\calO_{S}$, then $\til{S}$ is the `scheme of dual numbers' 
	over $S$.
	\item\label{RemDoubles6} In the local description, assume that $A$ is a field (of 
	characteristic $\neq2$, of course). Then:
	\begin{itemize}
		\item if $D\neq0$, then $B$ is either $A\times A$ or a quadratic 
		extension of $A$;
		\item if $D=0$, then $B\cong A[X]/(X^2)$. 
	\end{itemize}
	\item\label{RemDoubles7} It is well known that $\til{S}$ is a 
	regular scheme if and only if $S$ is regular and the branch locus 
	of $\pi$ is a regular divisor in $S$. For instance, in the local 
	description, if $A$ is a discrete valuation ring, then $A[\sqrt{D}]$ 
	is regular if and only if $D$ has valuation $0$ or $1$.
	
	In fact, the case of interest 
	to us will be when $S$ and $\til{S}$ are Dedekind schemes (mostly, 
	smooth curves over a field or localisations of such curves). 
\end{romlist}

\Subsection{Weil restriction. }\label{SecWeil}
Let us fix a scheme $S$ and a double cover $\pi:\til{S}\to S$. We 
denote by $\sigma=\sigma_{\pi}$ the canonical involution of $\til{S}$.
If $T$ is an $S$-scheme, we put $\til{T}:=\til{S}\times_{S}T$: this 
is a double cover of $T$, with canonical involution denoted by 
$\sigma_{T}$, or simply by $\sigma$.

For simplicity, we shall restrict our constructions (Weil 
restrictions and twists) to \emph{quasiprojective} $S$-schemes $X\to 
S$: this means, by definition, that $S$ can be covered by open 
subsets $U$ such that the restriction $X_{U}\to U$ is a (locally 
closed) subscheme of a projective $U$-space $\PP^n_{U}$.

If $X$ is a quasiprojective $S$-scheme, we shall denote by $\WR(X)$ 
the functor from $S$-schemes to sets defined by
\begin{equation}\label{EqDefWeil}
	T\longmapsto \WR(X)(T):=\mor_{S}\,(\til{T},X).
\end{equation}
This is known as the \emph{Weil restriction} of the $\til{S}$-scheme $\til{X}$, 
with respect to $\pi$ (in general we consider $\pi$ as fixed once and for all; if 
necessary we shall use the notation $\WR(X\to S,\pi)$). 

By functoriality, we have a canonical involution on $\WR(X)$, deduced from 
$\sigma$, and also denoted by $\sigma$ if no confusion arises.

We refer to section 7.6 of \cite{BLR} for general properties of the Weil 
restriction. Let us recall some of them:
\begin{romlist}
	\item\label{PropWeil1} $\WR(X)$ is representable by a 
	quasiprojective $S$-scheme (which, as is customary, we shall still denote 
	by $\WR(X)$). 
	\item\label{PropWeil2} $\WR$ `commutes with base change': if $S'\to S$ is an 
	$S$-scheme, and we denote by primes the objects over $S'$ obtained by base 
	change, 
	then $\WR(X'\to S',\pi')$ is canonically isomorphic to $\WR(X\to S,\pi)'$.
	\item\label{PropWeil2,5} The functor $X\mapsto \WR(X)$ commutes with 
	finite products (of $S$-schemes, i.e.~fibered over $S$), and in particular takes 
	$S$-group schemes to $S$-group schemes. 
	\item\label{PropWeil3} If $X$ is the affine $n$-space $\Aa^n_{S}$, 
	\emph{and} if 
	$\pi_{\ast}\calO_{\til{S}}$ is a free $\calO_{S}$-module (which is always 
	the case locally on $S$), then $\WR(X)\cong\Aa^{2n}_{S}$. 
	
	Explicitly, assume for simplicity that $S=\Spec(A)$ is affine, and 
	$\til{S}=\Spec(A[\sqrt{D}])$ for some $D\in A$. Then, for any $A$-algebra 
	$k$, we have a natural bijection
	$$\begin{array}{rcl}
	\Aa^{2n}_{S}(k)=k^n\times k^n & \ffl & (k[\sqrt{D}])^n=\WR(X)(k)\cr
	(\soul{x},\soul{y}) & \longmapsto & \soul{x}+\soul{y}\sqrt{D}.
	\end{array}
	$$
	The natural involution $\sigma$ on $\WR(X)$ is given by 
	$(\soul{x},\soul{y})\mapsto(\soul{x},-\soul{y})$.
	\item\label{PropWeil4} If $Y$ is a closed 
	(open) subscheme of 
	$X$, then $\WR(Y)$ is a closed (open) subscheme of $\WR(X)$.
	\item\label{PropWeil5} If $X$ is affine over $S$, then so is $\WR(X)$.
	\item\label{PropWeil6} If $\pi$ is the trivial double cover $S\somme S\to S$, 
	then $\WR(X)\cong X\times_{S}X$, with the involution exchanging 
	factors.
	\item\label{PropWeil7} If $\pi$ is the standard `zero discriminant' cover, 
	i.e.~$\pi_{\ast}\calO_{\til{S}}=\calO_{S}[x]/(x^2)$, then $\WR(X)$ is 
	the relative tangent bundle $T_{X/S}$ of $X$ over $S$. More 
	generally, if $\pi_{\ast}\calO_{\til{S}}=\calO_{S}\oplus L$, where $L$ 
	is an ideal of square zero, then $\WR(X)$ is 
	the vector bundle $T_{X/S}\otimes L^{-1}$ over $X$. The involution 
	$\sigma$ is the bundle automorphism given by multiplication by $-1$.
\end{romlist}

Note that \ref{PropWeil7} amounts to nothing but the usual scheme-theoretic 
definition of the tangent bundle; in particular, the bundle projection 
$T_{X/S}\to X$ corresponds to the morphism $\WR(X)\to X$ deduced from 
the obvious section `$x=0$' of $\pi$. 

For the unfamiliar reader, let us make this explicit in the affine 
case: so, assume $S=\Spec(A)$ and $X=\Spec(R)$ (where $R$ is an 
$A$-algebra). For any $A$-algebra $k$, a $k$-valued point $\zeta$ of $\WR(X)$ is 
an $A$-algebra morphism $R\to k[\eps]=k\oplus k\,\eps$ (we adopt the 
traditional notation $\eps$ for the class of $x$ modulo $x^{2}$). 
Such a morphism has the form $f\mapsto \varphi(f)+\del(f)\,\eps$, 
where $\varphi: R\to k$ is a morphism of $A$-algebras (thus making $k$ 
into an $R$-module), and $\del: R\to k$ is an $A$-derivation. Now $\varphi$ 
defines a $k$-valued point $z\in X(k)$ (the projection of $\zeta$ on 
$X$) and $\del$ is a tangent vector at $z$, in the usual definition by 
derivations. More precisely, if we view $T_{X/S}$ as $\Spec{\rm 
Sym}_{R}\,\Omega^1_{R/A}$ (where $\Omega^1_{R/A}$ stands as usual for K\"ahler 
differentials), then we get a homomorphism ${\rm 
Sym}_{R}\,\Omega^1_{R/A}\to k$ (that is, a $k$-valued point of 
$T_{X/S}$) sending $f\,{\rm d}g_{1}\otimes\cdots\otimes {\rm 
d}g_{n}$ to $\varphi(f)\,\del(g_{1})\ldots\del(g_{n})$.

\Subsection{Twists. }\label{SecTw}
We keep the notations of \ref{SecWeil}. Let $(X,\tau)$ be a 
quasiprojective $S$-scheme with 
involution. Consider the functor from $S$-schemes 
to sets defined by
\begin{equation}\label{EqDefTwist}
	T\longmapsto X^{(\tau)}(T):=\modd_{S}(\til{T},X),
\end{equation}
the set of morphisms $\til{T}\to X$ compatible with the involutions. 
(Here again we omit $\pi$ from the notation).
This is a subfunctor of $\WR(X)$, which we shall call the \emph{(quadratic) 
twist} of $(X,\tau)$ by $\pi$. 

Let us list some properties of this construction (which are easily 
deduced from the properties of $\WR$ stated in \ref{SecWeil}):

\begin{romlist}
	\item\label{PropTwist1} In addition to the involution $\sigma$, 
	we now have on  $\WR(X)$ an involution deduced from $\tau$ by 
	functoriality (and also denoted by $\tau$), which commutes 
	with $\sigma$. It is clear from (\ref{EqDefTwist}) that 
	$X^{(\tau)}$ can be seen as the subfunctor of $\WR(X)$ defined by 
	`$\sigma=\tau$', or, equivalently, as the fixed point
	subfunctor for $\sigma\tau$. As a consequence, $X^{(\tau)}$ is 
	(representable by) a 
	closed subscheme of $\WR(X)$ (hence is also quasiprojective).
	\item\label{PropTwist2} $X^{(\tau)}$ commutes with any base change $S'\to S$, in a 
	sense analogous to \ref{SecWeil}\,\ref{PropWeil2}.
	\item\label{PropTwist2,5} The functor $X\mapsto X^{(\tau)}$ commutes with 
	finite products. If $X$ is an $S$-group scheme and $\tau$ is an 
	automorphism of $X$, then $X^{(\tau)}$ is an $S$-subgroup scheme of $\WR(X)$.
	\item\label{PropTwist3} Assume that $X$ is the affine $(m+n)$-space 
	$\Aa^{m+n}_{S}$, with coordinates $(\soul{x},\soul{x'})$, and $\tau$ acting 
	by $(\soul{x},\soul{x'})\mapsto(\soul{x},-\soul{x'})$. 
	
	Moreover, assume for simplicity that $S=\Spec(A)$ and 
	$\til{S}=\Spec(A[\sqrt{D}])$, as in \ref{DoublesLoc}.
	
	Then $X^{(\tau)}$ is isomorphic to $\Aa^{m+n}_{S}$. Explicitly, 
	with the notations of 	\ref{SecWeil}\,\ref{PropWeil3}, we have the 
	isomorphism
	$$\begin{array}{rcl}
	\Aa^{m+n}_{S}(k)=k^m\times k^n & \fflis & X^{(\tau)}(k)\subset 
	k[\sqrt{D}]^{m+n}\cr
	(\soul{x},\soul{y}) & \longmapsto & (\soul{x},\soul{y}\sqrt{D}).
	\end{array}
	$$
	\item\label{PropTwist4} If $Y$ is a closed (open) subscheme of 
	$X$, stable by $\tau$, then $Y^{(\tau)}$ is a closed (open) sub\-scheme of 
	$X^{(\tau)}$.
	\item\label{PropTwist5} If $X$ is affine over $S$, then so is $X^{(\tau)}$.
	\item\label{PropTwist6} If $\pi$ is the trivial double cover 
	$S\somme S\to S$, then 
	$X^{(\tau)}\cong X$. Consequently, if $\pi$ is \'{e}tale then 
	$X^{(\tau)}$ is an \'{e}tale twist of $X$, i.e.~locally isomorphic to 
	$X$ for the \'{e}tale topology on $S$; explicitly, $X\times_{S}\til{S}$ 
	is $\til{S}$-isomorphic to $X^{(\tau)}\times_{S}\til{S}$.
	\item\label{PropTwist7}Assume that $\pi$ is the standard zero 
	discriminant cover, as in \ref{SecWeil}\,\ref{PropWeil7}, and let 
	$Y\subset X$ be the closed subscheme of fixed points of $\tau$. Then 
	$X^{(\tau)}$ is isomorphic to the normal bundle $N_{Y/X}$ of $Y$ in $X$. 
	(More generally, if $\pi_{\ast}\calO_{\til{S}}=\calO_{S}\oplus L$, where $L$ 
	is an ideal of square zero, then $X^{(\tau)}$ is 
	the vector bundle $N_{Y/X}\otimes L^{-1}$ over $Y$).
\pauseromlist
Let us explain \ref{PropTwist7}. In this case we know from 
\ref{SecWeil}\,\ref{PropWeil7} that $\WR(X)$ is 
the tangent bundle $T_{X/S}$. On this bundle, $\tau$ acts as the tangent map 
to $\tau:X\to X$, and $\sigma$ acts by multiplication by $-1$ on the 
fibres. Using, \ref{PropTwist1}, this means that $X^{(\tau)}$ sits 
above $Y$, and is in fact the subbundle of $(T_{X/S})_{\vert Y}$ on 
which $\tau$ (which is now a vector bundle endomorphism) acts by $-1$. 
Now since all residue characteristics are $\neq2$, the action on $\tau$ 
on $(T_{X/S})_{\vert Y}$ splits it into its $+1$ and $-1$-subbundles, 
the former being $T_{Y/S}$ and the latter being canonically 
isomorphic to the normal bundle of $Y$ in $X$.
\smallskip

We also have a projective analogue of \ref{PropTwist3} 
in the \'{e}tale case:
\finpauseromlist
\item\label{PropTwist8} Assume that $X$ is the projective $(m+n-1)$-space 
	$\PP^{m+n-1}_{S}$, with homogeneous coordinates written as 
	$(Y_{1}:\ldots:Y_{m}:Y'_{1}:\ldots:Y'_{n})=(\soul{Y}:\soul{Y}')$, 
	and $\tau$ acting 
	by $(\soul{Y}:\soul{Y}')\mapsto(\soul{Y}:-\soul{Y}')$. 
	
	Moreover, assume (for simplicity) that $S=\Spec(A)$ and 
	$\til{S}=\Spec(A[\sqrt{D}])$, as in \ref{DoublesLoc}, and 
	that \emph{$D$ is invertible in $A$}.
	Then $X^{(\tau)}$ is isomorphic to $\PP^{m+n-1}_{S}$ (with similar
	homogeneous coordinates denoted by $(\soul{V}:\soul{V}')$), via the map
	$$\begin{array}{rcl}
	\PP^{m+n-1}_{S} & \fflis & X^{(\tau)}\cr
	(\soul{V}:\soul{V}') & \longmapsto & (\soul{V}:\soul{V}'\,\sqrt{D}).
	\end{array}
	$$
\end{romlist}

\begin{subrem}\label{RemProjDeg}
	In case \ref{PropTwist8}, the situation is more complicated if $D$ 
	is not invertible. Assuming for instance that 
	$D=0$, we can determine $X^{(\tau)}$ by using \ref{PropTwist7}. Now 
	the fixed locus of $\tau$ is the disjoint union of two linear 
	subspaces $F$ and $F'$ of $X$, given respectively by $\soul{Y}=0$ 
	and $\soul{Y}'=0$. Hence $X^{(\tau)}$ is a disjoint union of their 
	normal bundles, isomorphic to $X\setminus F'$ and $X\setminus F$ 
	respectively.
\end{subrem}

\Subsection{The case of elliptic curves.}\label{TwEll}

With $\pi:\til{S}\to S$ as in \ref{RevDoubles}, we consider the elliptic 
curve $E$ in the projective $S$-plane $\PP^2_S$, defined in homogeneous 
coordinates $(X,Y,Z)$ by
\begin{equation}\label{EqEll1}
	\llap{$E:\qquad$} Y^2\,Z=P(X,Z)
\end{equation}
where $P$ is homogeneous of degree $3$ with coefficients in $\calO_S$, 
monic in $X$, with nonvanishing discriminant. As usual, we give $E$ the 
group structure with origin $\omega=(0:1:0)$; in this way, $E$ becomes 
a smooth commutative $S$-group scheme. 

The involution $\tau=[-1]_E$ 
sends $(X:Y:Z)$ to $(X:-Y:Z)$ (or equivalently to $(-X:Y:-Z)$), and the 
subgroup $E[2]$ is defined by 
$Y=0$: this is a finite \'etale group scheme of degree $4$ over $S$, 
the disjoint union of $\omega$ and the closed subscheme $E[2]^*$
of `points of exact order $2$' (see \ref{ParInvol}). Since 
$E[2]^*$ is also the intersection 
of $E$ with the line $Y=0$, its complement in $E$ is affine, and can 
be identified (via the affine coordinates $x=X/Y$, $z=Z/Y$) 
with the closed subscheme 
\begin{equation}\label{EqEaff}
	\llap{$E\aff:\qquad$} z=P(x,z)
\end{equation}
of the affine plane $\Aa^2$. Note that the origin $\omega$ conveniently 
has coordinates 
$(0,0)$ in $E\aff$, and that $E\aff$ is invariant under $\tau$, which sends 
$(x,z)$ to $(-x,-z)$. This affine model will turn out to be much more 
useful to us than the `usual' one (the complement of $\omega$ in $E$).

Our goal is to study the twist $E^{(\tau)}$ of $E$ by $\pi$. For 
simplicity, we shall always assume, as in \ref{DoublesLoc}, that 
$S=\Spec(A)$ and $\til{S}=\Spec(B)$ with $B= A[\sqrt{D}]$, for some 
$D\in A$. We put
\begin{equation}\label{EqLocDisc}
	\begin{array}{rcl}
		S\nd & := & \Spec(A[1/D]) \subset S \cr
		\Delta & := & \Spec(A/DA) \subset S.
	\end{array}
\end{equation}
Thus, $\Delta$ is the `branch locus' of $\pi$, and $S\nd$ is its open 
complement. The subscript `nd' stands for `nondegenerate' and will 
also  be used to denote restriction of $S$-schemes to $S\nd$. 
\smallskip

We can easily give a crude `qualitative' description of $E^{(\tau)}$: 
we already know that it is a smooth quasiprojective $S$-group scheme, 
whose restriction to $S\nd$ is an \'{e}tale twist of $E\nd$ (that is, an 
$S\nd$-elliptic curve locally isomorphic to $E\nd$ for the \'{e}tale 
toplogy). On the other hand, if $s$ is a point of $\Delta$, then the 
fibre $E_{s}^{(\tau)}$ of $E^{(\tau)}$ at $s$ is isomorphic to the 
normal bundle of 
the fixed locus of $\tau$. This fixed locus is $E_{s}[2]$ which is 
\'{e}tale over the residue field $\kappa(s)$, hence $E_{s}^{(\tau)}$ 
is the restriction to $E_{s}[2]$ of the tangent bundle of $E$, which 
is trivial. We conclude that there is a canonical isomorphism
\begin{equation}\label{EqFibreDeg}
	E_{s}^{(\tau)}\cong E_{s}[2]\times\Gak{\kappa(s)}.
\end{equation}
so that, geometrically, $E_{s}^{(\tau)}$ is a disjoint union of four 
affine lines.
\smallskip

Describing $E^{(\tau)}$ as a scheme is harder, especially if one 
wants a description by equations and inequations in some projective 
space (recall that $E^{(\tau)}$ is quasiprojective). What we shall 
do is describe by equations an open subgroup scheme of $E^{(\tau)}$, 
sufficiently big for our needs. Specifically, $E$ has two open 
subschemes whose twists are easy to see, namely, $E\nd$ and $E\aff$. 
So let us twist these first:

\Subsubsection{The twist of $E\nd$ }\label{SsecTwEnd} is isomorphic to the 
elliptic curve $\mathcal{E}\nd\subset\PP^2_{S\nd}$ given in 
homogeneous coordinates $(U,V,W)$ by
\begin{equation}\label{EqTwEllnd}
		 V^2\,W=D\,P(U,W).
\end{equation}
The isomorphism is given by
\begin{equation}\label{EqIsomTwEll2}
	\begin{array}{rcl}
		{\mathcal{E}\nd} & \fflis & E\nd^{(\tau)}\cr
		(U:V:W) & \longmapsto & (X:Y:Z)=(\sqrt{D}\,U:V:\sqrt{D}\,W)
	\end{array}	
\end{equation}
which is just a special case of \ref{SecTw}\,\ref{PropTwist8}, restricted 
to the appropriate curve.

\Subsubsection{The twist of $E\aff$ }\label{SsecTwEaff} is isomorphic to the 
affine curve $\mathcal{E}\aff\subset\Aa^2_{S}$ given in affine 
coordinates $(u,w)$ by
\begin{equation}\label{EqEaffTord}
	w=D\,P(u,w).
\end{equation}
The isomorphism is given by
\begin{equation}\label{IsomEaffTord}
	\begin{array}{rcl}
	\mathcal{E}\aff & \fflis & (E\aff)_{\mathstrut}^{(\tau)}\cr
	(u,w) & \longmapsto & (\sqrt{D}^{\mathstrut}\,u, \sqrt{D}\,w)
\end{array}
\end{equation}
(recall that $\tau$ is induced by multiplication by $-1$ on the plane).

If $s$ is a point of $\Delta$, then by (\ref{EqEaffTord}) the fibre  
$(\mathcal{E}\aff)_{s}$ is the affine line $w=0$. On the other hand, 
by definition of $E\aff$, the only point of order two in $(E\aff)_{s}$ is 
the origin, so the fibre of $(E\aff)^{(\tau)}$ at $s$ is the tangent 
line to $E_{s}$ at the origin; from the description (\ref{EqFibreDeg}) 
we see that this 
is the connected component of the $\kappa(s)$-group scheme 
$(E_{s})^{(\tau)}$. 
\medskip

Recall that, as a smooth $S$-group scheme, $E^{(\tau)}$  has a 
\emph{connected component} 
$E^{(\tau)\circ}$, which is the largest open subgroup scheme of $E^{(\tau)}$ 
with connected fibres. We shall now describe $E^{(\tau)\circ}$, 
first (in \ref{CompConnexeTordue})
as a subgroup scheme of $E^{(\tau)}$, and then (in \ref{PropDescrEzero})  
as an $S$-scheme in its own right (that is, by equations).
\begin{subprop}\label{CompConnexeTordue}
	\begin{romlist}
		\item\label{CompConnexeTordue1} $E^{(\tau)\zo}$ is the open subset of 
		$E^{(\tau)}$ given by
		$$E^{(\tau)\circ}=(E\aff)^{(\tau)}\cup (E\nd)^{(\tau)}.$$
		\item\label{CompConnexeTordue2} The multiplication by two in 
		$E^{(\tau)}$ factors through $E^{(\tau)\zo}$. Equivalently, 
		for any $A$-algebra $k$, we have
		$$2\,E^{(\tau)}(k)\subset E^{(\tau)\zo}(k).$$
		\item\label{CompConnexeTordue3} Let $x\in E^{(\tau)}(A)$ 
		correspond to the odd morphism $\til{x}:\til{S}\to E$. Then 
		the following are equivalent:
		\begin{subromlist}
			\item\label{CompConnexeTordue3a} $x\in E^{(\tau)\zo}(A)$;
			\item\label{CompConnexeTordue3b} for each point $s\in\Delta$, 
			$\til{x}$ sends the only point $\til{s}$ of $\til{S}$ above $s$ to the 
			origin of $E_{s}$;
			\item\label{CompConnexeTordue3c} for $s$ and $\til{s}$ as in 
			\rref{CompConnexeTordue3b}, $\til{x}(\til{s})$ is not a nontrivial 
			$2$-division point of $E_{s}$ (equivalently, $\til{x}(\til{s})\in 
			E\aff$).
		\end{subromlist}
	\end{romlist}
\end{subprop}
\dem \ref{CompConnexeTordue1} Since both sides are open subschemes of 
 $E^{(\tau)}$ we need only show that they coincide set-theoretically, 
 or that they have the same fibre over $S$.  
 If $s$ is a point of $S\nd$, then 
 $E_{s}^{(\tau)}=(E\nd)_{s}^{(\tau)}$ is an elliptic curve, hence 
 connected and equal to  $E_{s}^{(\tau)\zo}$. If $s\in\Delta$, we 
 have seen in \ref{SsecTwEaff} that 
$E^{(\tau)\zo}$ and 
$(E\aff)^{(\tau)}$ have the same fibres at $s$. 
\smallskip

\noindent\ref{CompConnexeTordue2} Again, since $\mathcal{E}\szo$ is an open 
subscheme it suffices to see that the set-theoretic image of 
$[2]_{\mathcal{E}}$ is contained in $\mathcal{E}\szo$. This is clear 
above $S\nd$, and above $\Delta$ it follows from the description 
(\ref{EqFibreDeg}). 
\smallskip

\noindent\ref{CompConnexeTordue3} Clearly, the restriction of $x$ to $S\nd$ 
is in $E^{(\tau)\zo}(S\nd)$ by \ref{CompConnexeTordue1}.
Thus, $x\in E^{(\tau)\zo}(S)$ if and only if, 
for each $s\in\Delta$, we have $x(s)\in E_{s}^{(\tau)\zo}$. So we 
may, by base change, 
assume that $A$ is a field (with spectrum $S=\{s\}$) and $D=0$. But by 
\ref{CompConnexeTordue1} (or by \ref{SsecTwEaff}), $E^{(\tau)\zo}$ is then 
equal to 
$(E\aff)^{(\tau)}$, so this is equivalent to the condition that 
$\til{x}$ factors through $E\aff$, which in turn is equivalent to 
\ref{CompConnexeTordue3c} because $E\aff$ is an open subscheme of $E$.
Hence, we have \ref{CompConnexeTordue3a}$\Iff$\ref{CompConnexeTordue3c}.

Obviously, 
\ref{CompConnexeTordue3b} implies \ref{CompConnexeTordue3c}; 
conversely, since $\til{x}$ is odd the point $\til{x}(\til{s})$ must 
be a $2$-division point of $E$, hence \ref{CompConnexeTordue3c} implies 
\ref{CompConnexeTordue3b} because the only $2$-division point of $E\aff$ 
is the origin.\qed

\Subsubsection{Gluing. }\label{SsecColl} It follows from 
\ref{CompConnexeTordue}\,\ref{CompConnexeTordue1}
and the descriptions \ref{SsecTwEnd} and \ref{SsecTwEaff} that 
$E^{(\tau)\circ}$ can be constructed by gluing $\mathcal{E}\nd$ and 
$\mathcal{E}\aff$ in some way. The gluing is in fact the obvious one: 
the affine equation (\ref{EqEaffTord}) defining $\mathcal{E}\aff$ is the 
affine form  (obtained by putting $u=U/V$, $w=W/V$) of the homogeneous equation 
(\ref{EqTwEllnd}) defining $\mathcal{E}\nd$ (observe, however, that 
$\mathcal{E}\nd$ is restricted to $D\neq0$). So we are led to define 
closed subschemes $\ol{\mathcal{E}\szo}$ and $\mathcal{F}$ of 
$\PP^2_{S}$ (in homogeneous coordinates $U,V,W$) by
\begin{equation}\label{EqEbF}
	\begin{array}{ll}
		\ol{\mathcal{E}\szo}:\qquad & V^2\,W=D\,P(U,W) \cr
		{\mathcal{F}}:\qquad & V=D=0.
	\end{array}
\end{equation}
Clearly, $\mathcal{F}\subset\ol{\mathcal{E}\szo}$. Now define 
$\mathcal{E}\szo\subset\PP^2_{S}$ by
\begin{equation}\label{EqTwEll1}
	\mathcal{E}\szo:=\ol{\mathcal{E}\szo}\setminus\mathcal{F}.
\end{equation}
It is clear from the equations (\ref{EqEbF}) that:
\begin{sitemize}
	\item over $S\nd$, we have 
	$\mathcal{E}\szo\times_{S}S\nd=
	\ol{\mathcal{E}\szo}\times_{S}S\nd={\mathcal{E}}\nd$, and
	\item $\mathcal{E}\szo\cap(V\neq0)=\mathcal{E}\aff$.
\end{sitemize}
Moreover the isomorphisms (\ref{EqIsomTwEll2}) and 
(\ref{IsomEaffTord}), viewed as open immersions from ${\mathcal{E}}\nd$ 
and ${\mathcal{E}}\aff$ into $E^{(\tau)}$, obviously glue together to form a 
morphism
\begin{equation}\label{EqMorphTwEll}
	\begin{array}{rcl}
		j:\quad  \mathcal{E}\szo & \ffl & E^{(\tau)}\cr
		(U:V:W) & \longmapsto & (X:Y:Z)=(\sqrt{D}\,U:V:\sqrt{D}\,W).
	\end{array}
\end{equation}
\begin{subprop}\label{PropDescrEzero}
	The morphism $j$ of {\rm(\ref{EqMorphTwEll})} is an isomorphism of 
	$\mathcal{E}\szo$ with the connected component $E^{(\tau)\zo}$ of 
	$E^{(\tau)}$.
\end{subprop}
\dem it is clear from 
\ref{CompConnexeTordue}\,\ref{CompConnexeTordue1} and the 
construction of $j$ that the image of $j$ is $E^{(\tau)\zo}$. We know 
that $j$ induces an isomorphism of ${\mathcal{E}}\nd$ with 
$E\nd^{(\tau)}$ (namely, (\ref{EqIsomTwEll2})) and an isomorphism of 
${\mathcal{E}}\aff$ with $E\aff^{(\tau)}$ (namely, 
(\ref{IsomEaffTord})). To conclude, we only have to check
that (\ref{EqIsomTwEll2}) maps 
${\mathcal{E}}\nd\cap{\mathcal{E}}\aff$ \emph{onto} 
$E\nd^{(\tau)}\cap E\aff^{(\tau)}$, 
which is immediate if one observes that 
$E\nd^{(\tau)}\cap E\aff^{(\tau)}=(E\nd\cap E\aff)^{(\tau)}$.\qed

\medskip

We now prove a useful property of $E^{(\tau)}$, which 
follows directly from its definition:

\begin{subprop}\label{Neron} \emph{(`N\'eron mapping property')} Assume 
that $S$ is integral with generic point $\eta$, and  that $\til{S}$ is regular. 
Then the natural map
$$E^{(\tau)}(S)\ffl E^{(\tau)}(\eta)$$
is an isomorphism.
\end{subprop}
\dem We have to show that the natural map
$$\modd(\til{S},E)\ffl\modd(\til{\eta},E)$$
is bijective. Injectivity is clear because $\til{\eta}$ is dense in 
$\til{S}$ (indeed, $\til{S}$ is flat over $S$ and $\eta$ is dense in 
$S$). Next, take any odd $S$-morphism $x:\til{\eta}\to E$. It suffices 
to extend $x$ to an $S$-morphism $\til{S}\to E$ 
(automatically odd by density). Now, $x$ certainly extends to an 
$S$-morphism $x_{1}:U\to E$ where $U$ is a dense open subset of $\til{S}$. 
We can view $x_{1}$ as a section over $U$ of the $\til{S}$-elliptic 
curve $E\times_{S}\til{S}$. But since $\til{S}$ is regular, the fact 
that $x_{1}$ extends to a section over $\til{S}$ is a standard 
property of elliptic curves (and, more generally, abelian schemes): 
see (\cite{BLR}, 1.2, Proposition 8).

(Remark: we shall only use \ref{Neron} when 
$\dim S=1$. In this case, it is easy to replace the reference to \cite{BLR} by 
the valuative criterion of properness).\qed

\part{Proof of the main theorem}
\section{Twisted elliptic curves over function fields}\label{SecMD} 
\Subsection{Notations.}\label{NotSecMD}
In this section we apply the constructions of Section \ref{SecRevDouble} 
to the following situation:
\begin{itemize}
	\item $k$ is a field of characteristic different from $2$;
	\item the base scheme $S$ is $\PP^1_{k}$ (with standard coordinate 
	$t$);
	\item the double cover $\pi: \til{S}\to S$ is 
\begin{equation}\label{EqPiBis}
	\pi:\Gamma\ffl\PP^1_{k}
\end{equation}
	as in \ref{NotGamma}; thus $\Gamma$ is a smooth 
	curve over $k$, and $\pi$ is \'{e}tale above $\infty$ and ramified 
	at $0$. The unique point above $0$ is denoted by 
	$0_{\Gamma}$, and the natural involution on 
	$\Gamma$ by $\sigma$. 
	
	(Note that $S$ is not affine, while results of Section 
	\ref{SecRevDouble} are often presented in the affine case; however 
	the extension to this case is usually obvious). 
	
	The function field of $\Gamma$ is 
	isomorphic to $k(t,\sqrt{R(t)})$. It will be convenient (and 
	compatible with the notations of \ref{RevDoubles})
	to denote by $\delta$ one 
	of the square roots of $R(t)$ in this field; the affine curve 
	$\pi^{-1}(\Aa^1_{k})$ is then 
	$\Spec k[t,\delta]$, where $k[t,\delta]\cong k[t][s]/(s^2-R(t))$, with $\delta$ 
	corresponding to the class of $s$. 
	
	\item $E$ denotes an elliptic curve over $k$, with  
	equation
	\begin{equation}\label{EqEll2}
	\llap{$E:\qquad$} Y^2\,Z=P(X,Z)
\end{equation}
in projective coordinates $X, Y,Z$, and $E_{S}\to \PP^1_{k}$ denotes the constant 
$\PP^1_{k}$-elliptic curve $E\times_{k}\PP^1_{k}\to \PP^1_{k}$, with its natural 
involution $\tau=[-1]_{E_{S}}$.
\end{itemize}
From these data, we deduce a $\Pu{k}$-group scheme
\begin{equation}\label{EqTwellP1}
	\calE=E_{S}^{(\tau)}\ffl\Pu{k},
\end{equation}
the twist of $E_{S}$ by $\pi$, defined in \ref{SecTw} and studied in \ref{TwEll}. 
This is a smooth commutative $\Pu{k}$-group scheme of relative 
dimension $1$; over the generic point of $\Pu{k}$ it can be given by 
the homogeneous equation (in $\PP^2_{k(t)}$, with projective 
coordinates $U,V,W$)
\begin{equation}\label{EqTwEllnd2}
		 V^2\,W=R(t)\,P(U,W)
\end{equation}
as in (\ref{EqTwEllnd}); this description is in fact valid above 
$\Aa^1_{k}\cap S\nd$ where $S\nd$ (notation of (\ref{EqLocDisc})) 
is the complement of the zeros of $R$ in $S$.

\Subsection{Points of $\calE$.}
Let $T\to S=\Pu{k}$ be any $S$-scheme. By definition of a twist (see 
(\ref{EqDefTwist})), the 
group $\calE(T)=\mor_{S}\,(T,\calE)$ can be described as 
$\modd_{S}(\til{T},E_{S})$, where $\til{T}=\til{S}\times_{S}T$ is the 
double cover of $T$ deduced from $\pi$.

But here, by definition of $E_{S}$, an $S$-morphism $\til{T}\to 
E_{S}=E\times_{k}S$ is the 
same thing as a $k$-morphism $\til{T}\to E$, so we get a canonical
isomorphism
\begin{equation}\label{PtsTwEllP1}
	\calE(T) \fflis {E(\til{T})}\odd :=\modd_{k}(\til{T},E).
\end{equation}
Explicitly, assume that, say, $T=\Spec(L)$ where $L$ is an extension of 
$k(t)$: we can describe an element of $\calE(T)$ as a nontrivial 
solution $(U,V,W)\in L^3$ of (\ref{EqTwEllnd2}). The 
isomorphism (\ref{PtsTwEllP1}) maps this element to 
$(X:Y:Z)=(\delta\,U:V:\delta\,W)$: this is indeed an $L[\delta]$-valued 
point of $E$, which is odd since changing $\delta$ to 
$-\delta$ amounts to applying $[-1]_{E}$.

\Subsection{The twisted elliptic curve: points over function fields. }
\label{SsecPtsMDCourbes}
\Subsubsection{The curve $\til{C}_{g}\,$. }
\label{DefCg}
Consider now our $k$-curve $C$ from \ref{NotC}, and assume given a 
nonconstant rational function $g$ on $C$, which we view as a 
$k$-morphism $C\to S$. We denote by $C_{g}$ the $S$-scheme thus 
obtained, and accordingly by $K_{g}$ (as in \ref{NotC}) 
the function field $K$ of $C$ viewed as a finite extension of $k(t)$, via $g$. 

Denote by $\til{C}_{g}$ the curve $C\times_{g,\PP^1_{k},\pi}\Gamma$ 
(recall that our plan is to let $g$ vary). Thus, we 
have a Cartesian diagram of projective $k$-curves
\begin{equation}\label{DiagCart}
	\begin{array}{ccc}
		\widetilde{C}_{g} & \varflman{1.5cm}{\varphi} & \Gamma \cr
		\mapdown{\pi'} &  & \mapdown{\pi} \cr
		\llap{$C=$ }C_{g}\ & \varflman{1.5cm}{g} &  S\rlap{ $=\PP^1_{k}$}
	\end{array}
\end{equation}
\noindent with all maps finite and flat. The ring of rational 
functions on $\widetilde{C}_{g}$ (its function field, if it is 
irreducible) is $\widetilde{K}_{g}=K[\delta]$.  

In any case, $\pi'$ is a double cover and $\widetilde{C}_{g}$ 
has a canonical involution, which we denote by $\til{\sigma}$ (of 
course, this $\til{\sigma}$ induces the involution on 
$\widetilde{K}_{g}$ sending $\delta$ to $-\delta$).

As explained in Section \ref{intro}, we are interested in the group 
$\calE(K_{g})$ of $K_{g}$-rational points of $\calE$. 
\begin{subprop}
	\label{PtsMDCourbes} We keep the notations and assumptions of 
	\rref{DefCg}.
	\begin{romlist}
		\item\label{PtsMDCourbes1} $\widetilde{C}_{g}$ is geometrically 
		connected over $k$.
		\item\label{PtsMDCourbes2} If the 
		ramification loci of $g$ and $\pi$ in $\PP^1_{k}$ are disjoint, 
		then 
		$\widetilde{C}_{g}$ is smooth over $k$.
		\item\label{PtsMDCourbes3} If $g$ has only simple ramification, 
		then $\widetilde{C}_{g}$ is semistable (see 
		\rref{Pic}\,\rref{Pic5}).
		\item\label{PtsMDCourbes3,5} If $g$ is admissible in the sense of 
		\rref{DefGood}, then $\widetilde{C}_{g}$ is smooth, and for all 
		(resp.~all but finitely many) $\lambda\in k^\ast$, the curve 
		$\widetilde{C}_{\lambda g}$ is semistable (resp.~smooth).
		\item\label{PtsMDCourbes3,75} There is a canonical isomorphism of 
		groups
		\begin{equation}\label{PtsTwEllC}
			\mor_{S}\,(C_{g},\calE) \fflis \modd_{k}(\til{C}_{g},E).
		\end{equation}
		\item\label{PtsMDCourbes4} If $\widetilde{C}_{g}$ is smooth, the 
		natural map $\mor_{S}\,(C_{g},\calE)\to\calE(K_{g})$ is an 
		isomorphism. 
		
		In particular, by {\rm(\ref{PtsTwEllC})}, we have a canonical 
		isomorphism
		\begin{equation}\label{PtsTwEllK}
			\calE(K_{g})\cong\modd_{k}(\widetilde{C}_{g},E).
		\end{equation}
	\end{romlist}
\end{subprop}
\pf \ref{PtsMDCourbes1} is clear because $C$ is geometrically 
connected and $\pi'$ is a \emph{ramified\/} double cover since $\pi$ 
is.
\smallskip

\noindent\ref{PtsMDCourbes2} Let $q$ be a point of 
$\widetilde{C}_{g}$. The assumption implies that either $g$ is 
\'{e}tale 
at $\pi'(q)$, or $\pi$ is \'{e}tale at $\varphi(q)$. In the first 
(resp.~second) case, $\varphi$ (resp. $\pi'$) is \'{e}tale at $q$ by 
base 
change, hence  $\widetilde{C}_{g}$ is smooth at $q$ because $\Gamma$ 
(resp.~$C$) is smooth.
\smallskip

\noindent\ref{PtsMDCourbes3} We may assume $k$ algebraically 
closed. The proof of \ref{PtsMDCourbes2} shows that if $q$ is a 
singular point of $\widetilde{C}_{g}$, then $c=\pi'(q)$ and 
$e=\varphi(q)$ must be branch points of $g$ and $\pi$ respectively, 
mapping to the same point $x$ of $\PP^1_{k}$. If $z$ denotes a local 
coordinate at $x$, the completed local ring of $x$ in $\PP^1_{k}$ is 
isomorphic to $k[[z]]$, and the completed local rings of $c$ and $e$ 
in $C$ and $\Gamma$ are respectively isomorphic, as $k[[z]]$-algebras, to 
$k[[z]][u]/(u^{2}-z)$ and $k[[z]][v]/(v^{2}-z)$. So, the completed 
local ring of $q$ is isomorphic to $k[[z]][u,v]/(u^{2}-z,v^{2}-z)$, 
hence to $k[[u]]/((u^{2}-v^{2})$, which proves that $q$ is an ordinary 
double point.
\smallskip

\noindent\ref{PtsMDCourbes3,5} is an immediate consequence of 
\ref{PtsMDCourbes2}  and \ref{PtsMDCourbes3}, and \ref{PtsMDCourbes3,75} 
is a special case of (\ref{PtsTwEllP1}) (applied with $T=C_{g}$).
\smallskip

\noindent\ref{PtsMDCourbes4} Since $\widetilde{C}_{g}$ 
is a smooth curve and $E$ is projective, every $k$-morphism 
$u:\Spec\widetilde{K}_{g}\to 
E$ extends uniquely to a $k$-morphism $\widetilde{C}_{g}\to 
E$, and the condition on involutions is of course preserved. 
(Alternatively, we can invoke the N\'{e}ron property \ref{Neron}).
\qed
\begin{subrem}\label{RemStrucPtsMDCourbes} It follows from 
(\ref{PtsTwEllK}) and Remark \ref{RemStructMorodd} that if 
$\widetilde{C}_{g}$ is smooth, then $\calE(K_{g})$ is a finitely generated abelian group, 
	with torsion subgroup $E[2](k)$.
\end{subrem}
\medskip

We can now prove assertion \ref{PropGood2} of \ref{PropGood}:

\begin{subcor}\label{DemPropGood2}
	With the same notations and assumptions, let $k'$ be an
	extension of $k$. If $k$ is separably closed in $k'$, then 
	$\calE(k'(t))=\calE(k(t))$ and $\calE(k'(C)_{g})=\calE(K_{g})$. 
\end{subcor}
\dem It suffices to prove the second equality. By 
\ref{PtsMDCourbes}\,\ref{PtsMDCourbes4} this is equivalent to 
$$\modd_{k'}(\til{C}_{g},E)=\modd_{k}(\til{C}_{g},E)$$
which follows from the rigidity property  \ref{RigOdd2} (see 
Remark \ref{RemRigOdd}).\qed

\begin{subrem} From \ref{PtsMDCourbes}\,\ref{PtsMDCourbes4} we obtain a 
	description of the group $\calE(K_{g})$, entirely in terms of 
	(morphisms of) projective curves {\it over $k$}.
	This will be our viewpoint in the subsequent sections, where no 
	`geometry over $K$'  will be involved; in 
	particular we shall then forget about $\calE$ completely.
\end{subrem}
\medskip

Let us now compare $\calE(k(t))$ and $\calE(K_{g})$: 
\begin{prop}\label{Trad1} Assume that $g:C\to \PP^1_{k}$ is 
admissible. Then the following conditions are equivalent:
	\begin{romlist}
		\item\label{Trad11} $g$ is good for $E$ and $\Gamma$, 
		i.e.~$\calE(K_{g})=\calE(k(t))$ (see {\rm\ref{DefGood}});
		\item\label{Trad12} the homomorphism
		\begin{equation}\label{EqMoroddMorodd}
			\begin{array}{rcl}
				\modd_{k}(\Gamma,E) & \ffl & 
				\modd_{k}(\widetilde{C}_{g},E)\cr
				j & \longmapsto & j\circ\varphi
			\end{array}
		\end{equation}
		is almost bijective ({\rm see \ref{DefAlmostOnto}});
		\item\label{Trad13} the homomorphism 
		\begin{equation}\label{InclEnd}
			\modd_{k}(\Gamma,E)/E[2](k)\hookrightarrow 
			\modd_{k}(\widetilde{C}_{g},E)/E[2](k)
		\end{equation}
		(deduced from composition with $\varphi$) is almost bijective. 
	\end{romlist}
\end{prop}
\pf	Of course, the subgroups $E[2](k)$ appearing in (\ref{InclEnd}) are simply 
the groups of {\it constant\/} odd morphisms from $\Gamma$ 
	(resp.~$\widetilde{C}_{g}$) to $E$. Since they are also the torsion 
	subgroups of $\modd_{k}(\Gamma,E)$ and 
	$\modd_{k}(\widetilde{C}_{g},E)$, the equivalence of 
\ref{Trad12} and \ref{Trad13} follows.

	The equivalence \ref{Trad11}$\Leftrightarrow$\ref{Trad12} is clear 
because the maps (\ref{InclPtsMD}) and (\ref{EqMoroddMorodd}) 
correspond to each other via the isomorphism of 
\ref{PtsMDCourbes}\,\ref{PtsMDCourbes4}.\qed
\medskip

\Subsection{Reduction to a problem about Jacobians. 
}\label{ParReduction} 
We shall now reformulate the `goodness' property of \ref{DefGood} (or 
rather, its equivalent form \ref{Trad1}\,\ref{Trad13}) in 
terms of Jacobians. So let $g:C\to \PP^1_{k}$ be an admissible 
morphism. By functoriality of Jacobians, 
$\varphi:\til{C}_{g}\to\Gamma$ gives rise to 
$\varphi^\ast:\jac(\Gamma)\to\jac(\widetilde{C}_{g})$, hence 
to a morphism of abelian groups
\begin{equation}\label{EndToJac}
	\begin{array}{rcl}
		\Hom_{k}^{\rm odd}(E,\jac(\Gamma)) & \ffl & 
		\Hom_{k}^{\rm odd}(E,\jac(\widetilde{C}_{g}))\cr
		j & \longmapsto & j\circ\varphi^\ast
	\end{array}
\end{equation}
where $\varphi^\ast:\jac(\Gamma)\to\jac(\widetilde{C}_{g})$ is deduced from 
$\varphi$, and the `odd' superscripts refer to the involution $[-1]$ 
on $E$, and the involutions induced by the double cover structures 
on $\Gamma$ and $\til{C}_{g}$.

Clearly, the map (\ref{EndToJac}) is connected to (\ref{InclEnd}) via 
the morphism (\ref{MorJacOdd}) of \ref{Jac}: we have a commutative 
diagram 
\begin{equation}\label{DiagMorOddToJac}
	\begin{array}{ccc}
		\modd_{k}(\Gamma,E)/E[2](k) & \fflis & 
		\Hom_{k}^{\rm odd}(E,\jac(\Gamma))\\
		\noalign{\vskip.5ex}
		\Mapdown{\text{(\ref{InclEnd})}} & & \Mapdown{\text{(\ref{EndToJac})}}\\
		\noalign{\vskip.5ex}
		\modd_{k}(\widetilde{C}_{g},E)/E[2](k) & \fflis &
		\Hom_{k}^{\rm odd}(E,\jac(\widetilde{C}_{g}))
	\end{array}
\end{equation}
where all maps are injective, and the horizontal maps are isomorphisms by 
\ref{OddEtJac} (indeed, $\Gamma$ has a rational point fixed by $\sigma$, namely 
$0_{\Gamma}$). So, \ref{Trad1} immediately implies:
\begin{subprop}\label{CritGoodBis} 
	With the notations and assumptions of \rref{Trad1}, $g$ is good if and only 
	if {\rm(\ref{EndToJac})} is almost bijective.\qed
\end{subprop}

Let us now get rid of the annoying $^{\rm odd}$ superscript in 
(\ref{EndToJac}). By its very definition as a fibre product, $\widetilde{C}_{g}$ is 
contained in the product surface $C\times\Gamma$; let us denote by 
$i_{g}=(\pi',\varphi):\widetilde{C}_{g}\inj C\times\Gamma$ the inclusion. 
This induces a morphism of $k$-abelian varieties
\begin{equation}\label{PicInj}
	i_{g}^\ast: \soul{\rm Pic}^0_{C\times\Gamma/k}
	\ffl \jac(\widetilde{C}_{g})
\end{equation}
and, in turn, a group homomorphism
\begin{equation}\label{EndToJacTer}
	\begin{array}{rcl}
		\Hom_{k}\,(E,i_{g}^\ast):  	
		\Hom_{k}\,(E,\soul{\rm Pic}^0_{C\times\Gamma/k})
		& \ffl & 
		\Hom_{k}\,(E,\jac(\widetilde{C}_{g}))\cr
		u & \longmapsto & i_{g}^\ast\circ u.
	\end{array}
\end{equation}
\begin{sublem}\label{InvJacGamma}
	The involution $\sigma$ of $\Gamma$ induces $[-1]$ on $\jac(\Gamma)$.
\end{sublem}
\dem We may assume $k$ algebraically closed. Then a point of $\jac(\Gamma)$ 
corresponds to a divisor $\xi$ of degree $0$. Now $\xi+\sigma(\xi)$ is
the pullback of a divisor of degree $0$ on $\Pu{k}$ (namely 
$\pi_{\ast}(\xi)$). Such a divisor is principal, hence so is $\xi+\sigma(\xi)$.\qed

\begin{subprop}\label{RedPicSurf} 
	Let  $g:C\to \PP^1_{k}$ be an admissible morphism. 
	
	Assume that the homomorphism $\Hom_{k}\,(E,i_{g}^\ast)$ 
	of {\rm(\ref{EndToJacTer})} is almost bijective. Then $g$ is good.
\end{subprop}
\pf Clearly, the canonical involution $\til{\sigma}$ on $\widetilde{C}_{g}$ 
is induced by $\tau:={\rm Id}_{C}\times\sigma$ on $C\times\Gamma$, 
hence $\til{\sigma}^\ast$ on $\jac(\widetilde{C}_{g})$ is compatible (via  
$i_{g}^\ast$) with the involution $\tau^\ast$ 
on $\soul{\rm Pic}^0_{C\times\Gamma/k}$. Now, the assumption of the 
proposition clearly implies (by restriction to the `odd' parts) that 
\begin{equation}\label{EndToJacBis}
	\begin{array}{rcl}
		\Hom_{k}^{\rm odd}(E,\soul{\rm Pic}^0_{C\times\Gamma/k})
		 & \ffl & 
		 \Hom_{k}^{\rm odd}(E,\jac(\widetilde{C}_{g}))\cr
		 u & \longmapsto & i_{g}^\ast\circ u
	\end{array}
\end{equation}
is almost bijective. But we know from \ref{Pic}\,\ref{Pic8} 
that $\soul{\rm Pic}^0_{C\times\Gamma/k}$ is canonically isomorphic to 
$\jac(C)\times\jac(\Gamma)$; under this isomorphism, $\tau^\ast$ 
corresponds (by Lemma \ref{InvJacGamma}) to 
${\rm Id}_{\jac(C)}\times[-1]_{\jac(\Gamma)}$. It follows that the 
\emph{odd} homomorphisms from $E$ to $\soul{\rm Pic}^0_{C\times\Gamma/k}$ 
are those which factor through the inclusion 
$\jac(\Gamma)\varfl{{\rm pr}_{2}^\ast} \soul{\rm Pic}^0_{C\times\Gamma/k}$. 
Hence we have a chain of almost bijective homomorphisms
\begin{equation}\label{EndToJac4}
	\begin{array}{rcccl}
			\Hom_{k}\,(E,\jac(\Gamma)) & \fflis &
		\Hom_{k}^{\rm odd}(E,\soul{\rm Pic}^0_{C\times\Gamma/k})
		 & \sffl{\text{(\ref{EndToJacBis})}} & 
		 \Hom_{k}^{\rm odd}(E,\jac(\widetilde{C}_{g}))\cr
		 u & \longmapsto & {\rm pr}_{2}^\ast\circ u;\hfill 
		 v & \longmapsto & i_{g}^\ast\circ v.
	\end{array}
\end{equation}
Since the composite of these maps is just (\ref{EndToJac}), the 
proposition follows.\qed
\medskip

Of course, we can apply this `over $\kb\,$', to obtain:
\begin{subprop}\label{RedPicSurfBar} 
	Let  $g:C\to \PP^1_{k}$ be an admissible morphism. 
	Assume that the natural homomorphism 
	\begin{equation}\label{EndToJacTerBar}
		\begin{array}{rcl}
			\Hom_{\kb}\,(E,i_{g}^\ast):  	
			\Hom_{\kb}\,(E,\soul{\rm Pic}^0_{C\times\Gamma/k})
			& \ffl & 
			\Hom_{\kb}\,(E,\jac(\widetilde{C}_{g}))\cr
			u & \longmapsto & i_{g}^\ast\circ u.
		\end{array}
	\end{equation}
	is almost bijective. Then $g$ is very good.\qed
\end{subprop}
%
\section{Geometry of a pencil of curves}\label{SecPencil}
\Subsection{Notations. }
In this section we keep $k$, $C$, $Q$, $E$, $\Gamma$, $\pi$ as in \ref{DonnFond}, and 
an admissible $k$-morphism $f: C\to \PP^1_{k}$ as in \ref{Notf}. 
Recall that $0_{\Gamma}\in\Gamma(k)$ denotes the unique zero of $\pi$.
We put $d=\deg(f)$ and
$$c_{0}:=f^{-1}(0)\quad\hbox{and}\quad c_{\infty}:=f^{-1}(\infty).$$
These are viewed interchangeably as closed subschemes or as effective 
divisors on $C$; note that both are reduced of degree $d$, and 
that $c_{0}$ contains the finite set $Q$.

Similarly, we define divisors on $\Gamma$ by
$$\gamma_{0}:=\pi^{-1}(0)\quad\hbox{and}\quad 
\gamma_{\infty}:=\pi^{-1}(\infty).$$
Here, of course, $\gamma_{0}=2[0_{\Gamma}]$ as a divisor, while 
$\gamma_{\infty}$ is a reduced divisor of degree $2$ because of the assumptions on 
$\pi$ in \ref{NotGamma}.

For $\lambda\neq0$ in $k$ (or, more generally, in an extension of 
$k$) we put
$$\X_{\lambda}:=\widetilde{C}_{\lambda f}$$
where $\widetilde{C}_{\lambda f}$ is defined in \ref{DefCg}. Thus, 
$\X_{\lambda}$ is the curve in $C\times\Gamma$ defined by 
`$\lambda f(a)=\pi(b)$' ($a\in C$, $b\in\Gamma$). By abuse, we shall 
still denote by $f$ (resp.~$\pi$) the composed map 
$f\circ{\rm pr}_{1}$ (resp.~$\pi\circ{\rm pr}_{2}$) on $C\times\Gamma$. 

We want to  view $\X_{\lambda}$ as a 
`family of curves with parameter $\lambda$'. This leads to consider 
the rational map
$$\theta_0:=\pi/f:\; C\times\Gamma\;\cdots\ffl\PP^1_{k}.$$
Roughly speaking, $\X_{\lambda}$ is `$\theta_0^{-1}(\lambda)$'. We 
have to make this precise, including when $\lambda$ is $0$ or 
$\infty$; this amounts to `make the rational map $\theta_0$ into a 
morphism'. Now the divisors of zeros and poles of $f$ and $\pi$ on 
$C\times\Gamma$ are 
$$\begin{array}{ll}
	(f)_{0}=c_{0}\times\Gamma, & \quad (f)_{\infty}=c_{\infty}\times\Gamma,\cr
	(\pi)_{0}=C\times \gamma_{0}=2(C\times 0_{\Gamma}), & \quad
	(\pi)_{\infty}=C\times \gamma_{\infty}.
\end{array}
$$
Geometrically (i.e.\ over $\kb$) $(f)_{0}$ and $(f)_{\infty}$ 
are (disjoint) unions of $d$ copies of $\Gamma$; $(\pi)_{\infty}$ is a union of 
two 
copies of $C$, and $(\pi)_{0}$ is a `double $C$'. The function 
$\theta_0$ is undefined at the finite sets 
$S_{0}:=(f)_{0}\cap(\pi)_{0}=
c_{0}\times \gamma_{0}$ and $S_{\infty}:=(f)_{\infty}\cap 
(\pi)_{\infty}= c_{\infty}\times\gamma_{\infty}$. To make $\theta_0$ 
defined everywhere we have to perform some blowups. The situation is quite 
simple at $S_{\infty}$, because $\theta_0=\pi/f$ and $(1/\pi,1/f)$ is a local 
coordinate system at points of $S_{\infty}$. It is more complicated 
at $S_{0}$, where $\pi$ has a double zero: at each point of $S_{0}$ there are local 
coordinates of the form $(u,f)$ and a unit $\varepsilon$ such that 
$\pi=\varepsilon u^2$, so that $\theta_0=\varepsilon u^2/f$. So we need a 
two-step modification of $C\times\Gamma$, detailed below.
%
\Subsection{The blown-up surface $X$. }
\label{Eclatement}
%
%
Let $\rho_{1}:X'\to C\times\Gamma$ be the surface obtained by blowing up 
$S_{0}$ and $S_{\infty}$, viewed as reduced subschemes. 
This gives rise to exceptional divisors 
$D_{1,0}=\rho_{1}^{-1}(S_{0})$ and 
$D_{1,\infty}=\rho_{1}^{-1}(S_{\infty})$. 

If $Y$ is a curve on $C\times\Gamma$, let us denote by 
$\rho_{1}^{\bullet}(Y)$ 
its proper transform on $X'$. Then the divisor of 
$\theta_{1}:=\theta_{0}\circ\rho_{1}$ on $X'$ is given by
\begin{equation*}
	\begin{array}{rrcl}
		\text{poles:}\quad & (\theta_{1})_{\infty} & = & \rho_{1}^{\bullet}(C\times 
		\gamma_{\infty})+\rho_{1}^{\bullet}(c_{0}\times\Gamma),\\
		\text{zeros:}\quad & (\theta_{1})_{0} & = & \rho_{1}^{\bullet}(c_{\infty}\times\Gamma)
		+2\rho_{1}^{\bullet}(C\times 0_{\Gamma})+D_{1,0}.
	\end{array}
\end{equation*}
So $\theta_{1}$ is defined except at 
$S_{0}':=D_{1,0}\cap\rho_{1}^{\bullet}(c_{0}\times\Gamma)$ 
where the zeros 
and poles meet: this finite set maps isomorphically onto 
$S_{0}$ in $C\times\Gamma$. 

Now let $\rho_{2}:X\to X'$ be the blowup of $S_{0}'$, 
with exceptional divisor $D_{2,0}$, and let 
$\rho=\rho_{1}\circ\rho_{2}:X\to C\times\Gamma$. The divisor of 
$\theta:=\theta_{0}\circ\rho$ is given by
\begin{equation*}
	\begin{array}{rrcl}
	\text{poles:}\quad & (\theta)_{\infty} & = & \rho^{\bullet}(C\times 
	\gamma_{\infty})+\rho^{\bullet}(c_{0}\times\Gamma),\\
	\text{zeros:}\quad & (\theta)_{0} & = & \rho^{\bullet}(c_{\infty}\times\Gamma)
	+2\rho^{\bullet}(C\times 0_{\Gamma})+\rho^{\bullet}(D_{1,0}).
	\end{array}
\end{equation*}
Since these divisors have disjoint supports, $\theta$ is a 
morphism to $\PP^1$. Thus, we have a commutative diagram
\begin{center}
\unitlength=1mm
\begin{picture}(34,19)
\put(0,17){\makebox(0,0){$X$}}
\put(17,17){\makebox(0,0){$X'$}}
\put(32,17){\makebox(0,0)[l]{$C\times\Gamma$}}
\put(34,0){\makebox(0,0){$\PP^1_{k}$}}
\put(6,13.5){\vector(2,-1){22}}
\put(20,14){\makebox(0,0){.}}
\put(21,13){\makebox(0,0){.}}
\put(22,12){\makebox(0,0){.}}
\put(23,11){\vector(1,-1){8}}
\put(34,14){\makebox(0,0){.}}
\put(34,13){\makebox(0,0){.}}
\put(34,12){\makebox(0,0){.}}
\put(34,11){\vector(0,-1){8}}
\put(3,17){\vector(1,0){10}}
\put(20,17){\vector(1,0){10}}
\put(25,20){\makebox(0,0){$\rho_{1}$}}
\put(8,20){\makebox(0,0){$\rho_{2}$}}
\put(37,10){\makebox(0,0){$\theta_{0}$}}
\put(28,10){\makebox(0,0){$\theta_{1}$}}
\put(11,8){\makebox(0,0){$\theta$}}
\end{picture}
\end{center}
and the properties of $X$ and the family of curves $\theta$ are 
summarised in the following proposition, where we assume $k$ 
separably closed to alleviate notations; note that all our 
constructions commute with ground field extension. (The picture is 
of course simplified: for instance, $C\times\gamma_{\infty}$ in 
$C\times\Gamma$ consists in two copies of $C$, not one).
\begin{prop}\label{ThFamCourbes}  Assume that $k$ is separably 
closed.
	\begin{romlist}
		\item\label{ThFamCourbes1} $X$ is a smooth projective surface, and 
		$\rho$ is birational.
		\item\label{ThFamCourbes1,5} $\rho$ induces an isomorphism
		\begin{equation}\label{EqPicX}
			\rho^\ast: \soul{\rm Pic}^0_{C\times\Gamma/k}
			\fflis\soul{\rm Pic}^0_{X/k}
		\end{equation}
		of abelian varieties over $k$.
		\item\label{ThFamCourbes2} $\theta$ is projective and flat 
		with geometrically connected fibres. For 
		$\lambda\in k^*$, the fibre $\theta^{-1}(\lambda)$ maps 
		isomorphically via $\rho$ to the curve $\X_{\lambda}$ in 
		$C\times\Gamma$. 
		\item\label{ThFamCourbes3} The fibre $\theta^{-1}(\infty)$ maps 
		isomorphically to $(C\times \gamma_{\infty})\cup 
		(c_{0}\times\Gamma)\subset C\times\Gamma$. It is a union of $d$ 
		disjoint copies of $\Gamma$ and $2$ disjoint copies of $C$, each 
		$\Gamma$ of the first set meeting each $C$ of the second 
		transversally at one point.
		\item\label{ThFamCourbes4} The fibre 
		$\X_{0}=\theta^{-1}(0)$ 
		is a union of:
		\begin{itemize}
			\item a copy of $C$, with multiplicity two in the fibre, 
			mapping isomorphically to the divisor $2(C\times 0_{\Gamma})$ 
			in $C\times\Gamma$;
			\item $d$ disjoint copies of $\Gamma$, each  attached to the above 
			copy of $C$ by identifying $0_{\Gamma}\in\Gamma$ to one of the 
			poles of $f$ on $C$;
			\item $d$ disjoint copies of $\PP^1$, each attached to the 
			above copy of $C$ by identifying one point with a zero of $f$.
		\end{itemize}
		\item\label{ThFamCourbes5} The fibre $X_{\lambda}$ is semistable 
		for each $\lambda\in\PP^1_{k}\setminus\{0\}$.	
		\item\label{ThFamCourbes6} Every component of $D_{2,0}$ 
		(the exceptional divisor of $\rho_{2}$) maps isomorphically to 
		$\PP^1_{k}$ via $\theta$, 
		hence defines a section of $\theta$. The same holds for every
		component of 
		$\rho_{2}^{-1}(D_{1,\infty})=\rho_{2}^{\bullet}(D_{1,\infty})$.
	\end{romlist}
\end{prop}
%
%
\begin{center}
	{\epsfig{file=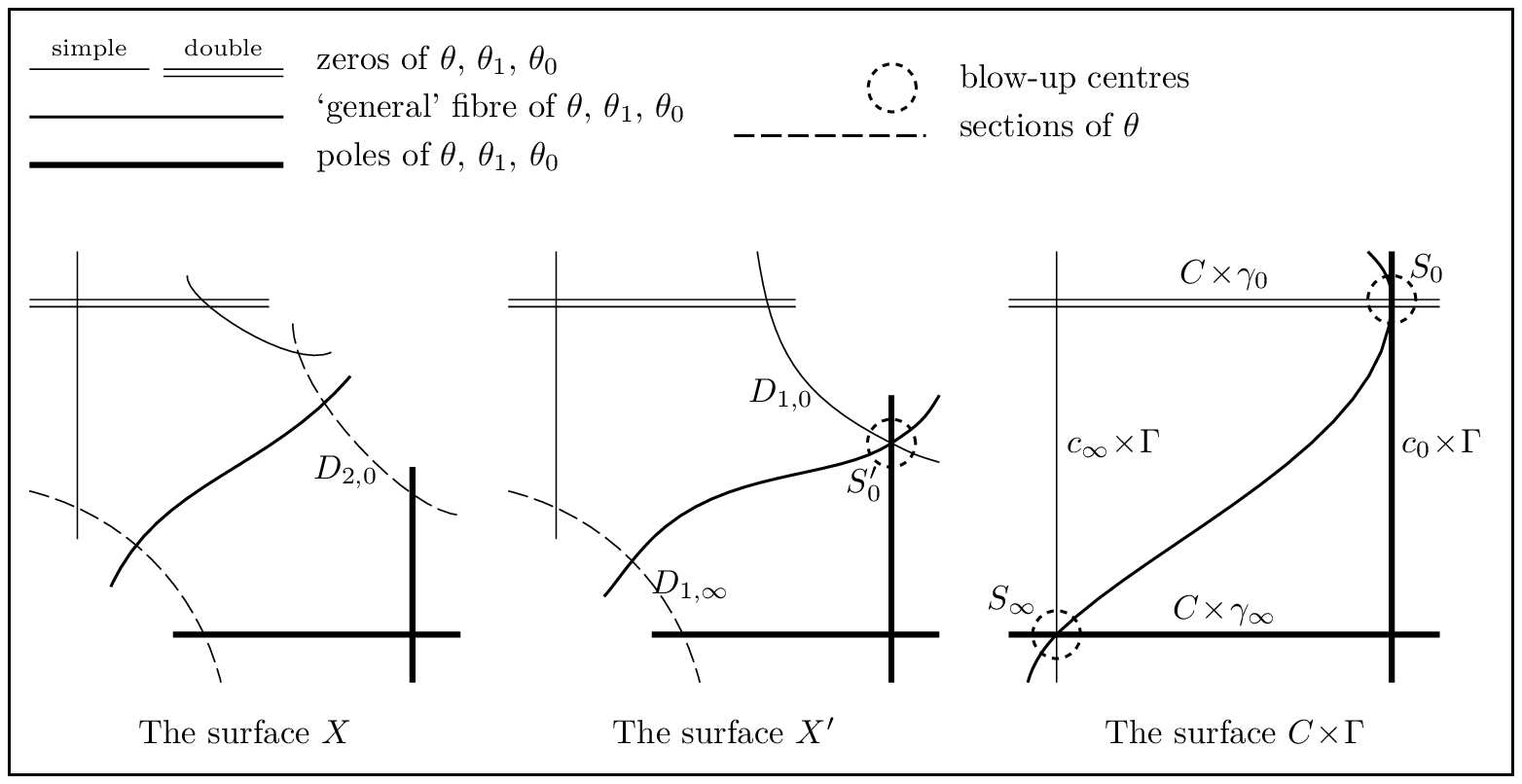, bbllx=69pt, bblly=539pt, bburx=521pt, 
	bbury=773pt}}
\end{center}
\pf Most assertions follow from a careful look at the construction of 
$X$, so we leave the details to the reader. 
For \ref{ThFamCourbes1,5}, use \ref{Pic}\,\ref{Pic7} twice. 
For \ref{ThFamCourbes5}, use \ref{ThFamCourbes3} for the fibre at 
$\infty$, and \ref{PtsMDCourbes}\,\ref{PtsMDCourbes3} for the other 
fibres. (Also, note that if $p>0$ we use the assumption of \ref{NotC} that the 
points of $Q$ are separable over $k$, hence, under our assumptions, 
rational).
\qed
\medskip

Proposition \ref{ThFamCourbes} 
allows us to apply the results of Section \ref{SecJacFib} to $\theta$:

\begin{subprop}\label{RegEtGood} Let $U\subset\Pu{k}$ be the 
smooth locus of $\theta$ (or any nonempty open subset of it).
We have the inclusions
	\begin{romlist}
		\item\label{RegEtGood1} $\reg({E_{U},J_{U},k})\subset\good(E,\Gamma,f,k)$,
		\item\label{RegEtGood2} $\Reg({E_{U},J_{U},k})\subset\vgood(E,\Gamma,f,k)$
	\end{romlist}
where the sets $\good$ and $\vgood$ (resp.~$\reg$ and $\Reg$) are defined 
	in \rref{DefGood} (resp.~{\rm(\ref{DefExc0})} and {\rm(\ref{DefExc})} of 
	Section \rref{SecThSpec}).
\end{subprop}
\dem \ref{RegEtGood1} Let $\lambda$ belong to $\reg({E_{U},J_{U},k})$. 
We have a commutative diagram 
\begin{equation}\label{Diag2Incl}
	\begin{array}{rcl} 
		\theta^{-1}(\lambda) & \sffl{j_{\lambda}} & X \\
		\mapdown{\cong} & & \mapdown{\rho} \\
		X_{\lambda}=\widetilde{C}_{\lambda f} & \sffl{i_{\lambda f}} & 
		C\times\Gamma
	\end{array}
\end{equation}
where $j_{\lambda}$ and $i_{\lambda f}$ are the natural inclusions, 
and the left vertical map is an isomorphism by 
\ref{ThFamCourbes}\,\ref{ThFamCourbes2}. Applying the 
`$\pics^0_{?/k}$' functor, we get a commutative diagram 
\begin{equation}\label{Diag2InclPic}
	\begin{array}{rcl} 
		 J_{\lambda} & \sfflgauche{j_{\lambda}^\ast} & \pics^0_{X/k} \\
		\mapup{\cong} & & \mapup{\overset{\rho^\ast}{\cong}} \\
		\jac(\til{C}_{\lambda f}) & \sfflgauche{i_{\lambda f}^\ast} & 
		\pics^0_{C\times\Gamma/k}
	\end{array}
\end{equation}
where $\rho^\ast$ is an isomorphism by 
\ref{ThFamCourbes}\,\ref{ThFamCourbes1,5}. 

By Proposition \ref{ThFamCourbes}, $\theta$ satisfies all the 
assumptions of \ref{NotFibSurf}. So we can apply Theorem 
\ref{ThSpecJac} (with $A=E$) to conclude that the group homomorphism
$$H(\lambda): \Hom_{k}\,(E,\pics^{0}_{X/k})\to{\rm 
Hom}_{k}\,(E,J_{\lambda}),$$
deduced from $j_{\lambda}^\ast$ by functoriality,
is almost bijective. But from diagram (\ref{Diag2InclPic}) we see 
that the same holds for the homomorphism 
$$\begin{array}{rcl}
		\Hom_{k}\,(E,\soul{\rm Pic}^0_{C\times\Gamma/k})
		& \ffl & 
		\Hom_{k}\,(E,\jac(\til{C}_{\lambda f}))\cr
		u & \longmapsto & i_{\lambda f}^\ast\circ u
\end{array}
$$
which is the morphism (\ref{EndToJacTer}) with $g=\lambda f$. Hence we conclude by 
Proposition \ref{RedPicSurf} that $\lambda f$ is good, 
i.e.~$\lambda\in\good(E,\Gamma,f,k)$.
\smallskip

The proof of \ref{RegEtGood2} is completely similar: just apply 
\ref{RedPicSurfBar} instead of \ref{RedPicSurf}, and
Theorem \ref{deltaIsom} instead of \ref{ThSpecJac}, observing that 
the extra assumptions of \ref{deltaIsom} are satisfied here (in 
particular the fibre at infinity does have a multiple component, but 
its multiplicity is $2$; this of course would ruin our argument in 
characteristic $2$).
\qed

\Subsection{Proofs of Main Theorem \ref{MainTh} and 
Theorem \ref{MainThEff}. }\label{FinDemMainTh}
Our Main Theorem now readily follows from Proposition \ref{RegEtGood} 
and the specialisation theorems \ref{ThSpec} and \ref{ThSpecEff}. Let 
us first prove \ref{MainTh}\,\ref{MainTh2}: let $k'$ be an extension of 
$k$, and let $\lambda\in k'$ be transcendental over $k$. Then 
$\lambda\in\Reg(E_{U},J_{U},k')$ by \ref{ThSpec}\,\ref{ThSpec2}, hence 
$\lambda\in\vgood(k')$ by \ref{RegEtGood}\,\ref{RegEtGood2}.
\smallskip

Assume now that $k$ is finitely generated over the prime field.  By 
\ref{ThSpec}\,\ref{ThSpec3}, there is a Hilbert subset of $k$ contained 
in $\Reg(k)$, hence in $\vgood(k)$ by 
\ref{RegEtGood}\,\ref{RegEtGood2}. This proves \ref{MainTh}\,\ref{MainTh3}.

The proof of Theorem \ref{MainThEff} is similar. Namely, we assume 
here that $k$ is presented over the prime field (\cite{FJ}, 17.2), 
and that $E$ and $C$ are explicitly given. By 
\ref{RegEtGood}\,\ref{RegEtGood1} it suffices to find an effective 
Hilbert set in $\reg({E_{U},J_{U},k})$. This is possible 
by the effective version \ref{ThSpecEff} of the specialisation 
theorem, provided that the rank of $\Hom_{\eta}\,(E_{\eta},J_{\eta})$ 
is known (here $\eta$ is the generic point of $\Pu{k}$) and that we 
have explicit equations for $E_{\eta}$ and $J_{\eta}$. For $E_{\eta}$, 
just take the equations of $E$. For $J_{\eta}$, a procedure for finding 
equations for the Jacobian of a curve is given in \cite{And}. 

For the rank of $\Hom_{\eta}\,(E_{\eta},J_{\eta})$, we have a 
chain of isomorphisms
\setbox0=\hbox{$\;\;{\scriptstyle\rm(\ref{Pic}\,\ref{Pic8})}\;\;\;$}%
\setbox1=\hbox to\wd0{$\;$\rightarrowfill$\;$}
\newcommand{\ffflis}[1]{%
\underset{\rm #1}{\overset{\cong}{\copy1}}}
$$\begin{array}{l}
 \Hom_{\eta}\,(E_{\eta},J_{\eta}) 
 \ffflis{(\ref{ThSpecJac})}
 \Hom_{k}\,(E,\pics\szo_{X/k}) 
\ffflis{(\ref{ThFamCourbes}\,\ref{ThFamCourbes1,5})}
 \Hom_{k}\,(E,\pics\szo_{C\times\Gamma/k})\\
\noalign{\vskip1ex}\ffflis{(\ref{Pic}\,\ref{Pic8})}
 \Hom_{k}\,(E,\jac(C)\times\jac(\Gamma))
 \ffflis{}
 \Hom_{k}\,(E,\jac(C))\times\Hom_{k}\,(E,\jac(\Gamma)) 
\end{array}
$$
which completes the proof since the rank of the right-hand side is 
known by assumption.\qed


\part{Applications to undecidability}\label{Part3}
\section{Self-twisted elliptic curves.}\label{SecSelfCan}

\Subsection{Notations. }\label{SsecNotSelfTw}
We denote by $\kappa$ a field of characteristic different from $2$; 
the important cases in applications will be $\kappa=\QQ$ and 
$\kappa=$ our ground field $k$ of characteristic zero.

We fix an elliptic curve $E$ over $\kappa$ (in our applications, $E$ will 
be defined over $\QQ$). We have a canonical double cover
\begin{equation}\label{EqRevDoubCan}
	\pi=E\to L
\end{equation}
which is the quotient of $E$ by the involution $[-1]$. The curve 
$L$ is smooth projective of genus zero, with a rational point (the 
image of the origin of $E$), hence is isomorphic to $\Pukp$. For the moment 
we refrain from fixing a coordinate on $L$ (or equations of $E$), to 
emphasise the intrinsic character of our constructions. However, we 
can safely denote by $0\in L(\kappa)$ the image of the origin. 

Thus, the branch locus of $\pi$ consists of the point $0$ plus 
three other (geometric) points of $L$, the images of the points of 
order two of $E$.

We denote by $\kappa(L)$ the function field of $L$, by 
$\calO\subset\kappa(L)$ the local ring of $0$ in $L$, and by 
$\mm$ its maximal ideal. Finally we put $S=\Spec(\calO)$. 

\Subsection{The self-twist $\calE$. }\label{SsecSelfTw}
We denote by $E_{L}=E\times_{\kappa}L\to L$ the constant $L$-elliptic 
curve deduced 
from $E$ by base change, and we consider the quadratic twist
\begin{equation}\label{EqTwSurP1}
	\calE\ffl L
\end{equation}
of $E_{L}$ by $\pi$, as defined in \ref{SecTw}: this is the \emph{self-twist} 
of $E$. It is a smooth quasiprojective group scheme over $L$, which 
induces an elliptic curve over the complement $L\nd$ of the branch locus of $\pi$ 
(notation of \ref{TwEll}, (\ref{EqLocDisc})).

Note that if, say, $\xi\in L(\kappa)$ is a rational point, then 
$\pi^{-1}(\xi)$ is a double cover of $\Spec(\kappa)$ (the spectrum of 
a two-dimensional $\kappa$-algebra), and the fibre $\calE_{\xi}$ is 
the twist of $E$ by $\pi^{-1}(\xi)$.

Recall from \ref{TwEll} that we have important open subschemes
$$\calE\aff\subset\calE\szo\subset\calE$$
where $\calE\aff$ is affine over $L$ and $\calE\szo$ is a subgroup 
scheme of $\calE$ with connected fibres.

\Subsection{Sections of $\calE$; the canonical section. 
}\label{SsecCan}%
From the definition (\ref{EqDefTwist}) of a twist, we have in particular a canonical 
isomorphism of groups
\begin{equation}\label{EqSecTwSurP1}
	\mathcal{E}(L)\fflis\modd_{L}(E,E_{L})
\end{equation}
where $E$ is viewed as an $L$-scheme via $\pi$. But by 
definition of $E_{L}$, this boils down to 
\begin{equation}\label{EqSecTwSurP1Bis}
	\begin{array}{rcl}
		\mathcal{E}(L) & \fflis & \modd_{\kappa}(E,E)\\
		 & \fflis & {\rm End}_{\kappa}(E)\times E[2](\kappa)
	\end{array}
\end{equation}
(clearly, odd morphisms $E\to E$ have the form $\tau\circ u$ 
where $u$ is a group scheme endomorphism and $\tau$ is translation 
by a $2$-division point).
In particular, we have a canonical section 
\begin{equation}\label{EqEltCan}
	\gamma:L\to\calE
\end{equation}
corresponding to ${\rm Id}_{E}$ under the first isomorphism  of
(\ref{EqSecTwSurP1Bis}) (and to $({\rm Id}_{E}, 0)$ under the 
second one). We call $\gamma$ the \emph{canonical} element
of $\calE(L)$. It has the following `tautological' description: 
if, say, $\xi\in L(\kappa)$ is a rational point, then we have an 
inclusion $\pi^{-1}(\xi)\inj E$ which obviously respects 
involutions (by definition of $\pi$). But this is precisely the 
definition of a $\kappa$-point of the fibre of $\calE$ at $\xi$, and this 
point is just $\gamma(\xi)$. 

We shall denote by $\gamma_{S}$ (or $\gamma_{\calO}$) the section 
of $\calE_{S}$ induced by $\gamma$; similarly we have 
$\gamma_{\kappa(L)}\in\calE(\kappa(L))$. 
\smallskip

Note that by construction $\gamma$ has infinite order in $\calE(S)$; 
we take this opportunity to prove the following result, which will be 
used in Section \ref{Secpadic}.

\begin{sublem}\label{LemCellRkPos} Assume $k=\QQ$. Then for all but finitely many 
$\xi\in L(\QQ)$, the fibre $\calE_{\xi}$ of $\calE$ at $\xi$ 
is an elliptic curve, and $\gamma(\xi)\in\calE_{\xi}(\QQ)$ has 
infinite order.

In particular, every elliptic curve over 
$\QQ$ has a quadratic twist with positive rank.
\end{sublem}
\dem Clearly, $\calE_{\xi}$ is an elliptic curve for almost every 
$\xi$. For such a $\xi$, let $F_{\xi}\subset\ol{\QQ}$ be the field of 
rationality of the two points $\pm\zeta$ of $\pi^{-1}(\xi)$: 
we have $[F_{\xi}:\QQ]\leq2$. By definition of $\gamma$, 
$\gamma(\xi)$ has finite order in $\calE_{\xi}(\QQ)$ if and only if 
$\zeta$ has finite order in $E(F_{\xi})$. So, all we have to show is 
that the set $T$ of torsion points of $E(\ol{\QQ})$ which are 
rational over some quadratic extension of $\QQ$ is finite. But this 
is an easy consequence of the theory of heights, for which we refer 
to \cite{lang}, Chapter 5 or to  \cite{Serre}, Chapters 2 and 3. 
Namely, the canonical height of any point of $T$ is zero, while 
bounding both the canonical height and the degree defines a finite 
subset of $E(\ol{\QQ})$.\qed
%
%
\begin{subprop}\label{PropEltCan} 
\begin{romlist}
	\item\label{PropEltCan1} The canonical homomorphisms
	$$\calE(L)\ffl\calE(\calO)\ffl\calE(\kappa(L))$$
	are isomorphisms. In particular, by {\rm(\ref{EqSecTwSurP1Bis})}, 
	we have an isomorphism
	\begin{equation}\label{PropEltCanIsom}
		\modd_{\kappa}(E,E)\fflis\calE(\calO).
	\end{equation}
	%
	\item\label{PropEltCan2} The isomorphism {\rm(\ref{PropEltCanIsom})} 
	above maps $\mathrm{End}_{\kappa}(E)$ onto $\calE\szo(\calO)$. In 
	particular, $\gamma_{\calO}\in\calE\szo(\calO)$.
	\item\label{PropEltCan4} The element 
	$\gamma(0)\in\calE\szo_{0}(\kappa)$ is nonzero. 
	\item\label{PropEltCan5} $\calE\aff(\calO)=\calE\szo(\calO)$. (Hence,
	by \rref{PropEltCan2}, $\ZZ\,\gamma_{\calO}\subset\calE\aff(\calO)$).
	\item\label{PropEltCan3} If $E$ does not have complex 
	multiplication over $\kappa$ (that is, $\mathrm{End}_{\kappa}(E)\cong\ZZ$), 
	then $\ZZ\,\gamma_{\calO}=\calE\aff(\calO)=\calE\szo(\calO)$.
\end{romlist}
\end{subprop}
\dem \ref{PropEltCan1} follows from the N\'eron property \ref{Neron}, 
since $E$ is a regular scheme, and \ref{PropEltCan3} is a trivial consequence 
of \ref{PropEltCan1}, \ref{PropEltCan2} and \ref{PropEltCan5}.

Let us prove \ref{PropEltCan2}. Let $u:E\to E$ be an odd $\kappa$-morphism, and 
let $\mu\in\calE(\calO)$ be the corresponding section of $\calE$. 
First, since $\calE\pp_{\kappa(L)}=\calE\szo_{\kappa(L)}$, we have 
$\mu\in\calE\szo(\calO)$ if and only if 
$\mu(0)\in\calE\szo_{0}$, the fibre of $\calE\szo$ at 
zero. Now, $\mu(0)$ is obtained as follows. Consider 
$j:\pi^{-1}(0)\inj E$. This is 
simply the first infinitesimal neighbourhood of the origin $\omega$ of 
$E$, isomorphic as a scheme to $\Spec(\kappa[\varepsilon])$ (with 
$\varepsilon^2=0$). The 
composition $u\circ j:\pi^{-1}(0)\to E$ is an odd morphism, 
hence by definition a point of $\calE_{0}(\kappa)$, which is precisely 
$\mu(0)$. But $u\circ j$ sends the closed point of $\pi^{-1}(0)$ 
to $u(\omega)$; hence, by the criterion 
\ref{CompConnexeTordue}\,\ref{CompConnexeTordue3}\,\ref{CompConnexeTordue3b}, 
$\mu(0)$ is in the connected component if and only if 
$u(\omega)=\omega$, that is, if and only if $u$ is a group 
endomorphism. This proves \ref{PropEltCan2}.

The previous computation, applied with $u={\rm Id}_{E}$ 
(or, equivalently, the tautological desription of $\gamma$) shows 
that $\gamma(0)$ is the point of $\calE\szo_{0}(\kappa)$ 
corresponding to the inclusion $j:\pi^{-1}(0)\inj E$; this 
is clearly nonzero, which proves \ref{PropEltCan4}.

It remains to prove \ref{PropEltCan5}. We know that $\calE\aff$ is an open 
subscheme of $\calE\szo$, and that they have the same fibre at the 
closed point $0$ of $\Spec(\calO)$ (cf.~\ref{SsecTwEaff}). It follows that if
$z:\Spec(\calO)\to\calE\szo$ is a section, then 
$z^{-1}(\calE\aff)$ is an open subscheme of $\Spec(\calO)$ which 
contains the closed point, hence is equal to $\Spec(\calO)$.
\qed

\begin{subrem}\label{RemNeron1}
	Assertion \ref{PropEltCan1} generalises (with the same proof) 
	in the following way: if 
	$\calO\subset\calO'$, where $\calO'$ is a regular semilocal domain of 
	dimension $1$, whose Jacobson radical is generated by the maximal 
	ideal of $\calO$, and 
	if $K$ denotes the fraction field of $\calO'$, then 
	$\calE(\calO')\ffl\calE(K)$ is an isomorphism. 
	
	This applies in particular if $C$ is a smooth curve over (some 
	extension of) $\kappa$, given with a morphism $g:C\to L$, and $\calO'$ 
	is the semilocal ring of $C$ at some set of simple poles of $g$.
\end{subrem}

\begin{subrem}\label{RemEltCan2}
	Assertion \ref{PropEltCan3} has the following consequence. Assume in 
	addition that $E$ does not have complex multiplication over the 
	algebraic closure of $\kappa$, and hence over any extension of 
	$\kappa$. For an extension $\kappa'$ of $\kappa$, let
	$\calO_{\kappa'}$ be the local ring at infinity on $L_{\kappa'}$.
	Then it follows from 
	\ref{PropEltCan3} that $\calE\szo(\calO)\flis\calE\szo(\calO_{\kappa'})$ 
	since both are generated by the same element $\gamma_{\calO}$. In 
	other words, the group $\calE\szo(\calO)$, which is isomorphic to 
	$\ZZ$, is essentially independent of the ground field $\kappa$.
\end{subrem}

\begin{subrem}\label{RemEltCan3}
	It follows from \ref{PropEltCan5} that $\calE\aff(\calO)$ is a 
	subgroup of $\calE(\calO)$, even though $\calE\aff$ is not a 
	subgroup scheme of $\calE$.
\end{subrem}

\section{The ring $\Lambda$ and its multiplication.}\label{SecLambda}
\Subsection{Notations, definition of $\Lambda$. }\label{SsecNotLambda}
In this section we keep the notations and assumptions of Section 
\ref{SecSelfCan} (including the self-twist $\calE$ and its canonical 
section $\gamma$), but now we assume $\car\kappa=0$. 

We fix a ring $\calO'$ 
containing $\calO$, with the aim of proving that the \dio\ theory of $\calO'$ 
(with constants $\calO$) is undecidable. 

We identifiy $\calE(\calO)$ with a subset of $\calE(\calO')$, and 
similarly for $\calE\aff$, $\calE\szo$, etc. In particular we have a subgroup
\begin{equation}\label{EqDefLambda}\Lambda:=\ZZ\,\gamma_{\calO}\subset 
\calE\szo(\calO)\subset \calE\szo(\calO')
\end{equation}
which is, in fact, contained in $\calE\aff(\calO)$ by 
\ref{PropEltCan}\,\ref{PropEltCan5}, and therefore also in 
$\calE\aff(\calO')$. 

The group isomorphism 
$\ZZ\flis\Lambda$ sending $n$ to $n\gamma_{\calO}$ defines a ring 
structure on $\Lambda$. To prove the \dio\ undecidability of $\calO'$, 
it suffices to prove that, for suitable $E$, this ring 
$\Lambda\subset\calE\aff(\calO')$ 
is \dio\ (this makes sense since $\calE_{{\rm aff},\calO}$ is an 
affine $\calO$-scheme of finite presentation; we shall be more explicit 
in \ref{SsecEqES0} below).
\smallskip

Proving that $\Lambda$ is a \dio\ ring involves two tasks:
\begin{sitemize}
	\item show that the ring structure on $\Lambda$ is relatively \dio, 
	in the sense of \ref{RelDio},
	\item show that $\Lambda\subset\calE\aff(\calO')$ is \dio.
\end{sitemize}
Concerning the second property, note that by \ref{PropEltCan}\,\ref{PropEltCan3}, 
it is true for $\calO'=\calO$ 
if $E$ has no complex multiplication, which we shall always assume in 
applications. In fact this extends to the situation mentioned in 
\ref{RemNeron1}. For other rings (specifically for function 
fields of curves) our standard weapon will be Theorem 1.8. 

But this will come later; this section is devoted to the first property, which 
involves checking several points. Here are the 
easy ones:

\begin{subprop}\label{PartieFacile}
	\begin{romlist}
		\item\label{PartieFacile1} The graph of the addition law on $\Lambda$ is 
		relatively \dio\ in $\Lambda^3$, with respect to $\calO$ (here 
		$\Lambda$ is seen as a subset of $\calE\aff(\calO')$).
		\item\label{PartieFacile2} The unit $\{\gamma_{\calO'}\}$ of 
		$\Lambda$ is a \dio\ subset of $\calE\aff(\calO')$.
	\end{romlist}
\end{subprop}
\dem \ref{PartieFacile2} is obvious since $\gamma\in\calE\aff(\calO)$. 

For \ref{PartieFacile1}, we have to be careful because $\calE\aff$ is 
not a group scheme; however, the graph $G$ of `addition' in $\calE_{{\rm 
aff},\calO}^3$ makes sense, as the intersection of $\calE_{{\rm 
aff},\calO}^3$ with the graph of addition on the $\calO$-group scheme 
$\calE_{\calO}^3$. Moreover, $G$ is clearly a closed $\calO$-subscheme of 
$\calE_{\calO}^3$, hence defines a ternary relation on 
$\calE\aff(\calO')$ which is \dio\ with respect to $\calO$. The conclusion 
follows by restriction to $\Lambda$.
\qed
\medskip

\noindent Thus (as always with this method) the hard part is the \dio\ character 
of multiplication in $\Lambda$, which will occupy the rest of this 
section. 

\Subsection{Evaluating at zero. }\label{SsecEvalInfty} 
Recall that by (\ref{EqFibreDeg}) the fibre $\calE\szo_{0}$ of $\calE\szo$ at 
zero is isomorphic to the additive group $\Gak{\kappa}$. Hence 
(once such an isomorphism is fixed, which we assume from now on), 
evaluating sections at $0$ defines a group homomorphism
\begin{equation}\label{EqEvInfini}
	\ev_{0}:\quad\calE\szo(\calO)\ffl\calE\szo(\calO/\mm)\flis\kappa.
\end{equation}
If we restrict this map to $\calE\aff(\calO)$, embedded, 
say, in the affine plane $\Aa^2_{\calO}$, this simply consists in reducing 
coordinates modulo $\mm$, and then applying the isomorphism with 
$\Gak{\kappa}$ (which is algebraic, hence must be given by a 
polynomial in the coordinates, with coefficients in $\kappa$).
\smallskip

The restriction of $\ev_{0}$ to $\Lambda$ necessarily has the form
\begin{equation}\label{EqEvInfiniLambda}
	n\gamma_{\calO}\longmapsto n\,\ev_{0}(\gamma_{\calO})=n\gamma(0).
\end{equation}
Since $\car\kappa=0$ and $\gamma(0)\neq0$ 
by \ref{PropEltCan}\,\ref{PropEltCan4}, this map is injective. 
Therefore we can `encode' the multiplication on $\Lambda$ as follows: 
if $z_{i}=n_{i}\gamma_{\calO}$ ($i=1,2,3$) are three elements of 
$\Lambda$, then:
\begin{equation}\label{EqCodeMult}
	\begin{array}{rcl}
		z_{3}=z_{1}z_{2}\;\text{ (in $\Lambda$)} & \Iff &
		n_{3}=n_{1}n_{2}\;\text{ (in $\ZZ$)}\\
		&\Iff & \ev_{0}(z_{3})\,\ev_{0}(\gamma)=
		\ev_{0}(z_{1})\,\ev_{0}(z_{2})\;
		\text{ (in $\kappa$)}.
	\end{array}
\end{equation}
The last condition involves the relation 
$t_{3}\,\ev_{0}(\gamma)=t_{1}\,t_{2}$ in $\kappa$. This is a 
polynomial relation (in which $\ev_{0}(\gamma)$ is a constant), 
which is good news. But it also involves $\ev_{0}$, i.e.~essentially 
a reduction modulo $\mm$, which is rather bad news. In fact, from now 
on, all the hard work will consist in showing, in various contexts, 
that in some sense reduction modulo $\mm$ has good \dio\ properties.

\Subsection{Explicit equations. }\label{EqExplicites} 
Assume now that $E\subset\PP^2_{\kappa}$ is given by an equation
\begin{equation}\label{EqEllSurk}
	Y^2Z=P(X,Z)=X^3+a\,X^2\,Z+b\,X\,Z^2+c\,Z^3
\end{equation}
in homogeneous coordinates $(X,Y,Z)$ (in our applications, $a,b,c$ will 
be in $\QQ$). We may, and will, identify $\Pukp$ with $L$ via the double cover 
$Z/X:E\to\Pukp$ (also called $\pi$); this is the inverse of the `usual' 
coordinate $x:=X/Z$, for which we shall have little use. We denote by $t$ the 
standard coordinate on $\Pukp$. Thus, the branch locus of $\pi$ consists of 
the three (geometric) 
zeros of $P(1,t)$ and the point $0$. The ring $\calO$ is then $\kappa[t]_{(t)}$, 
with maximal ideal $\mm=t\calO$. 

\begin{subrem}\label{RemChoixOrigine} Unlike $0$, the `point at infinity'
	of $L$ (the pole of $t$ in $L=\Pukp$) has no intrinsic meaning; in fact, 
	by a change of coordinates it can be chosen arbitrarily in 
	$L\setminus\{0\}$. In particular, assume that $E$ is defined over $\QQ$; 
	then, by \ref{LemCellRkPos}, 
	we can choose the equation (\ref{EqEllSurk}) (with $P\in\QQ[X,Z]$) in 
	such a way that the fibre $\calE_{\infty}$ of $\calE$ at $\infty$ is an 
	elliptic curve (this simply means $c\neq0$) and, moreover, that the 
	point $\gamma(\infty)\in\calE_{\infty}(\kappa)$ has \emph{infinite order}.
\end{subrem}

Let us now give equations for (some pieces of) $\calE$.

\Subsubsection{Equations for $\calE_{S}\szo$ and $\calE_{{\rm aff},S}$. }
\label{SsecEqES0}
We only give the results derived from Section \ref{SecRevDouble}, leaving 
details to the reader.

The traditional way of describing $\pi$ as a double cover of $\Pukp$ is by 
`extracting the square root of $P(x,1)$'; however, $P(x,1)=t^{-3}P(1,t)$ does not 
belong to $\calO$ (it has a triple pole at zero), so instead we put 
\begin{equation}\label{EqDInfini}
		\rho  :=  t/P(1,t) = t/(1+a\,t + b\,t^{2} + c\,t^{3});
\end{equation}
this is a uniformising parameter of $\calO$, such that the double cover $\pi$ 
is given above $S$ by $\Spec(\calO[\sqrt{\rho}\,])$. Accordingly, by 
\ref{PropDescrEzero}, $\calE_{S}\szo\to S$ can be described in $\PP^2_{S}$ by
\begin{equation}\label{EqAutoTw}
	\calE_{S}\szo=\ol{\calE\szo}\setminus\mathcal{F},\hbox{ 
	with }
	\left\{\begin{array}{ll}
		\ol{\calE\szo}:\quad & V^2\,W=\rho\,P(U,W)\\
		\mathcal{F}: \quad & V=\rho=0	
	\end{array}\right.
\end{equation}
in projective coordinates $U,V,W$. (In fact, this is not just a description 
of $\calE_{S}\szo$, but also 
of the restriction of $\calE\szo$ above the complement of $\infty$ in 
$\Pukp$.) The unit section of the group scheme $\calE_{S}\szo$ is $(0:1:0)$.

The open subscheme $\calE\aff$ corresponds to $V\neq0$; in affine 
coordinates $u=U/V$, $w=W/V$, it is given by
\begin{equation}\label{EqAutoTwAff}
	\calE\aff:\quad w=\rho\,P(u,w).
\end{equation}
The fibre  $\calE_{0}\szo$ at $0$ is the affine line 
$W=0$, $V\neq0$ in $\PP^2_{\kappa}$; it is isomorphic to $\Gak{\kappa}$ via the 
map
\begin{equation}\label{EqIdentifGa}
	\begin{array}{rcl}
		\calE_{0}\szo & \fflis & \Gak{\kappa} \\
		(U:V:0) & \longmapsto & U/V
	\end{array}
\end{equation}
or, using the affine coordinates of $\calE\aff$, via the coordinate $u$.
With the above identification, the evaluation at infinity is given on $\calE\aff$ (in 
these affine coordinates $u,v$) by the very simple formula, where $u$ 
is viewed as a function on $\calE\aff$:
\begin{equation}\label{EqEvInftyCoord}
	\begin{array}{rcl}
		\ev_{0}: \calE\aff(\calO) & \ffl & \kappa\\
		z & \longmapsto & u(z)\bmod\mm.
	\end{array}
\end{equation}
The canonical section $\gamma$ is given by 
$(U:V:W)=(1:1:t)$; in particular, its value at 
$0$ is $(1:1:0)$, which is indeed a nonzero element of 
$\calE_{0}\szo$, as predicted by \ref{PropEltCan}\,\ref{PropEltCan4}. 
In fact, using (\ref{EqIdentifGa}) to identify $\calE_{0}\szo(\kappa)$ 
with $\kappa$, we have $\ev_{0}(\gamma)=1$, hence the 
restriction of is $\ev_{0}$ to $\Lambda$ (identified with $\ZZ$) 
is just the natural inclusion of $\ZZ$ into $\kappa$. 

In particular, for the multiplication on $\Lambda$, property 
(\ref{EqCodeMult}) boils down to the following: if $z_{i}$ ($i=1,2,3$) 
are three elements of $\Lambda$, then
\begin{equation}\label{EqCodeMultCoord}
	z_{3}=z_{1}z_{2}\text{ (in $\Lambda$)}\; \Iff\;
	u(z_{3})\equiv u(z_{1})\,u(z_{2})\pmod\mm.
\end{equation}
This has the following consequence:
\begin{prop}\label{PropRedModM} Assume there exists an additive 
subgroup $X$ of $\calO$ with the following properties:
\begin{romlist}
	\item\label{PropRedModM1} $X$ contains all elements of the form 
	$u(z_{1})\,u(z_{2})$ with $z_{1},z_{2}\in\Lambda$;
	\item\label{PropRedModM2} the inclusion of $X_{+}:=X\cap\mm$ into $X$ is 
	relatively \dio\ (as subsets of $\calO'$, with respect to $\calO$).
\end{romlist}
Then the multiplication (hence the whole ring structure) on $\Lambda$ 
is relatively \dio.
\end{prop}
\dem By assumption, there is a \dio\ subset $\mathcal{D}\subset\calO'$ such that 
$\mathcal{D}\cap X=X_{+}$. It follows that if $z_{i}\in\Lambda$ ($i=1,2,3$) we 
have:
$$
z_{3}=z_{1}z_{2} \;\Iff \;u(z_{3})-u(z_{1})\,u(z_{2})\in\mm
\;\Iff\; u(z_{3})-u(z_{1})\,u(z_{2})\in \mathcal{D}
$$
since, by our assumptions,  $u(z_{3})-u(z_{1})\,u(z_{2})\in X$ (note 
that $u(z_{3})=u(z_{3})\,u(\gamma)$).\qed
\medskip

The simplest choice for $X$ is, of course, $X=\calO$, which gives:
\begin{subcor}\label{CorRedModM} Assume that $t$ is not 
invertible in $\calO'$ (in other words, $\mm\calO'\neq\calO'$). Then the 
ring structure on $\Lambda$ is relatively \dio.
\end{subcor}
\dem We have $\mm\calO'\cap\calO=\mm$ since it is a proper ideal of $\calO$ 
containing $\mm$. But of course $\mm\calO'=t\calO'$ is \dio\ in 
$\calO'$, hence we can apply \ref{PropRedModM} with $X=\calO$ (and 
$X_{+}=\mm$).\qed

\begin{subrem}\label{RemCorDioSemiloc} This of course applies to 
$\calO'=\calO$. In fact, at this point we can already conclude that 
\emph{$\calO$ is 
positive-existentially undecidable}; in other words, for any field $k$ 
of characteristic zero, the local ring $k[t]_{(t)}$ is positive-existentially 
undecidable with respect to $\QQ[t]_{(t)}$. Indeed, choose any 
$E$ over $\QQ$ without complex multiplication: then, from 
assertions \ref{PropEltCan3} and \ref{PropEltCan5} of \ref{PropEltCan} 
we have $\Lambda=\calE\aff(\calO)$, so $\Lambda$ is \dio, hence is 
a \dio\ ring by \ref{CorRedModM}. Of course this 
will be generalised later.
\end{subrem}

\Subsection{Description of $\calE$ at infinity. }\label{SsecDescr0}
The results below will be needed 
in Section \ref{Secpadic} to treat the $p$-adic case, because the 
Kim-Roush method involves controlling the order of certain functions 
at $\infty$.  

Denote by $\calR=\kappa[t^{-1}]_{(t^{-1})}$ the local ring of $\Pukp$
at $\infty$, and put $T=\Spec(\calR)$. Assume that $c\neq0$: then $\infty$ is not a 
branch point of $\pi$, and $\calE_{T}$ is an elliptic curve. 
Accordingly, $P(t^{-1},1)$ is a unit of $\calR$, so we can view $\pi$ 
(above $T$) as $\Spec(\calR[\sqrt{\rho'}])$ where 
$\rho'=P(t^{-1},1)^{-1}$ (this will give nicer coordinate changes than 
using $\sqrt{P(t^{-1},1)}$).

We can then describe $\calE_{T}$ by the homogeneous equation
\begin{equation}\label{EqAutoTwZero} \calE_{T}:\quad{V'}^2\,W'=\rho'\,P(U',W')
\end{equation}
in homogeneous coordinates $U',V',W'$. The canonical section is given 
by $(t^{-1}:1:1)$; the corresponding affine model is 
\begin{equation}\label{EqAutoTwZeroAff} (\calE\aff)_{T}:
	\quad w'=\rho'\,P(u',w')
\end{equation}
in affine coordinates $u'=U'/V'$, $w'=W'/V'$.

Of course, over $\Spec\kappa(t)$ (the intersection of $S$ and $T$ in 
$\Pukp$), equation (\ref{EqAutoTwZero}) 
defines the same curve as (\ref{EqAutoTw}); the isomorphism between 
the two models is readily checked to be given by
\begin{equation}\label{EqChCoord}U'=t^{-1}U,\quad V'=V,\quad W'=t^{-1}W.
\end{equation}
In particular, the rational functions $u=U/V$ and $u'=U'/V'$ on 
$\calE$ are related by
\begin{equation}\label{EqChCoordAff}u'=t^{-1}u.
\end{equation}
This implies:
\begin{subprop}\label{ValEnZero} 
	\begin{romlist}
		\item\label{ValEnZero1} Let $z\in\calE(\Pukp)$ be a section. 
		Assume that $z(\infty)\in\calE\aff$. 
		Then the value $u(z)\in\kappa(t)$ of the rational function $u$ at $z$ 
		has order $\geq-1$ at $\infty$.
		\item\label{ValEnZero2} Assume that the condition of \rref{RemChoixOrigine} 
		is satisfied, i.e.~$\gamma(\infty)$ has infinite order in $\calE_{\infty}$. Then 
		for every $z\in\Lambda$, the value $u(z)$ of the function $u$ at 
		$z$ has order $\geq-1$ at $\infty$.
	\end{romlist}
\end{subprop}
\dem \ref{ValEnZero1} The condition implies that $z$ maps 
$T=\Spec(\calR)$ into $\calE\aff$. In particular, $u'(z)$ belongs to $\calR$, 
i.e.~has nonnegative order at $\infty$: the assertion then follows from 
(\ref{EqChCoordAff}).
\smallskip

\noindent\ref{ValEnZero2} The condition in \ref{ValEnZero1} just 
means that $z(\infty)$ is not 
a point of order $2$ on $\calE_{\infty}$. With the assumption of 
\ref{ValEnZero2}, this will be satisfied for $z=n\gamma$ (any 
$n\in\ZZ$), so \ref{ValEnZero2} follows from \ref{ValEnZero1}.\qed

\begin{subrem}\label{RemSansCalcul} Without explicitly computing the 
	coordinate change (\ref{EqChCoord}), it 
	was a priori clear that $u$ must be a polynomial in $u',v'$ with 
	coefficients in $\kappa(t)$; it follows that \ref{ValEnZero} had to hold with 
	$-1$ possibly replaced by some unspecified integer independent of $z$ in 
	\ref{ValEnZero1}
	(resp.~of $n$ in \ref{ValEnZero2}). With some care, this `computation-free' 
	approach would be sufficient for our purposes.
\end{subrem}

\begin{subcor}\label{CorValEnZero} Let $X$ (resp.~$X_{+}$) be the set of 
rational functins in $\kappa(t)$ having order $\geq-2$ at $\infty$ and 
nonnegative (resp.~positive) order at $0$. Assume that $X_{+}$ is a 
relatively \dio\ subset of $X$ (in $\calO'$, with respect to $\calO$).

Then, if $\gamma(\infty)$ has infinite order in $\calE_{\infty}$, the ring 
structure on $\Lambda$ is \dio.
\end{subcor}
\dem Clearly, $X$ is a subgroup of $\calO$ and $X_{+}=X\cap\mm$. Also, 
it follows from \ref{ValEnZero} that $X$ contains all products 
$u(z_{1})u(z_{2})$ for $z_{1},z_{2}$ in $\Lambda$. So this is a 
special case of \ref{PropRedModM}.\qed

\section{\dio\ undecidability of semilocal rings
of curves.}\label{SecIndecSemiloc}
\Subsection{Notations. }\label{NotSemiloc}
\Subsubsection{The function field side. }\label{NotSecIndec}
We denote by $k$ a field of characteristic zero, by $C$ a smooth 
projective geometrically connected curve over $k$, and by $K$ the function 
field of $C$. 

Let $Q$ be a finite nonempty set of closed points of $C$. We 
choose $f$ in $K$ having simple ramification, simple zeros and 
simple poles on $C$, and vanishing at $Q$. 

We denote by $A$ the semilocal ring of $C$ at $Q$:
\begin{equation}\label{DefAnnSemiloc}
	A:=\bigcap_{q\in Q}\calO_{C,q}.
\end{equation}
Thus, $A$ is a regular semilocal domain of dimension $1$ with 
fraction field $K$. It contains $f$, which generates its Jacobson 
radical; since $Q\neq\emptyset$, the intersection $A\cap\QQ(f)$ is the ring
\begin{equation}\label{DefAZero}
	A_{0}:=A\cap\QQ(f)=\QQ[f]_{(f)}.
\end{equation}
All \dio\ sets (in some affine space over $A$) will be relative to 
$A_{0}$.
\Subsubsection{The elliptic curve side. }\label{NotSecIndecEll}
Let us fix an elliptic curve $E$ over $\QQ$. We choose an 
isomorphism $E/\{\pm\mathrm{Id}_{E}\}\flis\Pu{\QQ}$ sending the 
origin to $0$; we denote by $t$ the standard coordinate on 
$\PP^1$.

Applying the constructions of \ref{SsecNotSelfTw} and \ref{SsecSelfTw} 
with $\kappa=\QQ$, we obtain a group scheme
$$ \calE\ffl\Spec(\calO)$$
where $\calO=\QQ[t]_{(t)}$ is the local ring of $\Pu{\QQ}$ 
at $0$ (this $\calE$ was denoted by $\calE_{S}$ or 
$\calE_{\calO}$ in \ref{SsecCan} but we shall not need the 
original $\calE$, which was over $\Pu{\QQ}$). Inside $\calE$ we have 
open subschemes 
$$\calE\aff\subset\calE\szo\subset\calE$$
where $\calE\aff\subset\Aa^2_{\calO}$ is affine. Recall also from 
\ref{SsecCan} that we 
have a canonical section $\gamma\in\calE\aff(\calO)$, generating a 
subgroup $\Lambda=\ZZ\gamma$ of $\calE(\calO)$, which is contained in 
$\calE\aff(\calO)$. We give $\Lambda$ the ring structure 
deduced from the obvious isomorphism $\ZZ\flis\Lambda$.

\Subsubsection{Where both sides meet. }\label{SsecMeet}
For any $\lambda\in\QQ^\ast$, we can send $k[t]$ to $A$ by mapping
$t$ to $(\lambda f)$. This gives a diagram of injective 
ring homomorphisms
$$ \begin{array}{ccccc}
\calO && \calO_{k}\\
\Vert && \Vert\\
\noalign{\vskip 0.5ex}
\QQ[t]_{(t)} & \inj& k[t]_{(t)}& \inj& A\\
	\noalign{\vskip 0.5ex}
	\bigcap & & \bigcap &  & \bigcap \\
	\noalign{\vskip 0.5ex}
\QQ(t)& \inj &k(t) &\inj & K\\
 & & t & \longmapsto &\lambda f.
\end{array}
$$
I claim that 
\begin{equation}\label{EqIntersEaff}
	\calE\aff(\calO_{k})=\calE(k(t))\cap\calE\aff(A).
\end{equation}
This is in fact obvious: if we embed $\calE\aff$ into, say, $\Aa^2$ in the 
usual way, then a point of $\calE\aff(A)$ is in $\calE(k(t))$ if and only 
if its coordinates are in $k(t)$, hence in $k(t)\cap A=\calO_{k}$. 

Now, Theorem \ref{ThIntroTer}\,\ref{ThIntroTer2} (applied with 
$\Gamma=E$) implies that for suitable $\lambda$, we have
\begin{equation}\label{EqEgalite}
	\calE(k(t))=\calE(K)
\end{equation}
(just take $\lambda$ in $\good(k)\cap\QQ$). 

We choose $\lambda$ once and for all with 
this property, and identify $t$ with $\lambda f$, thus viewing all the 
maps in the above diagram as inclusions. Note that, independently of 
$\lambda$, the image of $\calO$ in $A$ is $A_{0}$.
\begin{prop}\label{LambdaDio} With the notations and assumptions of 
\rref{NotSemiloc}, assume 
that $E$ has no complex multiplication over $\CC$. 

Then $\Lambda$ is a \dio\ subset of $A^2$, and of $K^2$ (with respect 
to $A_{0}$).
\end{prop}
\dem By \ref{PropEltCan}\,\ref{PropEltCan3} we have 
$\Lambda=\calE\aff(\calO)=\calE\aff(\calO_{k})$. By (\ref{EqIntersEaff}) 
and (\ref{EqEgalite}) we have 
$\Lambda=\calE(K)\cap\calE\aff(A)=\calE\aff(A)$ hence $\Lambda$ is 
\dio\ in $A^2$. 

Let us show that $\Lambda$ is \dio\ in $K^2$. Recall that 
$\calE(k(t))=\Lambda\times E[2](k)$ by (\ref{EqSecTwSurP1Bis}), 
hence $2\Lambda=2\calE(k(t))=2\calE(K)$. This is also equal to 
$2\calE\aff(K)$ because the complement of $\calE\aff$ in $\calE$ 
consists of the nontrivial $2$-torsion points. But the graph of addition 
is \dio\ in $\calE\aff(K)^3$ (as 
in the proof of \ref{PartieFacile}\,\ref{PartieFacile1}). 
Hence $2\Lambda=2\calE\aff(K)$ is \dio\ in $K^2$, and so is 
$\Lambda=(2\Lambda)\cup(\gamma+2\Lambda)$.\qed
\medskip

We can now prove part (1) of Theorem \ref{ThVague} (more precisely 
the $1$-dimensional case, which implies the general case as explained 
in the introduction):

\begin{thm}\label{ThIndecSemiloc} With the notations and assumptions 
of \rref{NotSecIndec}, there is a \dio\ ring $\Lambda\subset A^2$, 
		isomorphic to $\ZZ$. In particular, 
		the positive-existential theory of $A$ in $\lr(A_{0})$ 
		is undecidable.
\end{thm}
\dem Choose  any elliptic curve $E$ over $\QQ$, without complex 
multiplication over $\CC$. Applying the constructions of \ref{NotSecIndecEll} 
and \ref{SsecMeet}, we conclude from \ref{LambdaDio} that 
$\Lambda$ is \dio\ in $\calE\aff(A)$, hence is a \dio\ 
ring by \ref{CorRedModM}, applied with $\calO'=A$ (the fact that 
$Q\neq\emptyset$ is used here!).\qed

\begin{subrem}\label{ThIndecSemilocExpli}
	The \dio\ ring in \ref{ThIndecSemiloc} has a very simple explicit 
	definition, following from the computations in \ref{EqExplicites}. Let 
	$E\aff$ be given by the equation
	$$z=P(x,z)=x^3+a\,x^2\,z+b\,x\,z^2+c\,z^3$$
	(which is the affine form of (\ref{EqEllSurk})), with $a$, $b$, 
	$c\in\QQ$; of course we assume that the discriminant of $P(x,1)$ is 
	nonzero, and also that $E$ has no complex multiplication over $\CC$ 
	(this is easy to ensure; for instance it is true if  the 
	$j$-invariant of $E$ is not an integer). 
	
	We choose $\lambda\in\QQ^\ast$ as in \ref{SsecMeet}, i.e.~such that 
	(\ref{EqEgalite}) holds, and we define 
	$\calE\aff\subset\Aa^2_{A_{0}}$ by the equation (in affine 
	coordinates $(u,w)$):
	$$ w = \frac{\lambda f}{P(1, \lambda f)} \,P(u,w).$$
	Now, the \dio\ ring $\Lambda$ is defined as follows: 
	\begin{romlist}
		\item the underlying set is 
		the subset $\calE\aff(A)$ of $A^2$ (that is, the 
		set of solutions $(u,w)\in A^2$ of the above equation),
		\item the zero element is $0_{\Lambda}=(0,0)$,
		\item the ring unit is $1_{\Lambda}=\gamma=(1,\lambda f)$,
		\item the addition is given by the elliptic curve law, or equivalently by:
		$$(u'',w'')=(u,w)+(u',w')\; \Iff\;
		\exists\,\alpha, \;u''=u+u'+\alpha f,$$
		\item the multiplication is given by:
		$$(u'',w'')=(u,w)(u',w')\;\Iff\;
		\exists \,\alpha, \;u''=uu'+\alpha f.$$
	\end{romlist}
\end{subrem}

\begin{subrem}\label{RemThIndecSemiloc1}
	Assume that there is a \dio\ subset $\cal D$ of $A$ such 
	that $\ZZ\subset{\cal D}\subset k$ (or, more generally, that 
	$\ZZ\subset{\cal D}$ and the composite ${\cal 
	D}\inj A\to A/fA$ is injective). Then we have the stronger 
	property that $\ZZ$ is \dio\ in $A$: indeed, for 
	$\alpha\in A$, we have $\alpha\in\ZZ$ if and only if $\alpha\in{\cal 
	D}$ and there exists 
	$(u,w)\in\calE\aff(A)$ such that $\alpha-u\in fA$.
\end{subrem}

\begin{subcor}\label{CorThIndecSemiloc1}
	Assume that $A$ (or equivalently, its Jacobson radical $fA$) 
	is \dio\ in $K$. Then there is a \dio\ ring in $K^2$, isomorphic to $\ZZ$.
	
	In particular, the positive-existential theory of $K$ in $\lr(A_{0})$ 
	is undecidable.\qed
\end{subcor}

\begin{subcor}\label{CorCasReelClos}
	Let $k$ be a \emph{real closed} field, and let $C$ be a smooth 
	$k$-curve  having a rational point. Then the function field $K$ of 
	$C$ is positive-existentially undecidable.
\end{subcor}
\dem Let $q\in C(k)$ be a rational point. One can find $\varphi\in K$ 
having only simple zeros, and vanishing at $q$. If $Q$ is the set of 
$k$-rational zeros of $\varphi$, then Theorem (1.8) of \cite{Za}, 
Chapter V shows that the 
semilocal ring of $Q$ is \dio\ in $K$. Hence we can apply 
\ref{CorThIndecSemiloc1}.\qed
\medskip

Apart from this case (which will be superseded by \ref{ThReel}), the 
\dio\ definability of 
$A$ in $K$ seems, in general, rather difficult 
to prove. 

\Subsection{Remarks on effectivity. }\label{SsecRemEff}
As explained in the introduction, the choice of $\lambda$ 
satisfying (\ref{EqEgalite}) is not effective in general. Let us 
describe a procedure for finding such a $\lambda$ if $k$ is finitely 
generated over the prime field. 

More precisely, we assume here that:
\begin{sitemize}
	\item $k$ is presented over $\QQ$ (in the sense of \cite{FJ}, Section 17.2),
	\item $K$ is presented over $k$ (which essentially means that the 
	curve $C$ is explicitly described),
	\item $f$ is explicitly given. 
\end{sitemize}
Assume we have effectively constructed an elliptic 
$E$ over $\QQ$, without complex multiplication over $\CC$, and such 
that $\Hom_{k}\,(E,\Jac(C))=0$. Then by Theorem \ref{MainThEff}, we 
can effectively find $\lambda\in\good(k)\cap\ZZ$: simply list all 
integers until such an element is found, which can be checked 
effectively since it reduces to deciding whether a given polynomial 
in $k[y]$ has no root in $k$. The rest of the proof of \ref{ThIndecSemiloc} 
involves only effective constructions.

Let us now construct $E$ with the required properties. For an 
indeterminate $z$, fix an elliptic curve $E_{\QQ(z)}$ over $\QQ(z)$ with 
$j$-invariant equal to $z$. By ground field 
extension to $k(z)$ we obtain an elliptic curve $E_{k(z)}$ with the 
property that $\Hom_{k(z)}\,(E_{k(z)},\Jac(C)_{k(z)})$ is zero: indeed, by 
\ref{RigSousVarAb}, any abelian subvariety of $\Jac(C)_{k(z)}$ 
is defined over $\kb$, while $E_{k(z)}$ (or any nontrivial quotient 
of it) is not. Hence we can apply Theorem \ref{MainThEff} to find 
$\zeta\in\QQ$ such that $E_{\QQ(z)}$ specialises to an elliptic 
curve $E_{\zeta}$ over $\QQ$ with the property that 
$\Hom_{k}\,(E_{\zeta},\Jac(C))=0$. Moreover we can certainly find such a
$\zeta$ which is not an integer, which implies that the resulting 
$E_{\zeta}$ (whose $j$-invariant is $\zeta$) has no complex 
multiplication.

\section{\dio\ undecidability of real function fields.}\label{SecReel}
The following lemma combines results of \cite{Denef} (for the real case) 
and \cite{KR} (for the $p$-adic case, which will be considered later):

\begin{lem}\label{LemDioDense} Let $k$ be a field of 
	characteristic zero, and let $K$ be a finitely generated 
	regular extension of $k$. 
	\begin{romlist}
		\item\label{LemDioDense1} There is an elliptic curve $E_{1}$ 
		over $\QQ$ with the following properties:
		\begin{subromlist}
			\item\label{cond1} $E_{1}(\QQ)$ is infinite (i.e.~$E$ has positive rank 
			over $\QQ$);
			\item\label{cond2} $E_{1}(K)=E_{1}(k)$.
		\end{subromlist}
		\item\label{LemDioDense2} Let $\Sigma$ be a finite set of 
		independent absolute values on $\QQ$. Denote by
		$\QQ_{\Sigma}=\prod_{v\in\Sigma}\QQ_{v}$ the $\Sigma$-completion of $\QQ$.
		Then there is a $\QQ$-\dio\ subset $\con$ of $K$ such that 
		$\con\subset k$ and $\con\cap\QQ$ is dense in $\QQ_{\Sigma}$.
	\end{romlist}
\end{lem}
\dem \ref{LemDioDense1} We may assume  $K$ transcendental aver 
$k$ (otherwise, $K=k$). By \ref{ProjectImpaire} there is a transcendence basis 
$(z_{1},\ldots,z_{n})$ of $K/k$ such that $K$ is a regular extension of 
$k(z_{1},\ldots,z_{n-1})$. For any elliptic curve $E_{1}$ over $k$, 
we have $E_{1}(k(z_{1},\ldots,z_{n-1}))=E_{1}(k)$ (immediate by induction 
on $n$ since there is no nonconstant rational map from an elliptic 
curve to $\Pu{}$). Thus, to prove \ref{LemDioDense1}, we may replace 
$k$ by $k(z_{1},\ldots,z_{n-1})$ and assume $n=1$. 

Let $C$ be the (projective, smooth, geometrically connected) 
$k$-curve with function field $K$. The elliptic curves with a nonconstant 
morphism from $C$ 
are those appearing (up to isogeny) as factors of the Jacobian of 
$C$, which form a finite set of isogeny classes. So we can choose an 
elliptic curve $E_{0}$ over $\QQ$ which is not $\ol{\QQ}$-isogenous 
to any of these, ad then apply \ref{LemCellRkPos} to find a twist 
$E_{1}$ of $E_{0}$ with positive rank over $\QQ$, thus satisfying 
both conditions. 
\smallskip

\noindent\ref{LemDioDense2} Choose $E_{1}$ as in \ref{LemDioDense1}, 
and write it in affine coordinates as 
$${E_{1,{\rm aff}}}:\qquad w=P(u,w)$$
with $P\in\QQ[u,w]$, homogeneous of degree $3$ and monic in $u$. Let 
$D\subset K$ be the set of $u$-coordinates of points of 
${E_{1,{\rm aff}}}(K)$, and let $\con$ be the set of quotients 
$u_{1}/u_{2}$ with $u_{1}\in D$ and $u_{2}\in D\setminus\{0\}$. Let us 
show that $\con$ has the required properties.

Clearly, $D$ and $\con$ are $\QQ$-\dio\ in $K$, and property 
\ref{cond2} implies that $D\subset k$, hence $\con\subset 
k$ as well. To prove the density property, it suffices to show that 
the closure of $D\cap\QQ$ in $\QQ_{\Sigma}$ contains a neighbourhood 
of $0$. Since $u$ is a local coordinate at the origin of 
${E_{1,{\rm aff}}}$, this will follow if we prove that the closure of 
${E_{1,{\rm aff}}}(\QQ)$ in ${E_{1,{\rm aff}}}(\QQ_{\Sigma})$ 
contains a neighbourhood 
of the origin. This is equivalent to the analogous statement with 
$E_{1}$ instead of ${E_{1,{\rm aff}}}$ since the latter is $E_{1}$ minus 
finitely many nonzero points.

Now, for each $v\in\Sigma$, $E_{1}(\QQ_{v})$ is a compact 
one-dimensional Lie group over $\QQ_{v}$, hence it has an 
open subgroup of finite index $U_{v}$ isomorphic to $\ZZ_{v}$ (if $v$ 
is $p$-adic) or to the circle group $S^1$ (if $v$ is real). Let 
$U\subset E_{1}(\Sigma)$ be the product of the $U_{v}$'s. By property 
\ref{cond1}, $E(\QQ)$ has an element $\gamma$ of infinite 
order; replacing it by some multiple we may assume that $\gamma\in U$. 
For each $v$, the projection of $\gamma$ in $U_{v}\subset 
E_{1}(\QQ_{v})$ still has infinite order, hence generates a subgroup 
whose closure is open. By weak approximation, it easily follows that 
the closed subgroup of $E_{1}(\Sigma)$ generated by $\gamma$ is also 
open. This completes the proof.\qed
\medskip

We can now prove part (2) of Theorem \ref{ThVague}.
\medskip

\begin{thm}\label{ThReel} Let $k$ be a formally real field., and let $K$ be a 
	finitely generated transcendental extension of $k$, which is also formally 
	real. Then there is a \dio\ ring in $K^2$, isomorphic to $\ZZ$. In 
	particular, the \dio\ theory of $K$ is undecidable.
\end{thm}
\dem By replacing $k$ by a bigger subfield of $K$, we may assume 
that $K$ has transcendence degree one over $k$, and that $k$ is 
algebraically closed in $K$. Then $K$ is the function field of a 
projective, smooth, geometrically connected $k$-curve $C$. Moreover, 
the assumption that $K$ is formally real means that $C$ has a 
closed point $q$ with formally real residue field. Putting $Q=\{q\}$, we 
are in the 
situation of \ref{NotSecIndec}. Choosing any elliptic curve $E$ over 
$\QQ$, without complex multiplication, we can perform the 
constructions of \ref{NotSecIndecEll} and \ref{SsecMeet}, and we keep 
the same notations. In particular we have a
ring $\Lambda\cong\ZZ$ in $A^2$, where $A$ is the local ring of $q$. Moreover, 
we know from \ref{LambdaDio} that $\Lambda$ is \dio\ in $K^2$, and it 
remains only to prove that the multiplication in $\Lambda$ is 
relatively \dio. 

To do this, we apply Proposition \ref{PropRedModM} with $\kappa=\QQ$ 
(hence $\calO=A_{0}=\QQ[f]_{(f)}$), $X=A_{0}$, and $\calO'=K$. 
Thus, all we have to prove is that the maximal ideal $\mm_{0}$ of $A_{0}$ is a 
relatively \dio\ subset of $A_{0}$ in $K$.

First take $\con\subset K$ as provided by \ref{LemDioDense}, applied 
with our $k$ and $K$, and with $\Sigma$ consisting of the ordinary absolute 
value (thus, $\QQ_{\Sigma}=\RR$). Consider the following formula in one variable $x$:
$$\varphi(x):\quad  \exists\,\alpha,\beta,x_{1},\ldots,x_{5}: \;
\alpha\in\con\wedge\beta\in\con\wedge(\alpha-f^{-1})\,x^2+\beta= 
x_{1}^2+\ldots+x_{5}^2.$$
I claim that the set $\mathcal{D}\subset K$ defined by $\varphi$ 
satisfies $\mathcal{D}\cap A_{0}=\mm_{0}=f\,A_{0}$.

First, let us show that $\mathcal{D}\subset f\,A$ (which implies 
that $\mathcal{D}\cap\QQ(f)\subset f\,A_{0}$). Indeed, for some 
$x\in K$, assume that $\varphi(x)$ holds and $x$ 
does not vanish at $q$. The elements $\alpha$ and $\beta$ in 
$\varphi(x)$ must be in $k$, hence $(\alpha-f^{-1})\,x^2+\beta$ has negative 
odd order at $q$. Since the residue field of $q$ is real, this cannot 
be a sum of squares in $K$, which contradicts $\varphi(x)$.

Let us now prove that $f\,A_{0}\subset\mathcal{D}$. We view $A_{0}$ 
as the local ring of $0$ in $\Pu{\QQ}$, with standard coordinate $f$. 
If $x\in f\,A_{0}$, then $x$ is a rational function on $\Pu{\QQ}$, 
vanishing at $0$; hence, so does $f^{-1}x^2$. In particular, 
$\vert f^{-1}x^2\vert\leq1$ on $I:=[-\varepsilon,\varepsilon]$ for 
some $\varepsilon>0$. We can choose $\alpha$ and $\beta$ in $\con$, and such 
that $\beta>1$ and $\alpha>1/\varepsilon$. Then, by our choices 
(recall that $f$ is the standard coordinate on $\Pu{\QQ}$, hence 
$\vert f\vert>\varepsilon$ on $\RR\setminus I$):
\begin{sitemize}
	\item on $\RR\setminus I$, we have $\vert 
	f^{-1}\vert<1/\varepsilon<\alpha$, hence 
	$(\alpha-f^{-1})\,x^2+\beta\geq\beta>0$;
	\item on $I$, we have $(\alpha-f^{-1})x^2+\beta=
	\alpha\,x^2-f^{-1}x^2+\beta\geq\alpha\,x^2-1+\beta>0$.
\end{sitemize}
This implies that $(\alpha-f^{-1})\,x^2+\beta$ is a sum of squares in 
$\QQ(f)$, and in fact a sum of five squares 
by \cite{Pour}. Hence $\varphi(x)$ is satisfied, and the proof is 
complete.\qed

\section{\dio\ undecidability of $p$-adic function fields.}\label{Secpadic}
In this section we prove part (3) of Theorem \ref{ThVague}:
\begin{thm}\label{Thpadic} Let $p$ be an odd prime. Let $\kappa$ be a subfield 
of a finite extension of $\QQ_{p}$, and let $K$ be a finitely 
generated transcendental extension of $k$.  

Then there is a \dio\ ring $\Lambda\subset K^d$, for some $d$, 
isomorphic to $\ZZ$ as a ring. In particular, $K$ is 
positive-existentially undecidable.
\end{thm}

\Subsection{Extending the ground field. }\label{ExtCorpsBase}
First, we may replace $\kappa$ by its algebraic closure in $K$ and 
assume that $K/k$ is a regular extension.

Next, assume there is a finite extension $L$ of $K$, of degree $n$, and a 
\dio\ ring isomorphic to $\ZZ$ in $L^2$. Then by using the Weil 
restriction (i.e.~by fixing a $K$-basis of $L$ and identifying $L^d$ 
with $K^{2n}$) we obtain a 
\dio\ ring isomorphic to $\ZZ$ in $K^{2n}$. In particular, to prove 
\ref{Thpadic} we may replace $\kappa$ by a finite extension (which we shall 
always view as embedded in some finite extension of $\QQ_{p}$, and in 
particular equipped with a $p$-adic valuation, normalised in such a 
way that its value group is $\ZZ$). 

Thus, replacing if necessary $\kappa$ by a finite extension $\kappa'$ 
and $K$ by $K\otimes_{\kappa}\kappa'$, we shall assume from now on that:
\begin{romlist}
	\item\label{ExtCorpsBase1} there is a transcendence basis 
	$(z_{1},\ldots,z_{n})$ of $K$ over $\kappa$ such that 
	$[K:\kappa(z_{1},\ldots,z_{n})]$ is odd and the extension 
	$K/\kappa(z_{1},\ldots,z_{n-1})$ is regular,
	\item\label{ExtCorpsBase2} $\kappa$ contains elements $i$, $a$, $\varpi$ 
	such that:
	\begin{subromlist}
		\item\label{ExtCorpsBase21} $i^2=-1$,
		\item\label{ExtCorpsBase22} $a$ is a root of unity,
		\item\label{ExtCorpsBase23} $\varpi$ is algebraic over $\QQ$, and 
		has odd $p$-adic valuation,
		\item\label{ExtCorpsBase24} the $4$-dimensional quadratic form
		\begin{equation}\label{EqFormeQuad}
			\langle1,a\rangle\langle1,\varpi\rangle=x^2+\varpi\,y^2+a\,z^2+ 
			a\varpi\,w^2
		\end{equation}
		is anisotropic over $\kappa$,
		\item\label{ExtCorpsBase25} the quadratic form (\ref{EqFormeQuad}) 
		is isotropic at all $2$-adic primes of the field $\QQ(i,a,\varpi)$.
	\end{subromlist}
	\pauseromlist
Indeed, \ref{ExtCorpsBase1} holds over some finite extension of 
$\kappa$ by \ref{ProjectImpaire}. The 
fact that $\kappa$  can be further enlarged to satisfy \ref{ExtCorpsBase2} 
is proved in \cite{KR}, Proposition 8; here we denote by $\varpi$ 
what was (somewhat confusingly) called $p$ in \cite{KR}.

(In the left-hand side of (\ref{EqFormeQuad}) we use the standard 
notation $\langle d_{1},\ldots,d_{n}\rangle$ for the 
diagonal quadratic form $\sum_{j=1}^n d_{j}\,x_{j}^2$, and the product is the 
`Kronecker product', or tensor product.)
\smallskip

From now on we fix $z_{1},\ldots,z_{n-1}$ as in \ref{ExtCorpsBase1}, $i$, 
$a$, $\varpi$ as in \ref{ExtCorpsBase2}, and we put 
$$k:=\kappa(z_{1},\ldots,z_{n-1}).$$
Thus, $K$ is a one-variable function field over $k$; we denote by $C$ 
the smooth, projective, geometrically connected $k$-curve with 
function field $K$. By condition \ref{ExtCorpsBase1}, $C$ is a 
cover of $\Pu{k}$ of odd degree, hence admits a divisor of odd 
degree. By \ref{ExistAdm}, this implies:
	\finpauseromlist
	\item\label{ExtCorpsBase3}
	there is an element $f$ of $K$ which, 
	viewed as a $k$-morphism $C\to\Pu{k}$, has simple ramification, 
	simple zeros, simple poles, and odd degree. 
\end{romlist}
From now on we fix $f$ as in \ref{ExtCorpsBase3}; in fact we may replace 
$z_{n}$ by $f$ and consider the tower of extensions
\begin{equation}\label{TourExt}
	\kappa\subset k\subset k(f)=\kappa(z_{1},\ldots,z_{n-1},f)\subset K
\end{equation}
in which the first two inclusions are purely transcendental and the 
last is finite and odd.
\smallskip

We have done all this to use our results on curves while ensuring the 
following property:
\begin{sublem}\label{LemFQAnis} Every anisotropic quadratic form 
over $\kappa(f)$ remains anisotropic over $K$.
\end{sublem}
\dem Let $\varphi$ be such a quadratic form. Clearly, $\varphi$ is 
still anisotropic over $k(f)$ which is purely transcendental over 
$\kappa(f)$ (\cite{lam}, Chapter 9, Lemma 
1.1). Since $K$ is finite of odd degree over $k(f)$, we 
conclude from Springer's theorem (\cite{lam}, Chapter 7, Theorem 2.3) 
that $\varphi$ is also  anisotropic over $K$.\qed

\Subsection{Defining the ring $\Lambda$. }\label{SsecDefLambdaPadic}
With $f:C\to\Pu{k}$ as in \ref{ExtCorpsBase}\,\ref{ExtCorpsBase3}, 
we choose a zero $q$ of $f$ on $C$ (not necessarily $k$-rational), 
put $Q=\{q\}$, and we adopt the notations of \ref{NotSecIndec}. 

We choose an elliptic curve $E$ over $\QQ$, without complex 
multiplication, and we identifiy $L:=E/\{\pm\mathrm{Id}_{E}\}$ with 
$\Pu{\QQ}$ in such a way that the origin of $E$ goes to $0$, 
and the condition of \ref{RemChoixOrigine} is satisfied; in other 
words, we choose the equation (\ref{EqEllSurk}) in such a way that 
$c\neq0$ and the points $(0:\pm\sqrt{c}:1)$ of $E(\ol{\QQ})$ have 
infinite order.

We then proceed with the constructions of \ref{NotSecIndecEll} and 
\ref{SsecMeet}. We obtain a subset $\Lambda\subset A^2$ with a 
ring structure isomorphic to $\ZZ$, which is a \dio\ subset of $K^2$ 
by \ref{LambdaDio}. To prove that the multiplication of $\Lambda$ is 
\dio, we need the following refinement of \ref{CorValEnZero} 
(cf.~\cite{KR}, Theorem 6, where $t$ 
corresponds to our $f$):

\begin{lem}\label{VarKR} Denote by 
$v_{\infty}$ (resp.~$v_{0}$) the valuation on $\QQ(f)\subset K$ such 
that $v_{\infty}(f)=-1$ (resp.~$v_{0}(f)=+1$). Define subsets $Y_{0}$, $Y_{1}$, 
$Y$ of $\QQ(f)$ by
\begin{equation}\label{DefX0X1}
	\begin{array}{rcl}
		Y_{i}& := & \{r\in\QQ(f)\,\mid\,v_{\infty}(r)=-2\text{ and } 
		v_{0}(r)=i\}\quad(i=0,\:1)\\
		Y_{\phantom{i}}& :=& Y_{0}\cup Y_{1}.
	\end{array}
\end{equation}
Assume that $Y_{1}$ is a relatively \dio\ subset of $Y$ (in $K$). 
Then the ring structure of $\Lambda$ is \dio. Hence $K$ is 
positive-existentially undecidable.
\end{lem}
\dem By assumption, there is a \dio\ set $\mathcal{D}\subset K$ such that 
$\mathcal{D}\cap Y=Y_{1}$. Put
$$\begin{array}{rcl}
X & := & \{r\in\QQ(f)\,\mid\,v_{\infty}(r)\geq-2\text{ and 
}v_{0}(r)\geq0\}\\
X_{+} & := & \{r\in\QQ(f)\,\mid\,v_{\infty}(r)\geq-2\text{ and 
}v_{0}(r)>0\}.
\end{array}
$$
We shall prove that 
$X_{+}$ is relatively \dio\ in $X$, which by \ref{CorValEnZero} 
implies the result.

If $r\in X$, then $\frac{1}{1+f^2}\,r$ has nonnegative $v_{\infty}$ and 
the same $v_{0}$ as $r$. It follows that if we put 
$$s:=f+{f^2}+\Bigl(\frac{1}{1+f^2}\,r\Bigr)^2,$$
then we have
$$\begin{array}{lcl}
v_{\infty}(s) & = & -2,\\
v_{0}(s) & = & 
\begin{cases}
	0 & \text{ if }v_{0}(r)=0\\
	1 & \text{ if } v_{0}(r)>0.
\end{cases}
\end{array}
$$
Hence, for any $r\in X$, we have $s\in Y$, and $s\in Y_{1}$ if and 
only if $r\in X_{+}$. Consequently, $X_{+}=X\cap \mathcal{D}_{1}$, where
$$\mathcal{D}_{1}=\Bigl\{\:r\in 
K\,\mid\,f+{f^2}+\Bigl(\frac{1}{1+f^2}\,r\Bigr)^2\in \mathcal{D}\:\Bigr\},$$
which proves the lemma.\qed

\Subsection{Isotropy of quadratic forms. }\label{IsotrFQ}
It remains to prove that the assumption of Lemma \ref{VarKR} 
is satisfied, i.e.~$Y_{1}$ is relatively \dio\ in $Y$; we follow 
\cite{KR}, indicating only the changes to be made. 

To stick to the notations of \cite{KR}, we put $t=f$ from 
now on. (Thus we forget the $t$ of \ref{NotSecIndecEll}, which 
corresponds to $\lambda\,f$ for some $\lambda\in\QQ$). 

Applying Lemma \ref{LemDioDense}, we fix a $\QQ$-\dio\ 
subset $\con$ of $K$, contained in $\kappa$ and such that $\QQ\cap\con$ is 
dense in $\QQ_{p}$.

To every $r\in K$ we associate two elements $u_{0}$, $u_{1}$ of $K$ 
and two quadratic forms $\varphi_{0}$, $\varphi_{1}$ over $K$ 
depending on parameters $c_{3}$, $c_{5}$, by the formulas
\begin{equation}\label{EqDefFQParam}
	\begin{array}{rcll}
		u_{e} & := & a^e\,((1+t)^3\,r+c_{3}\,t^{3}+c_{5}\,t^{5}) & 
		(e=0,1)\\
		\varphi_{e} & := & \langle 
		t,at,-1,-u_{e}\rangle\langle1,\varpi\rangle.
	\end{array}
\end{equation}
We define a \dio\ set $\mathcal{D}\subset K$ by
$$r\in \mathcal{D}\;\Iff\;\exists\,c_{3}, c_{5}\in\con \text{ such that 
}\varphi_{0}\text{ and }\varphi_{1}\text{ are isotropic over }K.$$
and claim that $\mathcal{D}\cap Y=Y_{1}$. This amounts to proving that 
$Y_{1}\subset\mathcal{D}$ and $Y_{0}\cap\mathcal{D}=\emptyset$.

\Subsubsection{The relation $Y_{0}\cap\mathcal{D}=\emptyset$. }
Assume that $r\in Y_{0}$. Then $r$ is in $\QQ(t)$ and has order 
$0$ at $0$, so the same holds for $u_{0}$, for any choice of 
$c_{3}$ and $c_{5}$ in $k$ (and in particular in $\con$). 
By the first assertion of  \cite{KR}, Proposition 7
(applied with $b=\varpi$, $g=\text{ our }u_{0}$, and $a=\text{ our 
}-a$), this implies that one of the forms $\varphi_{0}$, $\varphi_{1}$ is 
anisotropic over $\kappa(t)$, hence also over $K$ by Lemma \ref{LemFQAnis}.

\Subsubsection{The inclusion $Y_{1}\subset \mathcal{D}$. }
Assume that $r\in Y_{1}$. We have to show that for some choice of 
$c_{3}$ and $c_{5}$ in $\con$ (in fact, we can take them in 
$\con\cap\QQ$) both forms $\varphi_{0}$, $\varphi_{1}$ 
are isotropic over $K$ (and in fact, over $\kappa(t)$). We refer to 
\cite{KR} for the details: first, it is shown in the proof of \cite{KR}, Theorem 
9 that for suitable $c_{3}$ and $c_{5}$ in $\con\cap\QQ$, some 
condition on the Newton polygons of $u_{0}$ and $u_{1}$ is satisfied 
(the only thing that matters about $\con\cap\QQ$ is $p$-adic 
density). Then, the results of \cite{KR}, Section $3$ (in particular 
Theorem 21) imply that this 
Newton polygon condition in turn implies isotropy. This completes the 
proof.\qed


\begin{thebibliography}{SGA~7 }
%
	\bibitem[A]{And} \textsc{G.~W.~Anderson}, \textsl{Abeliants and their 
	application to an elementary construction of Jacobians}, preprint, 
	Univ.~of Minnesota, 2002.
	
%
	\bibitem[B-L-R]{BLR} \textsc{S.~Bosch},  \textsc{W.~L\"utkebohmert}, 
	and 
	\textsc{M.~Raynaud}, \textsl{N\'{e}ron Models}, Ergeb.\ Math.\ 
	Grenzgeb.~(3) Band 21, Springer (Berlin), 1990.
	
%
	\bibitem[D1]{Denef1} \textsc{J.~Denef}, 
	\textsl{Diophantine sets over $\ZZ[T]$\/}, 
	Proc.\ Amer.\ Math.\ Soc.\ 69 (1978), 148--150.
	
%
	\bibitem[D2]{Denef} \textsc{J.~Denef}, 
	\textsl{The Diophantine Problem for Polynomial Rings and Fields of 
	Rational Functions\/}, Trans.\ Amer.\ Math. Soc.\ 242 (1978), 391--399.
	
%
%
	\bibitem[E1]{Eis1} \textsc{K.~Eisentr\"ager}, \textsl{Hilbert's tenth 
	problem for function fields of varieties over $\CC$}, 
	Int.\ Math.\ Res.\ Notes, 59 (2004), 3191--3205.
%
	\bibitem[E2]{Eis} \textsc{K.~Eisentr\"ager}, \textsl{Hilbert's tenth 
	problem for function fields of varieties over number fields and 
	$p$-adic fields} (preliminary version, August 30, 2004).		
	%
%
	\bibitem[EGA~4]{ega4} \textsc{A.~Grothendieck} and 
	\textsc{J.~Dieudonn\'{e}}, \textsl{\'El\'{e}ments de 
	g\'{e}om\'{e}trie alg\'{e}brique, IV: \'Etude locale des 
	sch\'{e}mas et des morphismes de 
	sch\'{e}mas (quatri\`{e}me partie)}, Pub.~Math. I.H.\'E.S.~32 (1967).
	
%
	\bibitem[F-J]{FJ} \textsc{M.D.~Fried} and \textsc{M.~Jarden}, 
	\textsl{Field Arithmetic}, Ergeb.~Math.~Grenzgeb.~11, 
	Springer (1986).
%
	\bibitem[F-W]{FW} \textsc{G.~Faltings}, \textsc{G.~W\"ustholz}, et 
	al., \textsl{Rational Points}, Vieweg (1984). 
%
	\bibitem[K-R1]{KR1} \textsc{K.~H.~Kim} and \textsc{F.~W.~Roush}, 
	\textsl{Diophantine undecidability of $\CC(t_{1},t_{2})$}, 
	J.~of Algebra, 150 (1992), 35--44. 
%
	\bibitem[K-R2]{KR} \textsc{K.~H.~Kim} and \textsc{F.~W.~Roush}, 
	\textsl{Diophantine Unsolvability over $p$-Adic Function Fields}, 
	J.~of Algebra 176 (1995), 83--110.	
%
	\bibitem[Lam]{lam} \textsc{T.Y.~Lam}, 
	\textsl{The Algebraic Theory of Quadratic Forms}, Benjamin (1973).
	
%
	\bibitem[Lan]{lang} \textsc{S.~Lang}, 
	\textsl{Fundamentals of Diophantine Geometry}, Springer (1983).
	
	\bibitem[Ma]{mazur} \textsc{B.~Mazur}, 
	\textsl{Questions of Decidability and Undecidability in Number Theory}, 
	J.~of Symbolic Logic 59 (1994), 353--371.
	
%
	\bibitem[Mu]{Mum} \textsc{D.~Mumford}, \textsl{Abelian Varieties}, 
	Oxford University Press (1974).
	
%
	\bibitem[Mu-F]{GIT} \textsc{D.~Mumford} and \textsc{J.~Fogarty}, 
	\textsl{Geometric Invariant Theory}, 2nd enlarged edition, Springer 
	(1982).
	
%
	\bibitem[N]{Noot} \textsc{R.~Noot}, \textsl{Abelian 
	varieties---Galois representations and properties of ordinary 
	reduction}, Compositio Math.~97 (1995), 161--171.
	
%
 	\bibitem[P-Z]{PhZa} \textsc{T.~Pheidas} and \textsc{K.~Zahidi}, 
 	\textsl{Undecidability of existential theories of rings and fields: 
 	A survey}, in \textsl{Hilbert's Tenth Problem: Relations with 
 	Arithmetic and Algebraic Geometry}, Contemp.~Math.~270 (2000), 
 	49--105.
	
%
	\bibitem[Po]{Pour} \textsc{Y.~Pourchet}, 
	\textsl{Sur la repr\'{e}sentation en somme de carr\'{e}s des 
	polyn\^{o}mes \`{a} une 
	ind\'{e}termin\'{e}e sur un corps de nombres alg\'{e}briques}, 
	Acta Arith.~XIX (1971), 89--104.
	
%
	\bibitem[Sa]{Saito} \textsc{T.~Saito}, 
	\textsl{Vanishing Cycles and Geometry of Curves over a Discrete 
	Valuation Ring}, Amer.~J.~Math.~109 (1987), 1043--1085.
	
%
	\bibitem[Se1]{SerreRibet} \textsc{J.-P.~Serre}, 
	\textsl{Lettre \`a\ Ken Ribet du 1/1/1981}, \OE uvres (Collected 
	Papers), Volume IV, Springer (2000).
	
%
	\bibitem[Se2]{Serre} \textsc{J.-P.~Serre}, 
	\textsl{Lectures on the Mordell-Weil Theorem}, Vieweg (1997).
	
%
	\bibitem[SGA~7]{sga} \textsc{A.~Grothendieck} et al., 
	\textsl{Groupes de Monodromie en G\'{e}om\'{e}trie Alg\'{e}brique 
	(SGA~7~I)}, Lecture Notes in Math. 288, Springer (1972).
	
%
	\bibitem[Z]{Za} \textsc{K.~Zahidi}, \textsl{Existential 
	undecidability for rings of algebraic functions}, thesis, University 
	of Ghent (1999).
	
%
\end{thebibliography}
\end{document}